
\input amstex \documentstyle{amsppt}\catcode `\@=11 \def\logo@{}\catcode
`\@=12\NoBlackBoxes
\pagewidth{6.0in} \pageheight{8.5in} \hcorrection{0.5in}
\vcorrection{0in} \baselineskip=17pt \parindent=20pt \parskip=0pt
  \define\pd#1#2{\frac{\partial #1}{\partial #2}}
\newsymbol\zR 2052\newsymbol\zZ 205A\newsymbol\zQ 2051\newsymbol\zN
204E
\NoRunningHeads

\topmatter
\title Poisson-Lie Structures on Infinite-Dimensional Jet Groups
and Quantum Groups Related to Them\endtitle
\author Ognyan Stoyanov$^{*}$\endauthor
\thanks  $^{*}$ Department of Mathematics, Rutgers University, New Brunswick,
NJ 08903, USA,\newline e-mail: stoyanov\@math.rutgers.edu
 \endthanks
\email stoyanov@math.rutgers.edu\endemail
\abstract
We study the problem of
classifying all Poisson-Lie structures on the group $G_{\infty}$ of formal
diffeomorphisms of the real line $\zR^{1}$ which leave the
origin fixed, as well as the extended group of diffeomorphisms
$G_{0\infty}\supset G_{\infty}$ whose action on $\zR^{1}$ does not
necessarily fix the origin.
A complete local classification of all Poisson-Lie structures on the groups
$G_{\infty}$ and $G_{0\infty}$ is given.  This includes a classification of all
Lie-bialgebra structures on the Lie algebra $\Cal G_{\infty}$ of $G_{\infty}$,
which we prove to be all of coboundary type, and a classification of all
Lie-bialgebra strucutures on the Lie algebra
$\Cal G_{0\infty}$ (the Witt algebra) of $G_{0\infty}$ which also turned
out to be all of coboundary type.
A large class of Poisson structures on the space $V_{\lambda}$ of
$\lambda$-densities on the real line
is found such that $V_{\lambda}$ becomes a homogeneous Poisson space under the
action of the Poisson-Lie group $G_{\infty}$.
We construct a series of quantum semigroups whose quasiclassical
limits are finite-dimensional Poisson-Lie factor groups of $G_{\infty}$ and
$G_{0\infty}$.
\endabstract
\date June 8, 1995\enddate
\endtopmatter

\head Introduction \endhead

Quantum groups have been introduced in [D2,Ji] as deformations
of universal enveloping algebras of Lie groups and of the algebra of
functions on Poisson-Lie groups [FRT]. The latter are Lie groups equipped with
Poisson structures compatible with the group structure
(from where the term Poisson-Lie group originates). In this approach to
constructing quantum groups the first step is to analyze existence
of Poisson-Lie structures on the corresponding Lie group. The question of
classifying all Poisson-Lie structures on a given Lie group
(provided any exist) has been posed originally by Drinfel'd and
Belavin \cite{BD}. In the same paper they give a complete solution for
the
case of finite dimensional complex (semi)simple Lie groups. The
problem, in general, is very difficult. It has been solved for some
other
groups in low dimensions. Let us give a list of groups for which the solution
of the
classification problem is known at present:\newline
(a) Finite dimensional complex (semi)simple Lie groups [BD],\newline
(b) The groups $GL(2,\zR)$, $SL(2,\zR)$, $GL(1|1)$ [Ku1,2],\newline
(c) The 3-dimensional Heisenberg group [Ku3], and some higher
unipotent groups \cite{KS2}\newline
(d) The group of affine transformations of the line $Aff(1)$,\newline
(e) The group of motions of $\zR^1\times\zR^1$,\newline
(f) The Lorentz group considered as a realification of $SL(2,{\Bbb C})$
[Za].\newline
 Note that all the groups mentioned above are {\it
finite-dimensional}, and the only infinite {\it series} for which the
classification has been completed are {\it complex} not {\it real} groups.
\par In the work presented here we study the problem of local
classification (up to a local change of coordinates) of all the Poisson-Lie
structures on the group $G_{\infty}$ of formal diffeomorphisms of the real line
$\zR^{1}$ which leave the
origin fixed, as well as the extended group of diffeomorphisms
$G_{0\infty}\supset G_{\infty}$ whose action on $\zR^{1}$ does not
necessarily fix the origin. The latter we treat as a formal group.
\par The existence of Poisson-Lie structures on $G_{\infty}$ and
$G_{0\infty}$ is far from being obvious. For instance, since
$G_{\infty}$ is a projective limit of groups of finite jets (cf. Sec. 2), if we
consider the group of 3-jets leaving the origin fixed, then there
exists a Poisson-Lie structure on this group [KS1], which can not be extended
to $G_{\infty}$. Even though the above groups are infinite
dimensional, surprisingly, the classification problem has a complete solution.
\par The Lie algebras of the groups $G_{0\infty}$ and $G_{\infty}$ are
the Witt algebra
$$\Cal G_{0\infty}={\text{ span}_k}\{e_i\mid [e_i,e_j]=(i-j)e_{i+j},\ \
i,j\ge -1\},\ \ {\text{where}}\ \ k={\Bbb R}\text{ or }{\Bbb C},$$
and its principal
subalgebra
$\Cal G_{\infty}\subset\Cal G_{0\infty}$. We prove that there is a one-to-one
correspondence between the Poisson-Lie
structures on $G_{\infty}$ and the Lie-bialgebra structures on $\Cal
G_{\infty}$. The latter are shown to be {\it all} of coboundary
type. {\it All} Lie-bialgebra structures on $\Cal G_{0\infty}$ are of
coboundary type, they are all classified, and there is a one-to-one
correspondence
between them and an explicitly listed family of the Poisson-Lie structures on
$G_{0\infty}$. Thus a complete classification of all
Lie-bialgebra structures on the Witt algebra $\Cal G_{0\infty}$
and its principal subalgebra $\Cal G_{\infty}$ is given.
\par With an arbitrary fixed Poisson-Lie structure the group $G_{\infty}$ acts
naturally on the space $V_{\lambda}$ of $\lambda$-densities on $\zR^{1}$. For
each Poisson-Lie
structure on the group $G_{\infty}$ we determine a Poisson structure on
$V_{\lambda}$ such that $V_{\lambda}$ becomes a homogeneous
Poisson $G_{\infty}$-space under the action of the Poisson-Lie group
$G_{\infty}$.
\par Finally, the quantization problem is addressed. We construct a
series of finitely generated non-commutative non-cocommutative
bialgebras (quantum semigroups) whose quasi-classical
limits are finite-dimensional Poisson-Lie factor groups of $G_{\infty}$ and
$G_{0\infty}$. The Poisson-Lie structures on these
finite-dimensional groups are restrictions of the Poisson-Lie
structures on $G_{\infty}$ and $G_{0\infty}$.
\par We give now a brief guide to the organization of the text.
 In Section 1 we introduce the basic concepts related to the Poisson-Lie theory
and formulate the fundamental theorem of Drinfel'd
relating Poisson-Lie groups and Lie-bialgebras.
 In Section 2 we introduce the infinite-dimensional group $G_{\infty}$ and a
smooth structure on it.
In Section 3 we find all bialgebra structures on the Lie algebras
$\Cal G_{0\infty}$ and  $\Cal G_{\infty}$.
In Section 4 we find a class of Poisson-Lie structures on $G_{\infty}$.
In Section 5 we show that there is a one-to-one correspondence between the
Lie-bialgebra structures on $\Cal G_{\infty}$ and the
Poisson-Lie structures found on $G_{\infty}$, and prove that the latter give a
complete list of all Poisson-Lie structures on
$G_{\infty}$.
In Section 6 we study Poisson-Lie structures on $G_{0\infty}$ which correspond
to all Lie-bialgebra structures on $\Cal G_{0\infty}$.
Section 7 is devoted to elements of representation theory for the Poisson-Lie
group $G_{\infty}$ on the homogeneous spaces
$V_{\lambda}$.
In Section 8 we describe a series of finitely generated quantum
semigroups. These, we believe to be precursors of quantizations of $G_{\infty}$
and $G_{0\infty}$ which are presently unknown.
\head 1. Poisson-Lie Theory\endhead
\par In this section we review the basic objects to be studied:
Poisson manifolds, Poisson-Lie groups, Lie bialgebras, and
some basic results about them.
\par Let $\Cal M$ be a finite-dimensional smooth manifold. A Poisson structure
(bracket) on
$\Cal M$ is defined as a bilinear map $\{\  ,\  \}\:{C}^{\infty}(\Cal M)\times
{C}^{\infty}(\Cal M)\to {C}^{\infty}(\Cal M)$, which
makes ${C}^{\infty}(\Cal M)$ into a Lie  algebra, and is a derivation with
respect to each argument. That is,
there exists a section $\omega \in {\wedge }^2 T_{\Cal M}$, where $T_{\Cal M}$
is the tangent bundle of $\Cal M$, such that for any $f,g,h\in
{C}^{\infty}(\Cal M)$ we have
$( f,g) \mapsto \{ f,g\} =\Bigl< {\omega},df\wedge dg\Bigr>,$
and\newline
 (i)   $\{\{ f,g\} ,h\} +\{\{ g,h\} ,f\} +\{\{ h,f\} ,g\} =0$ (Jacobi
identity);\newline
 (ii)  $\{ f,gh\} =\{ f,g\} h+\{ f,h\} g$ (derivation property);\newline
 (iii) $\{ f,g\}=-\{ g,f\}$ (antisymmetry),\newline
 where $\langle\ ,\ \rangle$ denotes the natural pairing between
${\wedge }^2 T_{\Cal M}$ and ${\wedge }^2 T^{*}_{\Cal M}$, where
$T^{*}_{\Cal M}$ is the tangent bundle of $\Cal M$. The second
property, (ii), amounts to compatibility between the Lie algebra
structure defined by $\{\ ,\ \}$ and the multiplication in ${C}^{\infty}(\Cal
M)$.
 In local coordinates,
$$\{ f,g\}(x) =\omega _{ij}(x)\pd {f}{x_i} \pd {g}{x_j},$$
where ${\omega}_x=\omega _{ij}(x)\pd {}{x_i}\wedge\pd {}{x_j}\in \wedge^2T_x$
is a bi-vector field at the point $x\in \Cal M$, and $\bigl\{\frac
{\partial}{\partial x_i}\bigr\}$ is a basis of the tangent space $T_x$ at $x\in
\Cal M$ in the local coordinates $(x_i)$.
\par Here and throughout this text a summation is understood over
repeated nonfixed indices unless stated otherwise. Note also that our
convention about the position of indices of tensors is the opposite to
the standard one. Namely, all contravariant(covariant) tensors have
lower(upper) indices. We found this notation more convenient when
working with power series, and hope that it will not create confusion.
\par The Jacobi identity (i) is equivalent to the following system of equations
for the components $\omega_{ij}(x)=-\omega_{ji}(x)$:
$$\omega_{ij}\frac {\partial\omega_{kl}}{\partial x_i}+\omega_{ik}\frac
{\partial\omega_{lj}}{\partial x_i}+\omega_{il}\frac
{\partial\omega_{jk}}{\partial x_i}=0.\tag 1.1$$
\definition{Definition 1.1 {\rm (Poisson Manifold)[L]}} A Poisson manifold is a
smooth manifold with a Poisson structure.\enddefinition
\par A smooth map $F\: {\Cal M}_1\to {\Cal M}_2$, of two Poisson manifolds
${\Cal M}_1$ and ${\Cal M}_2$, is said to
be {\it Poisson} if
$F^*(\{ g,h\}_{{\Cal M}_2}) =\{ {F^*}(g),{F^*}(h)\}_{{\Cal M}_1}$, for all
$g,h\in {C}^{\infty}({\Cal M}_2),$
 where $({F^*}(g))(x){:=} g(F(x))$, for any $x\in\Cal M_1$, and
$\{\ ,\ \}_{{\Cal M}_1}$, $\{\ ,\ \}_{{\Cal M}_2}$ are the Poisson brackets on
${\Cal M}_1$ and ${\Cal M}_2$ respectively.
Thus the above condition is equivalent to $\{ g,h\}_{{\Cal M}_2}\circ
F=\{g\circ F,h\circ F\}_{{\Cal M}_1}$.
\par If ${\Cal M}_1$ and ${\Cal M}_2$ are two Poisson manifolds with Poisson
structures defined by ${\omega}_1\in {\wedge}^2 T_{{\Cal M}_1}$ and
${\omega}_2\in {\wedge}^2 T_{{\Cal M}_2}$, respectively, we define the direct
product Poisson
structure on ${\Cal M}_1\times {\Cal M}_2$ as
$${\omega}_1\times {\omega}_2 {:=}{\omega}_1\times 1+ 1\times {\omega}_2 ,\tag
1.2$$
 which is a map $\: {\wedge}^2 T_{{\Cal M}_1}\oplus {\wedge}^2 T_{{\Cal
M}_2}\hookrightarrow {\wedge}^2 T_{{\Cal M}_1\times {\Cal M}_2}$.
Here the space $C^{\infty}({\Cal M}_1\times {\Cal M}_2)$ is identified with the
space $C^{\infty}({\Cal M}_1)\otimes C^{\infty}({\Cal M}_2)$(the reason being
that a Poisson structure on $C^{\infty}({\Cal M}_1\times {\Cal M}_2)$ is
uniquely defined by the one on
$C^{\infty}({\Cal M}_1)\otimes C^{\infty}({\Cal M}_2)$) under
appropriate completion of the tensor product.
In more detail, for any function $f\in C^{\infty}({\Cal M}_1\times {\Cal
M}_2)$, and for each $x\in \Cal M_1$ and $y\in \Cal M_2$ let us define the
functions
$f^x$ on $\Cal M_2$ and $f^y$ on $\Cal M_1$ as follows:
$$f^x(y)=f(x,y) \hskip 0.5in \text{and} \hskip 0.5in f^y(x)=f(x,y).$$
Then (1.2) means
$$\{f_1,f_2\}_{\Cal M_1\times\Cal M_2}(x,y)=\{f_1^x,f_2^x\}_{\Cal
M_2}(y)+\{f_1^y,f_2^y\}_{\Cal M_1}(x), \tag 1.3$$
for any two functions $f_1,f_2\in C^{\infty}(\Cal M_1\times\Cal M_2)$.
\definition{Definition 1.2 {\rm (Poisson-Lie group)[D1]}} Let $G$ be a Lie
group. Let $\omega$ be a Poisson structure
on $G$. The pair $(G,{\omega})$ is said to be a Poisson-Lie group if the
multiplication map $m\: G\times G\to G$ is Poisson, where the manifold $G\times
G$ is equipped
with the direct product Poisson structure ${\omega}\times {\omega}$.
\enddefinition
\par Let $L_x\:G\to G$ and $R_x\:G\to G$ be the left and right actions of $G$
on itself
defined by $y\mapsto xy$ and $y\mapsto yx$ respectively, where $x,y\in G$. Then
for any two functions
$f_1,f_2\in C^{\infty}(G)$ the compatibility between the product Poisson
structure on $G\times G$ introduced
by (1.2) and the Poisson structure on $G$ can be witten as
$$\eqalign{\{f_1,f_2\}_{G}(xy)&=m^{*}(\{f_1,f_2\}_{G})(x,y)\cr
&=\{m^{*}f_1,m^{*}f_2\}_{G\times G}(x,y)\cr
&=\{(m^{*}f_1)^x,(m^{*}f_2)^x\}_G(y)+\{(m^{*}f_1)^y,(m^{*}f_2)^y\}_G(x)\cr
&=\{f_1\circ L_x,f_2\circ L_x\}_G(y)+\{f_1\circ R_y,f_2\circ R_y\}_{G}(x).}\tag
1.4$$
If $(L_x)_{*y}$ and $(R_y)_{*x}$ are the
tangent maps to $L_x$ and $R_y$ evaluated at the points $y$ and $x$
respectively, we deduce
$${\omega}_{xy}=(L_x)_{*y}{\omega}_{y}+(R_y)_{*x}{\omega}_{x}  . \tag 1.5$$
In local coordinates
$$\omega _{ij}({z})=\omega _{kl}(x)\pd {z_i}{x_k} \pd {z_j}{x_l} +
\omega_{kl}{(y)}\pd {z_i}{y_k} \pd {z_j}{y_l} , \tag 1.6$$
where $z =xy $.
\par If $e\in G$ is the identity of $G$, then (1.5) yields
$2{\omega}_{e}={\omega}_{e}$. Therefore ${\omega}_e=0$.
This implies that $\omega$ is not a symplectic structure since the rank of
$\omega$ at the identity of $G$ is zero,
and we are dealing with more general Poisson manifolds.
\par Thus, locally, (1.1) and (1.6) are the defining equations of a
Poisson-Lie group. Many examples have been worked out.
They are mostly examples of Poisson-Lie structures on algebraic groups.
Presently there is no general method of classifying all Poisson-Lie structures
on a given
Lie group. Some results are known for finite dimensional Lie groups. The most
general is the classification of all Poisson-Lie structures for complex simple
Lie groups, due to Belavin and Drinfel'd [BD].
Little is known in the infinite dimensional case.
\remark{Remark} In the definitions above all manifolds were finite-dimensional
($\Cal M$, respectively the group $G$). To extend these to the
infinite-dimensional case, one needs two objects:
$T_{\Cal M}$ and $C^{\infty}(\Cal M)$. Since we shall study
infinite-dimensional groups in this text, the infinite-dimensional aspects
will be addressed at the moment they are introduced.\endremark
\par We now proceed with the definition of a Lie-bialgebra and formulate a
theorem (again due to Drinfel'd) relating the concept of a Lie-bialgebra to the
concept of a Poisson-Lie group.
\definition{Definition 1.3} A Lie-bialgebra $\Cal G$ is a Lie algebra $\Cal G$
equipped with a coalgebra map  $\alpha\:\Cal G\to\Cal G\wedge\Cal G$ such that
$\alpha$ is a 1-cocycle of $\Cal G$
with values in the $\Cal G$-module ${\wedge}^2\Cal G$, where $\Cal G$ acts on
${\wedge}^2\Cal G$ by means of the adjoint representation, and $\alpha$
satisfies the co-Jacobi identity. Thus,
$(\Cal G,\alpha)$ is a Lie bialgebra iff
$$\eqalign{(i)\qquad & \tau\circ\alpha=-\alpha\cr
(ii) \qquad & \alpha\bigl([X,Y]\bigr)=ad_X\alpha(Y)-ad_Y\alpha(X),\ \ \
X,Y\in\Cal G,\cr
(iii)\qquad & [1\otimes 1\otimes 1+(\tau\otimes
1)(1\otimes\tau)+(1\otimes\tau)(\tau\otimes 1)](1\otimes\alpha)\circ\alpha=0
,}$$
where $\tau$ is the transposition map $\tau\:\Cal G\otimes\Cal G\to\Cal
G\otimes\Cal G$ defined by $\tau(a\otimes b)=b\otimes a$, for any $a,b\in\Cal
G$.
\enddefinition
This definition encompasses the case when $\Cal G$ is
infinite-dimensional.
Condition (ii) means that $\alpha$ is a 1-cocycle in the
Chevalley-Eilenberg cohomology of Lie algebras.
Therefore we will refer to (ii) as the 1-cocycle condition in the
sequel. In the case when $\alpha =\delta\alpha^{0}$ is a
$1$-coboundary, $\alpha(X)=ad_{X}r$, where $r\in \Cal G\wedge\Cal G$
is a 0-cochain, which is referred to as the classical
$r$-matrix [STS2].
\par Let $\{e_i\}$ be a basis of $\Cal G$ and let us write $\alpha$ in this
basis as $\alpha(e_n)=\alpha_{ij}^n e_i\wedge e_j$. Let $C_{n}^{ij}$ be the
structure constants of $\Cal G$ defining the Lie structure on $\Cal G$
by $[e_i,e_j]=C_{n}^{ij}e_n$.
Property (i) in the definition of $\alpha$ implies that
$\alpha_{ij}^n=-\alpha_{ji}^n$.
Then the equation (iii) written in terms of $\alpha_{ij}^n$ becomes
$$\alpha_{ij}^n\alpha_{sp}^j+\alpha_{pj}^n\alpha_{is}^j+\alpha_{sj}^n\alpha_{pi}^j=0.\tag 1.7$$
Similarly, equation (ii) expressed in terms of $\alpha_{ij}^n$ and the
structure constants $C_{n}^{ij}$ of $\Cal G$ becomes
$$C_n^{ij}\alpha_{kl}^n=\alpha_{ml}^jC_k^{im}+\alpha_{km}^jC_l^{im}-\alpha_{ml}^iC_k^{jm}-\alpha_{km}^iC_l^{jm}. \tag 1.8$$
Thus, these two systems of equations, (1.7) and (1.8), plus the Jacobi identity
for the structure constants of the Lie algebra $\Cal G$
$$C_m^{ij}C_n^{mk}+C_m^{jk}C_n^{mi}+C_m^{ki}C_n^{mj}=0,$$
define a Lie-bialgebra structure on $\Cal G$.
Then we have the following result.
\proclaim{Theorem 1.4 {\rm \cite{D1}}} The category of connected, simply
connected finite dimensional Poisson-Lie groups is equivalent to the category
of  finite dimensional
Lie-bialgebras.\endproclaim
\par For a proof see, e.g., [LW,St]. We now prove a property of
(finite-dimensional) Poisson-Lie groups, which is usually assumed to be
part of the definition.
\proclaim{Theorem 1.5}  Let $G$ be a Poisson-Lie group. Then the map
$\varphi\:G\to G$ defined by  $\varphi(x)=x^{-1}$  is an anti-Poisson
map.\endproclaim
\demo{Proof} We prove the statement in a neighbourhood of the identity element
of $G$. Let $z_i=z_i(x,y)$, for $i=1,\ldots,n$, be the coordinate functions of
$z=xy$. The multiplicativity condition reads
$$\omega_{ij}{(z)}=\omega_{kl}{(x)}\pd {z_i}{x_k} \pd
{z_j}{x_l}+\omega_{kl}{(y)}\pd {z_i}{y_k} \pd {z_j}{y_l}.\tag 1.9$$
After solving $z_i=z_i(x,y)$ with respect to the coordinates of $y$ we have
$y_i=y_i(x,z)$. We differentiate the identities
$$y_i\equiv y_i(x,z(x,y)),\ \ \ \text{for}\ \ \ i=1,\ldots,n,$$
with respect to $y_k$ for each $k=1,\ldots,n$ to obtain
$$\delta^k_i=\pd {y_i}{z_l}\biggm|_{(x,z)}\pd {z_l}{y_k}\biggm|_{(x,y)}.\tag
1.10$$
\par Let $\varphi\:G\to G$ be the map defined by $\varphi(x)=x^{-1}$, which is
given in coordinates by the functions $\varphi_i=\varphi_i(x)$. The we have
$$0=z_i(x,\varphi(x)),\ \ \ \text{for}\ \ \ i=1,\ldots,n.\tag 1.11$$
We differentiate (1.11) with respect to $x_k$ to obtain
$$0=\pd {z_l}{x_k}\biggm|_{(x,\varphi(x))}+\pd
{z_i}{y_l}\biggm|_{(x,\varphi(x))}\pd {\varphi_l}{x_k}.\tag 1.12$$
After multiplying both sides of (1.9) by $\pd {y_m}{z_i}\Bigm|_{(x,z)}\pd
{y_n}{z_j}\Bigm|_{(x,z)}$ and summing over $i,j$ we get
$$\pd {y_m}{z_i}\biggm|_{(x,z)}\pd
{y_n}{z_j}\biggm|_{(x,z)}\omega_{ij}{(z)}=\pd {y_m}{z_i}\biggm|_{(x,z)}\pd
{y_n}{z_j}\biggm|_{(x,z)}\omega_{kl}{(x)}\pd {z_i}{x_k} \pd
{z_j}{x_l}+\omega_{mn}{(y)}, \tag 1.13$$
where we used (1.10).
\par We now set $z=e=x\varphi(x)$ in (1.13), and obtain
$$0=\pd {y_m}{z_i}\biggm|_{(x,e)}\pd
{y_n}{z_j}\biggm|_{(x,e)}\omega_{kl}{(x)}\pd {z_i}{x_k}
\biggm|_{(x,\varphi(x))}\pd
{z_j}{x_l}\biggm|_{(x,\varphi(x))}+\omega_{mn}{(\varphi(x))}.\tag 1.14$$
Using (1.12) the above equality is equivalent to
$$0=\pd {y_m}{z_i}\biggm|_{(x,e)}\pd
{y_n}{z_j}\biggm|_{(x,e)}\omega_{kl}{(x)}\pd
{z_i}{y_p}\biggm|_{(x,\varphi(x))}\pd {\varphi_p}{x_k}\pd
{z_j}{y_s}\biggm|_{(x,\varphi(x))}\pd
{\varphi_s}{x_l}+\omega_{mn}{(\varphi(x))}.$$
Now using again (1.19) as
$$\delta^k_i=\pd {y_i}{z_l}\biggm|_{(x,e)}\pd
{z_l}{y_k}\biggm|_{(x,\varphi(x))},$$
we finally conclude that
$$\omega_{mn}{(\varphi(x))}=-\omega_{kl}{(x)}\pd {\varphi_m}{x_k}\pd
{\varphi_n}{x_l}.\qed$$\enddemo
\par In the following sections we will adapt the above given definitions for
the case where $G$ will stand for particular infinite-dimensional groups, and
address the classification problem. We will
show that theorems analogous to Theorems 1.4 and 1.5 above also hold in this
case.
\head 2. The group of infinite-jets $G_{\infty}$ and its Lie algebra\endhead
\par Let $G_{\infty}=\left\{x=(x_1,x_2,\ldots)\in
\Bbb R^{\infty}\mid x\ne 0\right\}\subset\Bbb R^{\infty}$ be a subset
of the set of
infinite
sequences of real numbers. (We may adopt a purely formal point of view
and
take sequences of letters (indeterminates) which we interpret as
generators of an algebra of ``functions'' on the group. This is done in Sec\. 6
where the group of diffeomorphisms of the line is treated as a formal
group. For the group of diffeomorphisms of a point, introduced below,
both
points of view are possible, and they lead to the same results, since
our treatment is mostly algebraic in nature.) For each $x\in
G_{\infty}$ consider
 the formal power series, ${x}(u)=\sum_{i=1}^{\infty}{x_i u^i}$,  in the
variable $u$. This defines a bijective map
from $G_{\infty}$ into the group of $\infty$-jets, maps from $\Bbb
R^1\to\Bbb
R^1$ at $0\in\Bbb
R^1$, as follows. Define a multiplication $m_{\infty}\: {G_{\infty}}\times
{G_{\infty}}\to {G_{\infty}}$ on $G_{\infty}$
induced by the substitution of formal power series. For any
${x},{y}\in G_{\infty}$ define $xy\in G_{\infty}$ by
$({x}{y})(u):={x}({y}(u))$.
 The induced multiplication makes $G_{\infty}$ a group with an
identity $e=(1,0,0,\ldots)$ or $e(u)=u$. That is, $e\:\Bbb R^1\to\Bbb
R^1$ is the identity map. The associativity of multiplication is implied
by the associativity of substitution of power series. In coordinates
$$z _k =\sum_{i=1}^k
x_i\sum_{(\sum_{\alpha=1}^{i}{j_{\alpha}})=k}{y_{j_1}\ldots y_{j_i}},\tag 2.1$$
where $z(u)=(xy)(u)=\sum_{i=1}^{\infty}z_iu^i$.
The first several formulae are given below
$$\eqalign{{z}_1&=x_1 y_1\cr {z}_2&=x_1y_2+x_2y^2_1\cr
{z}_3&=x_1y_3+x_2 2y_1y_2+x_3y^3_1\cr
{z}_4&=x_1y_4+x_2(y^2_2+2y_1y_3)+x_3 3y^2_1y_2+x_4y^4_1\cr\vdots\cr
{z}_n&=x_1y_n+x_ny^n_1+y_{n-1}2y_1x_2+x_{n-1}(n-1)y^{n-2}_1y_2+{
O}(<n-1),\hskip 0.1in n>3\cr
\vdots .\cr}$$
\par The group so obtained is
the group of formal diffeomorphisms ($\infty$-jets) of the line,
leaving the origin fixed.
It can be viewed as a projective limit of a family of
finite dimensional Lie groups with a smooth structure introduced as follows.
Consider the family of Lie groups and maps
$(G_n,{\pi}_{n+1,n})_{n\in\zN}$, where $G_n=\zR^n\setminus M_n$ and
$M_n=\{\,(x_1,\ldots,x_n)\in\zR^n\mid x_1=0\,\}$. The  multiplication
$m_n\:G_n\times G_n\to G_n$ is induced by substitution
$\bigl(\Cal X_n\Cal Y_n\bigr)(u)=\Cal X_n\bigl(\Cal Y_n(u)\bigr)\mod{u^{n+1}}$,
where $\Cal X_n$ and $\Cal Y_n$ are
polynomials in $u$ of degree $n$.
That is, the group $G_n$ is an open subset of $\zR^n$,
that carries the structure of a finite dimensional $C^{\infty}$ manifold
modeled on $\zR^{n}$.
The  maps $\pi_{n+1,n}\:G_{n+1}\to G_n$ defined by
$\pi_{n+1,n}(x_1,\ldots,x_{n+1})=(x_1,\ldots,x_n)$
 are homomorphisms, i.e\.,
$\pi_{n+1,n}\circ m_{n+1}=m_n\circ (\pi_{n+1,n}\times\pi_{n+1,n}).$
 This follows from the definition of $\pi_{n+1,n}$
and (2.1). The family $(G_n,{\pi}_{n+1,n})_{n\in\zN}$ has a projective limit
$(G_{\infty},(\pi_{\infty,n})_{n\in\zN})$, where
$G_{\infty}=\{x\in\zR^{\infty}\mid x_1\ne 0\}$ is an open subset of
$\zR^{\infty}$. The maps $\pi_{\infty,n}\:G_{\infty}\to G_n$  are defined by
$\pi_{\infty,n}(x_1,\ldots,x_n,x_{n+1},\ldots)=(x_1,\ldots,x_n).$
 Obviously, these maps satisfy
$\pi_{\infty,n}=\pi_{n+1,n}\circ\pi_{\infty,n+1}$, and are homomorphisms,
$\pi_{\infty,n}\circ m_{\infty}=m_n\circ
(\pi_{\infty,n}\times\pi_{\infty,n})$, where $m_{\infty}\:G_{\infty}\times
G_{\infty}\to G_{\infty}$ is defined by (2.1).
\par Let us consider now the family of spaces and maps
$(C^{\infty}(G_n),\pi^*_{n,n+1})_{n\in\zN}$, where the maps
$\pi^*_{n,n+1}\:C^{\infty}(G_n)\to C^{\infty}(G_{n+1})$ are defined
by $(\pi^*_{n,n+1}(f))(x_1,\ldots,x_n,x_{n+1})=f(x_1,\ldots,x_n)$, for any
$f\in C^{\infty}(G_n)$.
Then the above family has an inductive limit
$(C^{\infty}(G_{\infty}),\pi^*_{n,\infty})$, where
$\pi^*_{n,\infty}\:\allowbreak C^{\infty}(G_n)\to C^{\infty}(G_{\infty})$ is
defined by
$((\pi^*_{n,\infty})(f))(x)=f(x_1,\ldots,x_n)$
 for any $x\in G_{\infty}$ and $f\in C^{\infty}(G_n)$. Thus, by
definition, the space $C^{\infty}(G_{\infty})$ of smooth functions on
$G_{\infty}$ is the space of
 smooth functions (of {\it finite} number of variables) on $\zR^{\infty}$,
restricted to $\zR^{\infty}\setminus M$, where $M=\{x\in\zR^{\infty}\mid
x_1=0\}$.
\par One can define the Lie algebra of $G_{\infty}$ in different ways. Probably
the most efficient one is as the Lie algebra of
derivations (smooth vector fields) of the algebra $C^{\infty}(G_{\infty})$.
These are of the form
$$X=\sum_{i=1}^{\infty}v_i\frac {\partial}{\partial x_i},\ \ \ \ v_i\in
C^{\infty}(G_{\infty}).\tag 2.2$$
Note that if $f\in C^{\infty}(G_{\infty})$, then
$X(f)=\sum_{i=1}^{n}v_i\frac {\partial f}{\partial x_i}$
is a finite sum, for some $n\in\zN$, since $f$ depends only upon a finite
number of variables. We also have $X(f)\in C^{\infty}(G_{\infty})$. Every
automorphism $\varphi\:G_{\infty}\to G_{\infty}$ acts on the space of
derivations by $(\varphi_*X)=(\varphi_*^{-1})X \varphi^*$, and
on $C^{\infty}(G_{\infty})$ it acts by $(\varphi^*f)(x)=f(\varphi(x))$. Since
the functions $f\in C^{\infty}(G_{\infty})$ are functions
of finite number of variables it is enough to describe the map $\varphi_*$ on
vector fields restricted to $C^{\infty}(G_{n})$ for each
$n\in\zN$.
\proclaim{Lemma 2.11} The set $\{X_n\}_{n\ge 1}$ of left-invariant vector
fields on $G_{\infty}$ is given by
$$X_n=\sum_{i=1}^{\infty}ix_i\frac {\partial}{\partial x_{i+n-1}}.\tag 2.3$$
\endproclaim
\demo{Proof} From (2.1) the map $y\mapsto xy$ is given by
$$z_n=(xy)_n=x_1y_n+y_1^nx_n+\sum_{i=2}^{n-1}
x_i\sum_{(\sum_{\alpha=1}^{i}{j_{\alpha}})=n}{y_{j_1}\ldots y_{j_i}},\ \ \
\text{for each}\  n\ge 1.\tag 2.4$$
The matrix of the tangent to the map defined by (2.4) is $\frac {\partial
z_n}{\partial y_m}\Bigm|_{y=e}$. The only terms in
$$\sum_{i=2}^{n-1} x_i\sum_{(\sum_{\alpha=1}^{i}{j_{\alpha}})=n}{y_{j_1}\ldots
y_{j_i}},\ \ \
\text{for each}\  n\ge 1,$$
that would contribute to the tangent map are the ones for which the product
$y_{j_1}\ldots y_{j_i}$ has exactly $(i-1)$ multiples
equal to $y_1$ and the one remaining equal to $y_{j_{\alpha}}$ for some
$\alpha$, $2\le\alpha\le i$. There are exactly $\binom {i}{i-1}=i$
terms of this form. Therefore we rewrite (2.4) as
$$z_n=\sum_{i=1}^{n}ix_iy_1^{i-1}y_{n-i+1}+\ldots,$$
where the dots indicate terms that do not contribute to $\frac {\partial
z_n}{\partial y_m}\Bigm|_{y=e}$. Hence,
$$\frac {\partial z_n}{\partial
y_m}\biggm|_{y=e}=\sum_{i=1}^{n}ix_i\delta_{n-i+1}^m=(n-m+1)x_{n-m+1}.$$
If $\{\frac {\partial}{\partial y_i}\}$ is a basis of vector fields at the
identity, then
$$\align \biggl(\varphi_*\frac {\partial}{\partial
y_m}\biggr)_x=&\sum_{i=1}^{n}\frac {\partial z_i}{\partial
y_m}\biggm|_{y=e}\frac {\partial}{\partial x_i}\\
=&\sum_{i=m}^{n}(i-m+1)x_{i-m+1}\frac {\partial}{\partial x_i}\\
=&\sum_{i=1}^{n-m+1}ix_{i}\frac {\partial}{\partial x_{i+m-1}}.\endalign$$
Therefore for each $n\in\zN$, the set of vector fields $\{X_k\}_{k=1}^n$, where
$X_k=\sum_{i=1}^{n-k+1}ix_{i}\frac {\partial}{\partial x_{i+k-1}},\ \ \
\text{for}\ \ 1\le k\le n,$
forms a basis of left-invariant vector fields on $G_n$. Therefore the set
$\{X_n\}_{n\ge 1}$, where
$X_n=\sum_{i=1}^{\infty}ix_i\frac {\partial}{\partial x_{i+n-1}},$
forms a basis of left invariant vector fields on $G_{\infty}$.\qed\enddemo
\proclaim{Lemma 2.12} Every smooth vector field on $G_{\infty}$ is generated by
the set $\{X_n\}_{n\ge 1}$ of left-invariant vector fields
{\rm (2.3)} on $G_{\infty}$.\endproclaim
\demo{Proof} Let
$Y=\sum_{i=1}^{\infty}v_i\frac {\partial}{\partial x_i}$
be a smooth vector field on $G_{\infty}$. We define inductively the following
sequence of smooth vector fields. Let
$$\align Y_1=&Y-\frac
{v_1}{x_1}X_1=Y-\psi_1X_1,\qquad\text{where}\qquad\psi_1:=\frac {v_1}{x_1},\\
Y_2=&Y_1-\psi_2X_2,\qquad \text{where}\qquad\psi_2:=\frac
{1}{x_1}\biggl(v_2-2x_2\psi_1\biggr),\\
\vdots\tag 2.5\\
Y_n=&Y_{n-1}-\psi_nX_n,\qquad\text{where}\qquad\psi_n:=\frac
{1}{x_1}\biggl(v_n-\sum_{i=2}^nix_i\psi_{n-i+1}\biggr),\\
\vdots\endalign$$
Summing up the first $n$ equalities in (2.5) we get
$Y=\sum_{i=1}^{n}\psi_iX_i+Y_n.$
By construction $Y_n$ is such that $Y_n\Bigm|_{C^{\infty}(G_n)}=0$, for any
$n\in\zN$. Hence,
$Y=\sum_{i=1}^{\infty}\psi_iX_i.\qed$\enddemo
\par We now show that $\{X_n\}_{n\ge 1}$ forms a Lie subalgebra of the Lie
algebra of vector fields on $G_{\infty}$ with a Lie bracket
given by
$$\Bigl[X_n,X_m\Bigr]=(n-m)X_{n+m-1}.\tag 2.6$$
 For that we need to compute the commutator of two left-invariant vector fields
$X_n=\sum_{i=1}^{\infty}ix_i\frac {\partial}{\partial x_{i+n-1}}$ and
$X_m=\sum_{j=1}^{\infty}jx_j\frac {\partial}{\partial x_{j+m-1}}$. Namely,
$$\align
\Bigl[X_n,X_m\Bigr]=&\sum_{i=1}^{\infty}\sum_{j=1}^{\infty}ix_ij\delta_{i+n-1}^j\frac {\partial}{\partial x_{j+m-1}}-
\sum_{i=1}^{\infty}\sum_{j=1}^{\infty}jx_ji\delta_{j+n-1}^i\frac
{\partial}{\partial x_{i+n-1}}\\
=&\sum_{i=1}^{\infty}ix_i(i+n-1)\frac {\partial}{\partial
x_{i+n+m-2}}-\sum_{j=1}^{\infty}jx_j(j+m-1)\frac {\partial}{\partial
x_{j+n+m-2}}\\
=&(n-m)\sum_{i=1}^{\infty}ix_i\frac {\partial}{\partial x_{i+n+m-2}}\\
=&(n-m)X_{n+m-1}.\endalign$$
To make correspondence with the more familiar notation we shift the
indeces by 1, and introduce $\widetilde X_{n}:=X_{n+1}$. Then
$\left[\widetilde X_{n},\widetilde X_{m}\right]=(n-m)\widetilde
X_{n+m}$, for all $n,m\ge 0$. Thus, the algebra so obtained is the
maximal subalgebra of the Witt algebra.
\par Let us assume now that $G_n$ are equipped with Poisson-Lie
structures $\{\ ,\ \}_n\:C^{\infty}(G_n)\times C^{\infty}(G_n)\to
C^{\infty}(G_n)$. It is natural to require that the projection maps
$\pi_{n+1,n}\:G_{n+1}\to G_n$ are Poisson, i.e\.,
$\pi^*_{n+1,n}(\{f,g\}_n)=\{\pi^*_{n+1,n}(f),\pi^*_{n+1,n}(g)\}_{n+1}$, where,
as above, $(\pi^*_{n+1,n}(f))(x):=f(\pi_{n+1,n}(x))$ for every $f\in
C^{\infty}(G_n)$ and $x\in G_{n+1}$.
\par For a smooth bivector field $\omega$ on $G_{\infty}$ (i.e., a
rank 2 skew-symmetric tensor such that $\omega_{ij}\in
C^{\infty}(G_{\infty})$ and $\omega_x=\omega_{ij}(x)\frac
{\partial}{\partial x_i}\wedge\frac {\partial}{\partial x_j}$) let us assume
that the
map $\{\ ,\ \}\:C^{\infty}(G_{\infty})\times
C^{\infty}(G_{\infty})\to C^{\infty}(G_{\infty})$ gives a Poisson-Lie
structure on $G_{\infty}$ defined by $\{f,g\}(x)=\omega_{ij}(x)\frac
{\partial f}{\partial x_i}\frac
{\partial g}{\partial x_j}$ for any $f,g\in C^{\infty}(G_{\infty})$. Then it is
also natural to require that the maps $\pi_{\infty,n}$ are Poisson.
That is, we want the condition
$\pi^*_{n,\infty}\bigl(\{f,g\}_n\bigr)=\{\pi^*_{n,\infty}(f),\pi^*_{n,\infty}(g)\}$
to be satisfied, for any $f,g\in C^{\infty}(G_n)$. Then one defines
$(G_{\infty},\{\ ,\ \},(\pi_{\infty,n})_{n\in\zN})$ to be the projective limit
of the family of Poisson-Lie groups and maps $(G_n,\{\ ,\
\}_n,\pi_{n+1,n})_{n\in\zN}$. Later it will be clear
that for the Poisson-Lie groups studied in this text these conditions are
automatically satisfied.
\par Are there any Poisson-Lie structures on $G_{\infty}$?
If such structures exist, can they be classified?
 Also, since for any finite $n$ there is a one-to-one correspondence between
the Poisson-Lie structures on $G_n$ and the Lie-bialgebra structures on the Lie
algebra $\Cal G_n$ of $G_n$, one is led to enquire if there are any
Lie-bialgebra structures on the Lie algebra ${\Cal G}_{\infty}$
of $G_{\infty}$. The same questions exist for the group $G_{0\infty}$
of formal diffeomorphisms of the line without fixed points and its Lie algebra
$\Cal
G_{0\infty}$, which is the Witt algebra.
It turns out that all these questions can be fully answered.
Let us turn our attention to the Lie-algebraic picture first.
 Our first goal is to find all bialgebra structures on $\Cal
G_{\infty}$ and $\Cal
G_{0\infty}$. The next section is devoted to this.
\head 3. Bialgebra structures on the
Witt algebra $\Cal G_{0\infty}$ and its principal subalgebra $\Cal
G_{\infty}$
\endhead

{ In this section we use completely elementary methods to compute the
Shevalley-Eilenberg
cohomology groups $H^1(\Cal G_{0\infty},\Cal
G_{0\infty}\widehat\wedge\Cal G_{0\infty})$ and $H^1(\Cal G_{\infty},\Cal
G_{\infty}\widehat\wedge\Cal G_{\infty})$ with coefficients in the
completion of the tensor product under the natural grading induced from the
gradings on $\Cal G_{0\infty}$ and $\Cal G_{\infty}$. We show that both are
trivial.
Thus, we find all 1-cocycles on the Lie algebra $\Cal
G_{0\infty}$ (the Witt algebra), and all 1-cocycles on the Lie algebra $\Cal
G_{\infty}$ of the
group $G_{\infty}$. All of them are coboundaries, and they are all
explicitly enumerated.} The algebras $\Cal G_{0\infty}$ and $\Cal
G_{\infty}$ are taken over $\Bbb R$ (or $\Bbb C$).
\par Let \{${e_i}$\}$_{i{\ge}-1}$  be a basis of the Lie algebra
$\Cal G_{0\infty}$ (the Witt algebra).
We recall that the Lie algebra structure on $\Cal G_{0\infty}$ is defined by
$$ \left[e_n,e_m\right]=(n-m)e_{n+m}\qquad n,m\ge -1. \tag 3.1$$
 To find the 1-cocycles we turn to the 1-cocycle equation
$$\alpha(\left[e_n,e_m\right])=e_n.\alpha(e_m)-e_m.\alpha(e_n), \tag 3.2$$
 which we rewrite as
$$(n-m)\alpha(e_{n+m})=e_n.\alpha(e_m)-e_m.\alpha(e_n)\qquad
n,m\ge -1. \tag 3.3$$
\par Let $\alpha(e_n)=\sum_{i,j=-1}^{\infty}\alpha^n_{ij}
e_i\wedge e_j\in \Cal
G_{0\infty}\widehat\wedge\Cal
G_{0\infty}=\bigoplus_{k}\oplus_{i+j=k}\Cal G_{0i}\wedge\Cal G_{0j}$,
where
$\Cal G_{0i}$ is the 1-dimensional subspase of $\Cal G_{0\infty}$
generated by $e_i$. Then (3.3) is equivalent
to the following infinite system of linear equations
$$(n-m)\alpha^{n+m}_{ij}=(2n-i) \alpha^{m}_{i-n,j}+
(2n-j)\alpha^{m}_{i,j-n}-(2m-i)\alpha^{n}_{i-m,j}-
(2m-j)\alpha^{n}_{i,j-m},\tag 3.4$$
where $n,m,i,j\ge -1$. Therefore, to find all 1-cocycles, one has to
describe all solutions of this system.
We now formulate the first main theorem of this section.
\proclaim{Theorem 3.1} All 1-cocycles on the Lie algebra
$\Cal G_{0\infty}$ are coboundaries.\endproclaim
\demo{Proof} We observe that (3.4) consists of two independent
subsystems. Namely, one for $\alpha_{ij}^n$'s with $n\ne i+j$ and the
other for $\alpha_{ij}^n$'s with $n=i+j$.\enddemo
\par Set $m=0$ in (3.4). (Note that (3.4) is symmetric with respect to
$m\to n$, $n\to m$.) Then from (3.4) we deduce that
$$(n-i-j)\alpha^n_{ij}=(2n-i)\alpha^0_{i-n,j}+(2n-j)\alpha^0_{i,j-n},
\qquad \text{for every}\ \ \ n,i,j\ge -1. \tag 3.5$$
\par (a) In case $n\ne i+j$ we obtain from (3.5) immediately the
solution of the first subsystem
$$\alpha_{ij}^n=\frac {(2n-i)}{(n-i-j)}\alpha_{i-n,j}^0+\frac {(2n-j)}
{(n-i-j)}\alpha_{i,j-n}^0.$$
This means that the coalgebra structure constants $\alpha_{ij}^n$, where
$n\ne i+j$, are given in terms of $\alpha_{ij}^0$,
$i,j\in\zZ_{+}$, which are arbitrary.
\par (b) The case $n=i+j$ requires more thorough analysis. Set $j=n+m-i$
in (3.4). Then  we have
$$(n-m)\alpha_{i,n+m-i}^{n+m}=(2n-i)\alpha^m_{i-n,m+n-i}+
(n-m+i)\alpha^m_{i,m-i}-(2m-i)\alpha^n_{i-m,n+m-i}-(m-n+i)\alpha_{i,n-i}^n,
\tag 3.6$$
where $n,m,i\ge -1$.
It turns out that it is enough to investigate (3.6) for a few values
of $m$ and $i$
in order to obtain a complete set of reccurence relations sufficient
to find the general solution for $\alpha_{i,n-i}^{n}$. First we set
$m=0$ in (3.6). This implies that
$0=(2n-i)\alpha_{i-n,n-i}^0+(n+i)\alpha_{i,-i}^0$, for every $i,n\ge -1$.
In particular, if $n=i$ we obtain $0=i\alpha_{-i,i}^0$, for every $i\ge
-1$, since $\alpha_{0,0}^0=0$ by the antisymmetry of
$\alpha_{ij}^n(=-\alpha_{ji}^n)$.
Therefore $\alpha_{-i,i}^0=0$, for every $i\ge -1$.
Set $m=-1=i$ in (3.6) to obtain
$$(n+1)\alpha_{-1,n}^{n-1}=(2n+1)\alpha_{-1-n,n}^{-1}+n\alpha_{-1,0}^{-1}+
\alpha_{0,n}^n+(n+2)\alpha_{-1,n+1}^n.\tag 3.7$$
Then we have the following lemma.
\proclaim{Lemma 3.2} One has $\alpha_{-1,n+1}^n=\frac {1}{n+2}\bigl[
\frac {n(n+1)}{2}\alpha_{0,-1}^{-1}-\sum_{i=1}^n\alpha_{0,i}^i\bigr]$,
for every $n\ge 1$.\endproclaim
\demo{Proof} For $n=-1,0$, (3.7) is identically satisfied.
 For $n=1$ we have from (3.7):
$2\alpha_{-1,1}^0=\alpha_{-1,0}^{-1}+\alpha_{0,1}^1+3\alpha_{-1,2}^1,$
where we used that $\alpha_{-2,1}^{-1}=0$. From this it follows that
$\alpha_{-1,2}^1=\frac
{1}{3}\left[\alpha_{0,-1}^{-1}-\alpha_{0,1}^1\right]$,
Assume that
$\alpha_{-1,k+1}^k=\frac {1}{k+2}\left[\frac
{k(k+1)}{2}\alpha_{0,-1}^{-1}-\sum_{i=1}^k\alpha_{0,i}^i\right]$, for
all $k$,
 such that $1\le k\le n-1.$
Then from (3.7), using the induction hypothesis we have
$$\align \alpha_{-1,n+1}^n=&\frac
{1}{(n+2)}\left[n\alpha_{0,-1}^{-1}-\alpha_{0,n}^n+(n+1)\alpha_{-1,n}^{n-1}
\right]\\
=&\frac {1}{(n+2)}\left[n\alpha_{0,-1}^{-1}-\alpha_{0,n}^n+
\frac
{n(n-1)}{2}\alpha_{0,-1}^{-1}-\sum_{i=1}^{n-1}\alpha_{0,i}^i\right]\\
=&\frac {1}{(n+2)}\left[n\alpha_{0,-1}^{-1}+\frac {n(n-1)}{2}
\alpha_{0,-1}^{-1}-\sum_{i=1}^n\alpha_{0,i}^i\right]\\
=&\frac {1}{(n+2)}\left[\frac
{n(n+1)}{2}\alpha_{0,-1}^{-1}-\sum_{i=1}^n\alpha_{0,i}^i\right].\qed
\endalign$$
\enddemo
\proclaim{Lemma 3.3} One has
$\alpha_{0,2}^2=\alpha_{0,1}^1-\alpha_{0,-1}^{-1}$, and
$\alpha_{1,n-1}^n=\frac {(n+1)}{2}\left[-\alpha_{0,-1}^{-1}+
\alpha_{0,n-1}^{n-1}-\alpha_{0,n}^n\right]$, for every $n\ge 3.$
\endproclaim
\demo{Proof} From (3.6), after setting $m=-1,i=0$, we obtain
$$(n+1)\alpha_{0,n-1}^{n-1}=2n\alpha_{-n,n-1}^{-1}+(n+1)\alpha_{0,-1}^{-1}+
2\alpha_{1,n-1}^n+(n+1)\alpha_{0,n}^n.\tag 3.8$$
For $n=-1,0,1$, (3.8) is identically satisfied.
Next, for $n=2$ we obtain from (3.8) that
$$\alpha_{0,2}^2=\alpha_{0,1}^1-\alpha_{0,-1}^{-1},\tag 3.9$$
since $\alpha_{-2,1}^{-1}=0$. Finally, since $\alpha_{-n,n-1}^{-1}=0$
for $n\ge 2$, we obtain from (3.8) the second assertion of the Lemma.\qed
\enddemo
\proclaim{Lemma 3.4} The formula $\alpha_{2,n-2}^n=\frac
{n(n+1)}{6}\left[\alpha_{0,n-2}^{n-2}-
2\alpha_{0,n-1}^{n-1}+\alpha_{0,n}^n\right]$ is valid for every $n\ge 5.$
\endproclaim
\demo{Proof} In (3.6) set $m=-1,i=1$. Then we have
$$(n+1)\alpha_{1,n-2}^{n-1}=(2n-1)\alpha_{1-n,n-2}^{-1}+3\alpha_{2,n-2}^n+
n\alpha_{1,n-1}^n.\tag 3.10$$
For $n=-1,0$ this formula yields trivial identitites. For $n=1$ we
obtain that
$3\alpha_{-1,2}^1=\alpha_{0,-1}^{-1}-\alpha_{0,1}^1\
(=-\alpha_{0,2}^2).$
The case $n=2$ reduces to formula (3.9), and $n=3$ yields the
identity
$0=3\alpha_{2,1}^3+3\alpha_{1,2}^3.$
For $n=4$ we obtain the relation
$$5\alpha_{1,2}^3=4\alpha_{1,3}^4.\tag 3.11$$
If $n\ge 5$, then $\alpha_{1-n,n-2}^{-1}=0$, and using Lemma 3.3 we
obtain that
$$\align \alpha_{2,n-2}^n=&\frac
{1}{3}\left[(n+1)\alpha_{1,n-2}^{n-1}-n\alpha_{1,n-1}^n\right]\\
=&\frac {1}{3}\left[\frac {n(n+1)}{2}\left(-\alpha_{0,-1}^{-1}+
\alpha_{0,n-2}^{n-2}-\alpha_{0,n-1}^{n-1}\right)-\frac {n(n+1)}{2}
\left(-\alpha_{0,-1}^{-1}+\alpha_{0,n-1}^{n-1}-\alpha_{0,n}^n\right)\right]\\
=&\frac {n(n+1)}{6}\left[\alpha_{0,n-2}^{n-2}-2\alpha_{0,n-1}^{n-1}+
\alpha_{0,n}^n\right].\qed\endalign$$
\enddemo
 \par The result obtained in Lemma 3.4 is suggestive in finding a
general formula for $\alpha_{i,n-i}^n$ (i.e., $\alpha_{ij}^n$ with
$n=i+j$) in terms of $\alpha_{0,n}^n,\ldots,\alpha_{0,n-i}^{n-i}$,
for every $i\ge 2$ and $n\ge i+3$.
\proclaim{Lemma 3.5} For every $i\ge 2$, and $n\ge i+3$ the following
formula is valid:
$$\alpha_{i,n-i}^n=\frac {1}{i+1}\binom {n+1}{i}\sum_{k=0}^i
\binom {i}{k}(-1)^{i-k}\alpha_{0,n-k}^{n-k}.\tag 3.12$$
\endproclaim
\demo{Proof} We use an induction in $i$. For $i=2$ the statetment is
true by Lemma 3.4.
Assume that (3.12) is true
for all $\alpha_{j,n-j}^n$, where $2\le j\le i$. We now proceed in
proving that it is true for $j=i+1$ and all $n\ge i+4$. Setting $m=-1$
in (3.6) and using that $i\ge 2$ and $n\ge i+3$ we
obtain
$$(n+1)\alpha_{i,n-i-1}^{n-1}-(n+1-i)\alpha_{i,n-i}^n=(i+2)\alpha_{i+1,n-i-1}^n.\tag
3.13$$
Then we use induction hypothesis in the left-hand-side of (3.13) and deduce
that
$$\align \alpha_{i+1,n-i-1}^n=&\frac {1}{i+2}\left[\frac
{(n+1)}{i+1}\binom {n}{i}\sum_{k=0}^i
\binom {i}{k}(-1)^{i-k}\alpha_{0,n-k-1}^{n-k-1}-\frac
{(n+1-i)}{i+1}\binom {n+1}{i}\sum_{k=0}^i
\binom {i}{k}(-1)^{i-k}\alpha_{0,n-k}^{n-k}\right]\\
=&\frac {1}{i+2}\binom {n+1}{i+1}\left[\sum_{k=0}^i
\binom {i}{k}(-1)^{i-k}\alpha_{0,n-k-1}^{n-k-1}-\sum_{k=0}^i
\binom {i}{k}(-1)^{i-k}\alpha_{0,n-k}^{n-k}\right]\\
=&\frac {1}{i+2}\binom {n+1}{i+1}\left[\sum_{k=1}^{i+1}
\binom {i}{k-1}(-1)^{i-k+1}\alpha_{0,n-k}^{n-k}+\sum_{k=0}^i
\binom {i}{k}(-1)^{i-k+1}\alpha_{0,n-k}^{n-k}\right]\\
=&\frac {1}{i+2}\binom {n+1}{i+1}\left[\sum_{k=1}^{i}\left[\binom
{i}{k-1}+\binom {i}{k}\right](-1)^{i-k+1}\alpha_{0,n-k}^{n-k}+
\alpha_{0,n-i-1}^{n-i-1}+(-1)^{i+1}\alpha_{0,n}^{n}\right]\\
=&\frac {1}{i+2}\binom {n+1}{i+1}\sum_{k=0}^{i+1}\binom
{i+1}{k}(-1)^{i+1-k}\alpha_{0,n-k}^{n-k},\endalign$$
which concludes the proof. It is also an immediate corollary that the
relation (3.11) is identically satisfied.\qed\enddemo
\par Let $m=1,i=-1$ in (3.6). Then we have
$$(n-1)\alpha_{-1,n+2}^{n+1}=(2n+1)\alpha_{-1-n,n+2}^1+(n-2)\alpha_{-1,2}^1+
n\alpha_{-1,n+1}^n.\tag 3.14$$
For $n=-1$ (3.14) yields
$\alpha_{-1,2}^1=\frac
{1}{3}\left(\alpha_{0,-1}^{-1}-\alpha_{0,1}^1\right)$.
For $n=0,1$, it is identically satisfied.
For $n=2$, we obtain the relation
$$\alpha_{-1,4}^3=2\alpha_{-1,3}^2.\tag 3.15$$
Continuing further we have
$$\align 2\alpha_{-1,5}^4=&\alpha_{-1,2}^1+3\alpha_{-1,4}^3,\ \ \ \
(n=3),\\
3\alpha_{-1,6}^5=&2\alpha_{-1,2}^1+4\alpha_{-1,5}^4,\ \ \ \ (n=4),\\
&\vdots\endalign$$
which, as we shall see later, are all identically satisfied.
\proclaim{Lemma 3.6} The following formula is valid for every $n\ge
2$:
$$\alpha_{0,n}^n=a_n(\alpha_{0,1}^1-\alpha_{0,-1}^{-1}),\tag 3.16$$
where $a_2=1$, $a_3=3$, and the coefficients $a_n$, $(n\ge 4)$ are
computed by the recursive formula
$$a_{n+1}=\frac {2n}{(n-1)(n+2)}+\frac {2(n+1)}{(n+2)}a_n-\frac
{(n+1)(n-2)}{(n-1)(n+2)}a_{n-1},\tag 3.17$$
which is valid for all $n\ge 3$.\endproclaim
\demo{Proof} For $n=2$ the statement is true according to formula
(3.9). In order to compute $a_3$ we use the relation (3.15) and Lemma
3.2. From Lemma 3.2 we have
$$\alpha_{-1,3}^2=\frac
{1}{2}\left[2\alpha_{0,-1}^{-1}-\alpha_{0,1}^1\right],\quad
\text{and}\quad\alpha_{-1,4}^3=\frac
{1}{5}\left[7\alpha_{0,-1}^{-1}-2\alpha_{0,1}^1-\alpha_{0,3}^3\right].$$
Then from (3.15) it follows that
$\alpha_{0,3}^3=3(\alpha_{0,1}^1-\alpha_{0,-1}^{-1})$.
After setting $m=1,i=1$ in (3.6)
one obtains
$$(n-1)\alpha_{1,n}^{n+1}=-n\alpha_{0,1}^1-\alpha_{0,n}^n+(n-2)\alpha_{1,n-1}^n,\qquad \text{for }\ \ n\ge 3,\tag 3.18$$
since $\alpha_{1-n,n}^1=0$, for $n\ge 3$. Next, we use Lemma 3.3 to
write formulae for
$$\alpha_{1,n-1}^n=\frac
{(n+1)}{2}\left[-\alpha_{0,-1}^{-1}+\alpha_{0,n+1}^{n+1}-\alpha_{0,n}^n\right],\quad
\text{and}\quad
\alpha_{1,n}^{n+1}=\frac
{(n+2)}{2}\left[-\alpha_{0,-1}^{-1}+\alpha_{0,n}^{n}-\alpha_{0,n+1}^{n+1}\right].$$
Substituting the above formulae into (3.18) we obtain, after some
algebra, that
$$\alpha_{0,n+1}^{n+1}=\frac
{2n}{(n-1)(n+2)}\left(\alpha_{0,1}^1-\alpha_{0,-1}^{-1}\right)+
\frac {2(n+1)}{(n+2)}\alpha_{0,n}^{n}-\frac
{(n+1)(n-2)}{(n-1)(n+2)}\alpha_{0,n-1}^{n-1},\qquad\text{for all}\ \
n\ge 3.\tag 3.19$$
Now, we make an induction hypothesis. Namely, we assume that for all
$k$, $2\le k\le n$, it is true that
$\alpha_{0,k}^{k}=a_k(\alpha_{0,1}^1-\alpha_{0,-1}^{-1})$, for some
$a_k$. Then from (3.19) it follows that
$\alpha_{0,n+1}^{n+1}=a_{n+1}(\alpha_{0,1}^1-\alpha_{0,-1}^{-1})$,
where
$$a_{n+1}=\frac {2n}{(n-1)(n+2)}+\frac {2(n+1)}{(n+2)}a_n-\frac
{(n+1)(n-2)}{(n-1)(n+2)}a_{n-1}.$$
The first several coefficients are: $a_2=1$, $a_3=3$, $a_4=5$,
$a_5=\frac {64}{9}$, $a_6=\frac {28}{3}$, $a_7=\frac
{451}{45},\ldots$. This concludes the proof.\qed\enddemo
\proclaim{Corollary 3.7} From Lemma {\rm 3.6} and Lemma {\rm 3.5} it follows
that
$$\alpha_{i,n-i}^n=\frac {1}{i+1}\binom
{n+1}{i}\sum_{k=0}^i\binom
{i}{k}(-1)^{i-k}a_{n-k}(\alpha_{0,1}^1-\alpha_{0,-1}^{-1}),\tag
3.20$$
for every $i\ge 2$ and $n\ge i+3$.\endproclaim
\proclaim{Lemma 3.8} One has the following relation:
$\alpha_{0,-1}^{-1}=\alpha_{0,1}^1$.
\endproclaim
\demo{Proof} In (3.6) set $m=0,i=-1$. Then one has
$$(n-2)\alpha_{-1,n+3}^{n+2}=(2n+1)\alpha_{-1-n,n+3}^2+(n-3)\alpha_{-1,3}^2+(n-1)\alpha_{-1,n+1}^n.\tag
3.21$$
For $n=-1$ we obtain from (3.21) that
$-3\alpha_{-1,2}^1=-\alpha_{0,2}^2-4\alpha_{-1,3}^2-2\alpha_{-1,0}^{-1},$
which, after using Lemma 3.2, turns into an
identity. For $n=0$ we have an identity. For $n=1$ we obtain
$-\alpha_{-1,4}^3=-2\alpha_{-1,3}^2,$
which leads to an identity after using Lemma 3.2 and Lemma 3.6. The
case $n=2$ leads to an identity.
For $n=3$ we obtain the non-trivial relation
$$\alpha_{-1,6}^5=2\alpha_{-1,4}^3.\tag 3.22$$
We use now Lemma 3.2 and Lemma 3.6 to reduce both sides of the
relation (3.22) which leads to
$$\align \frac
{1}{7}\left[15\alpha_{0,-1}^{-1}-\sum_{i=1}^5\alpha_{0,i}^i\right]&=\frac
{2}{5}\left[6\alpha_{0,-1}^{-1}-\sum_{i=1}^3\alpha_{0,i}^i\right]\\
\frac
{1}{7}\left[15\alpha_{0,-1}^{-1}-\alpha_{0,1}^1-\left(1+3+5+\frac
{64}{9}\right)\left(\alpha_{0,1}^1-\alpha_{0,-1}^{-1}\right)\right]&=\frac
{2}{5}\left[6\alpha_{0,-1}^{-1}-\alpha_{0,1}^1-\left(1+3\right)\left(\alpha_{0,1}^1-\alpha_{0,-1}^{-1}\right)\right]\\
\frac
{1}{7}\left[15\alpha_{0,-1}^{-1}-\frac {154}{9}\alpha_{0,1}^1+\frac
{145}{9}\alpha_{0,-1}^{-1}\right]&=\frac
{2}{5}\left[6\alpha_{0,-1}^{-1}-5\alpha_{0,1}^1+4\alpha_{0,-1}^{-1}\right]\\
\frac
{1}{7}\left[\frac {280}{9}\alpha_{0,-1}^{-1}-\frac
{154}{9}\alpha_{0,1}^1\right]&=\frac
{2}{5}\left[10\alpha_{0,-1}^{-1}-5\alpha_{0,1}^1\right]\\
\frac {40}{9}\alpha_{0,-1}^{-1}-\frac
{22}{9}\alpha_{0,1}^1&=4\alpha_{0,-1}^{-1}-2\alpha_{0,1}^1\\
\frac {4}{9}\alpha_{0,-1}^{-1}&=\frac {4}{9}\alpha_{0,1}^1\\
\alpha_{0,-1}^{-1}&=\alpha_{0,1}^1.\endalign$$
This concludes the proof.\qed\enddemo
\par Immediately several corollaries follow.
\proclaim{Corollary 3.9} One has $\alpha_{0,n}^n=0$, for all $n\ge
2$.\endproclaim
\proclaim{Corollary 3.10} One obtains that $\alpha_{-1,n+1}^n=\frac
{1}{n+2}\left[\frac
{n(n+1)}{2}\alpha_{0,-1}^{-1}-\alpha_{0,1}^1\right]=\frac
{(n-1)}{2}\alpha_{0,1}^1$, for all $n\ge -1$.\endproclaim
\proclaim{Corollary 3.11} One has $\alpha_{i,n-i}^n=0$, for all $i\ge
2$ and $n\ge i+3$.\endproclaim
\proclaim{Corollary 3.12} One has $\alpha_{1,n-1}^n=-\frac
{(n+1)}{2}\alpha_{0,1}^1$, for all $n\ge 3$.\endproclaim
\par In the above formulae $\alpha_{0,1}^1$ is arbitrary. As a result
we are now able to write a formula for the general solution of (3.3).
Namely,
$$\aligned \alpha (e_n)&=\sum_{i,j=-1}^{\infty}\alpha_{ij}^ne_i\wedge
e_j\\
&=\sum_{i,j=-1,i+j\ne n}^{\infty}\alpha_{ij}^ne_i\wedge e_j+
\sum_{i,j=-1,i+j= n}^{\infty}\alpha_{ij}^ne_i\wedge e_j\\
&=\sum_{i,j=-1,i+j\ne n}^{\infty}\left[\frac
{(2n-i)}{(n-i-j)}\alpha_{i-n,j}^0+\frac
{(2n-j)}{(n-i-j)}\alpha_{i,j-n}^0\right]e_i\wedge e_j+
\sum_{i=-1}^{\infty}\alpha_{i,n-i}^ne_i\wedge e_{n-i}\\
&=e_n.\left(-\sum_{i,j=-1,i+j\ne 0}^{\infty}\frac
{1}{(i+j)}\alpha_{ij}^0e_i\wedge
e_j\right)+\alpha_{-1,n+1}^ne_{-1}\wedge
e_{n+1}+\alpha_{0,n}^ne_0\wedge e_n+\alpha_{1,n-1}^ne_1\wedge
e_{n-1}\\
&=e_n.\left(-\sum_{i,j=-1,i+j\ne 0}^{\infty}\frac
{1}{(i+j)}\alpha_{ij}^0e_i\wedge
e_j+\frac {1}{2}\alpha_{0,1}^1(1+\delta_{n,\pm 1})e_{-1}\wedge e_1\right),\ \ \
\text{for all}\ \ n\ge -1.
\endaligned\tag 3.23$$
This shows that the general solution of (3.3) is a coboundary. The proof of
Theorem 3.1 is completed.\qed
\par With different technique it has been obtained in [Dz] that the
first cohomology of $\Cal G_{0\infty}$ with coefficients in the
ordinary tensor product is trivial, which follows also from Theorem 3.1.
\par We now proceed with describing all bialgebra structures on the
Lie algebra $\Cal G_{\infty}$ of the group $G_{\infty}$. The
problem now
is to describe all solutions of
$$(n-m)\alpha(e_{n+m})=e_n.\alpha(e_m)-e_m.\alpha(e_n)\ \ \ \text{for
all}\ \ n,m\ge 0. \tag 3.24$$
Equation (3.24) is now equivalent to the following infinite system of
equations for the coalgebra structure constants $\alpha^n_{ij}$:
$$(n-m)\alpha^{n+m}_{ij}=(2n-i) \alpha^{m}_{i-n,j}+
(2n-j)\alpha^{m}_{i,j-n}-(2m-i)\alpha^{n}_{i-m,j}-
(2m-j)\alpha^{n}_{i,j-m},\tag 3.25$$
where $n,m,i,j\ge 0$.
\proclaim{Theorem 3.13} All 1-cocycles on the Lie algebra $\Cal
G_{\infty}$ are
coboundaries.
\endproclaim
\demo{Proof} The begining of the argument is similar to the one in
the proof of Theorem 3.1. Namely, the system of equations (3.25) is
split into two completely independent systems. One for the structure
constants $\alpha^n_{ij}$ with $n\ne i+j$, and one for the structure
constants $\alpha^n_{ij}$ with $n=i+j$. In the first case, after
setting $m=0$ in (3.25) we obtain:
$$(n-i-j)\alpha^n_{ij}=(2n-i)\alpha^0_{i-n,j}+(2n-j)\alpha^0_{i,j-n},
\hskip 0.5in \text{for every}\ \ \ n,i,j\ge 0. \tag 3.26$$
Since $n\ne i+j$ we obtain that
$$\alpha_{ij}^n=\frac {(2n-i)}{(n-i-j)}\alpha_{i-n,j}^0+\frac {(2n-j)}
{(n-i-j)}\alpha_{i,j-n}^0,\tag 3.27$$
i.e, the coalgebra structure constants $\alpha_{ij}^n$, where
$n\ne i+j$, are expressed in terms of $\alpha_{ij}^0$,
$i,j\ge 0$, which are arbitrary.
\par We now turn to the second case, which requires more complicated
analysis. The system of equations for the coalgebra structure
constants $\alpha^n_{ij}$ with $n=i+j$ is obtained from (3.25) by
setting $j=n+m-i$:
$$(n-m)\alpha_{i,n+m-i}^{n+m}=(2n-i)\alpha^m_{i-n,m+n-i}+
(n-m+i)\alpha^m_{i,m-i}-(2m-i)\alpha^n_{i-m,n+m-i}-(m-n+i)\alpha_{i,n-i}^n,
\tag 3.28$$
where $n,m,i\ge 0 $.
\enddemo
\par Again we split the rest of the proof into lemmas.
\proclaim{Lemma 3.14} One has
$\alpha^n_{0,n}=n\alpha_{0,1}^{1}$, for all $n\ge 0$.
\endproclaim
\demo{Proof} After setting $,m=1,i=0$ in (3.28) one obtains
$$\alpha^{n+1}_{0,n+1}=\alpha_{0,1}^{1}+\alpha^n_{0,n},\tag 3.29$$
since $n\alpha_{-n,n+1}^{1}=0$ for all $n\ge 0$. A simple inductive
argument now leads to
$\alpha^n_{0,n}=n\alpha_{0,1}^{1}$, for all $n\ge 0$.\qed
\enddemo
\proclaim{Lemma 3.15} One has
$\alpha_{1,n}^{n+1}=-(n+2)\alpha_{0,1}^{1}$
for all $n\ge 2$.
\endproclaim
\demo{Proof} Set $m=1,i=1$ in (3.28). Then one has
$$(n-1)\alpha_{1,n}^{n+1}=(2n-1)\alpha_{1-n,n}^{1}-2n\alpha_{0,1}^{1}+
(n-2)\alpha_{1,n-1}^{n}.\tag 3.30$$
Let us investigate (3.30) for small values  of $n$. For $n=0,1$ one
obtains identities.
For $n=2$ we obtain the first non-trivial relation. Namely,
$\alpha_{1,2}^3=-4\alpha_{0,1}^1$.
Clearly, since $\alpha_{1-n,n}^{1}=0$ for all $n\ge 2$, (3.30)
is equivalent to
$$(n-1)\alpha_{1,n}^{n+1}=-2n\alpha_{0,1}^{1}+
(n-2)\alpha_{1,n-1}^{n},\tag 3.31$$
for all $n\ge 2$. Again, an induction argument leads to
$\alpha_{1,n}^{n+1}=-(n+2)\alpha_{0,1}^{1}$, for all $n\ge 2$.\qed
\enddemo
\proclaim{Lemma 3.16} One has
$\alpha_{i,n}^{n+i}=-(n+2i)\alpha_{0,1}^{1}+\frac
{3(1-i^2)}{(n-i)}\alpha_{0,1}^{1}$ for all $n\ge 1$, and all $i\ge 1$, and
$n\ne i$.
\endproclaim
\demo{Proof} Set $m=i$ in (3.28) to obtain
$$(n-i)\alpha_{i,n}^{n+i}=(2n-i)\alpha_{i-n,n}^i+n\alpha_{i,0}^i-
i\alpha_{0,n}^n+(n-2i)\alpha_{i,n-i}^n.\tag 3.32$$
After using Lemma 3.14, (3.32) reduces to
$$(n-i)\alpha_{i,n}^{n+i}=-2ni\alpha_{0,1}^{1}+(n-2i)\alpha_{i,n-i}^n.\tag
3.33$$
If we set $n=i$ in (3.33) it reduces further to
$2i\alpha_{0,i}^i=2i^2\alpha_{0,1}^{1}$,
from where for $i\ne 0$ it is equivalent to the assertion of Lemma 3.14.
Let $n\ne i$.
Then the  general solution of (3.33) is given by
$$\alpha_{i,n}^{n+i}=-(n+2i)\alpha_{0,1}^{1}+f(n,i),\tag 3.34$$
where $f(n,i)$ satisfies the equation
$(n-i)f(n,i)=(n-2i)f(n-i,i)$.
If we define $\alpha (n,i):=(n-i)f(n,i)$ then this equation is equivalent to
$\alpha (n,i)=\alpha (n-i,i)$.
 From here it follows that $\alpha (n,i)=\alpha (i)$, i.e., $\alpha$ is
independent of $n$. Therefore the general solution (3.34) is given by
$$\alpha_{i,n}^{n+i}=-(n+2i)\alpha_{0,1}^{1}+\frac {\alpha
(i)}{(n-i)}.\tag 3.35$$
Our next goal is to determine the parameters $\alpha (i)$, for all
$i\ge 2$. From the conclusion of Lemma 3.14 we know that
 $\alpha (1)=0$. Set $n=1$ in (3.35). This yields
$\alpha_{i,1}^{i+1}=-(1+2i)\alpha_{0,1}^{1}+\frac {\alpha (i)}{(1-i)}$.
Using the result of Lemma 3.14, this leads to the equations
$-(1+2i)\alpha_{0,1}^{1}+\frac {\alpha
(i)}{(1-i)}=(i+2)\alpha_{0,1}^{1}$,
for all $i\ge 2$, from where, solving for $\alpha (i)$, we obtain
$\alpha (i)=3(1-i^2)\alpha_{0,1}^{1}$, for all $i\ge 2$.
Therefore,
$\alpha_{i,n}^{n+i}=-(n+2i)\alpha_{0,1}^{1}+\frac
{3(1-i^2)}{(n-i)}\alpha_{0,1}^{1}$, for all $n\ge 1$, and all $i\ge 1$, and
$n\ne i$.
\qed\enddemo
\par Quite similarly as in the case of the Witt algebra the following
observation helps to complete the argument.
\proclaim{Lemma 3.17} One has $\alpha_{0,1}^{1}=0$.
\endproclaim
\demo{Proof} Set $m=2,i=1$ in (3.28) to obtain
$$(n-2)\alpha_{1,n+1}^{n+2}=(2n-1)\alpha_{1-n,n+1}^2+
(n-3)\alpha_{1,n-1}^n.\tag 3.36$$
We investigate (3.36) for small values of $n$. For $n=0$ it yields
an identity. For $n=1$ it gives
$-\alpha_{1,2}^3=\alpha_{0,2}^2-2\alpha_{1,0}^1=4\alpha_{0,1}^1$,
a fact we learned from Lemma 3.16. For $n=2$ it yields a
trivial identity. For $n=3$ it gives
$\alpha_{1,4}^5=0$.
But from Lemma 3.15 we have
$\alpha_{1,4}^5=-6\alpha_{0,1}^{1}$, a clear contradiction.
Therefore we conclude that $\alpha_{0,1}^{1}=0$, which is what we set out to
show.
\qed\enddemo
\par The following two corollaries are direct consequences of the last
result and the preceding lemmas.
\proclaim{Corollary 3.18} One has $\alpha_{0,n}^n=0$, for all $n\ge
1$.
\endproclaim
\proclaim{Corollary 3.19} One has $\alpha_{i,n}^{n+i}=0$ for all
$n,i\ge 1$.
\endproclaim
\par Thus, all coalgebra structure constants $\alpha_{ij}^n=0$, whenever
$n=i+j$. Therefore the general solution of (3.24) is described by
$$\aligned \alpha (e_n)&=\sum_{i,j=0}^{\infty}\alpha_{ij}^ne_i\wedge
e_j\\
&=\sum_{i,j=0,i+j\ne n}^{\infty}\alpha_{ij}^ne_i\wedge e_j\\
&=\sum_{i,j=0,i+j\ne n}^{\infty}\left[\frac
{(2n-i)}{(n-i-j)}\alpha_{i-n,j}^0+\frac
{(2n-j)}{(n-i-j)}\alpha_{i,j-n}^0\right]e_i\wedge e_j\\
&=e_n.\left(-\sum_{i,j=0,i+j\ne 0}^{\infty}\frac
{1}{(i+j)}\alpha_{ij}^0e_i\wedge
e_j\right).\endaligned\tag 3.37$$
This shows that for $\Cal G_{\infty}$ all 1-cocycles are coboundaries.\qed

\head 4. Poisson-Lie structures on $G_{\infty}$\endhead

{ In this section we study Poisson-Lie structures on the group $G_{\infty}$. It
turns out that there exists a large class of such structures, which can be
described explicitly. In the next section we prove that in fact this class
exhausts all Poisson-Lie structures on $G_{\infty}$.}
{\par We recall the definition of a Poisson-Lie structure on the group
$G_{\infty}$ given in Sec\. 2.
It is defined as a skew-symmetric map $\{\ ,\
\}\:C^{\infty}(G_{\infty})\times C^{\infty}(G_{\infty})\to
C^{\infty}(G_{\infty})$ which is multiplicative, is a derivation in both
arguments, and satisfies the Jacobi identity. The derivation property implies
that there is a bi-vector field $\omega\in \wedge^2TG_{\infty}$ given locally
by
 $\omega_x={\omega}_{ij}(x)\pd {}{x_i}\wedge\pd {}{x_j},$
where $\omega_{ij}\in C^{\infty}(G_{\infty})$ are smooth function on
$G_{\infty}$. Then for every $f,g\in C^{\infty}(G_{\infty})$ we have
$\{ f,g\}(x) =\omega _{ij}(x)\pd {f}{x_i} \pd {g}{x_j}.$
 Note that on the right hand side we have in effect a finite sum since the
functions $f$ and $g$ depend only upon finite
number of arguments. In particular, $\{\Cal X_i,\Cal
X_j\}=\omega_{ij}$, where $\Cal X_i$, $i\in\Bbb N$, are the coordinate
functions of $x\in G_{\infty}$, i.e\., $\Cal X_i(x)=x_i$. Similarly, the
1-cocycle equation (1.6) for $\omega_{ij}$ is given by
$$\omega _{ij}({xy})=\omega _{kl}(x)\pd {z_i}{x_k} \pd {z_j}{x_l} +
\omega_{kl}{(y)}\pd {z_i}{y_k} \pd {z_j}{y_l}.\tag 4.0$$
Here again the sums on the right-hand-side are finite, since for every $n\in
\zN$ we have $z_n=z_n(x_1,\ldots,x_n;y_1,\ldots,y_n)$. The same is true
for the sums in the Jacobi identity (1.1) for the functions $\omega_{ij}$.
\par Let us introduce the formal series $\Cal
X(u):=\sum_{i=1}^{\infty}\Cal X_iu^i$. Then $x(u)=\Cal
X(u)(x)=\sum_{i=1}^{\infty}x_iu^i$ (cf\. Sec\. 2). Define the formal series
${\Omega}(u,v;{\Cal X}):=\sum_{i,j=1}^{\infty}{{\omega}_{ij}u^iv^j}$.
Thus ${\Omega}(u,v;{\Cal X})$ is a generating series for the brackets
$\omega_{ij}$. Evaluating at $x\in G_{\infty}$ we have
${\Omega}(u,v;{x})=\sum_{i,j=1}^{\infty}{{\omega}_{ij}(x)u^iv^j}$.
\proclaim{Lemma 4.0} In terms of $\Omega$ the cocycle condition {\rm (4.0)}
assumes the form
$${\Omega}(u,v;{z})={\Omega}({y}(u),{y}(v);{x})+\Omega(u,v;{y}){x}^{\prime}({y}(u)){x}^{\prime}({y}(v)), \hskip 0.2in z(u)=x(y(u)). \tag 4.1$$\endproclaim
\demo{Proof} Recall that
$z(u)=x(y(u))=\sum_{i=1}^{\infty}x_i\bigl[y(u)\bigr]^i=\sum_{i=1}^{\infty}z_iu^i$, where $z_i=(xy)_i$. From the last formula we obtain that
$$\aligned\pd {z(u)}{x_k}&=\bigl[y(u)\bigr]^k \qquad
\biggl(=\sum_{i=1}^{\infty}\pd {z_i}{x_k}u^i\biggr),\quad \text{and}\\\\
\pd {z(u)}{y_k}&=\sum_{i=1}^{\infty}ix_iu^k\bigl[y(u)\bigr]^{i-1}\\
&=u^k\sum_{i=1}^{\infty}ix_i\bigl[y(u)\bigr]^{i-1}\\
&=u^kx'\bigl(y(u)\bigr) \qquad \biggl(=\sum_{i=1}^{\infty}\pd
{z_i}{y_k}u^i\biggr).\endaligned$$
Here $x'(u)$ denotes the derivative of $x(u)$ with respect to its
argument $u$.
 If we multiply both sides of equation (4.0) by $u^iv^j$ and sum over $i$ and
$j$ we obtain
$$\eqalign{\sum_{i,j=1}^{\infty}\omega_{ij}(z)u^iv^j&=\sum_{k,l=1}^{\infty}\omega_{kl}(x)\sum_{i=1}^{\infty}\pd {z_i}{x_k}u^i\sum_{j=1}^{\infty}\pd {z_j}{x_l}v^j+\sum_{k,l=1}^{\infty}\omega_{kl}(y)\sum_{i=1}^{\infty}\pd {z_i}{y_k}u^i\sum_{j=1}^{\infty}\pd {z_j}{y_l}v^j\cr\cr
&=\sum_{k,l=1}^{\infty}\omega_{kl}(x)\bigl[y(u)\bigr]^k\bigl[y(v)\bigr]^l+x'\bigl(y(u)\bigr)x'\bigl(y(v)\bigr)\sum_{k,l=1}^{\infty}\omega_{kl}(y)u^kv^l.}$$
 Now, using the definition of $\Omega$ we finally obtain that
$${\Omega}(u,v;{z})={\Omega}({y}(u),{y}(v);{x})+\Omega(u,v;{y}){x}^{\prime}({y}(u)){x}^{\prime}({y}(v)).$$
{\it Notice also that both sides of the above equation are divisible by}
$uv$.\qed\enddemo
\par Equation (4.1) has a large class of solutions $\Omega (u,v;{x})$. Namely,
we have the following theorem.
\proclaim{Theorem 4.1} For any $\varphi=\varphi(u,v)$ with
the properties\newline
{(i)} $\varphi(u,v)$ is divisible by $u$ and $v$;\newline
{(ii)} $\varphi(u,v)=-\varphi(v,u)$,
\newline
we have the following solution of {\rm(4.1)}:
$${\Omega}(u,v;{x})=\varphi(u,v){x}^{\prime}(u){x}^{\prime}(v)-\varphi({x}(u),{x}(v)). \tag 4.2$$\endproclaim
\demo{Proof} Indeed, in terms of (4.2), the left hand side of the equation
(4.1) reads
$$\aligned{\Omega}(u,v;{z})&=\varphi(u,v){z}^{\prime}(u){z}^{\prime}(v)-\varphi({z}(u),{z}(v))\\
&=\varphi(u,v){x}^{\prime}(y(u))y^{\prime}(u){x}^{\prime}(y(v))y^{\prime}(v)-\varphi({z}(u),{z}(v)).\endaligned$$
The right hand side of (4.1) gives
$$\multline{\Omega}({y}(u),{y}(v);{x})+\Omega(u,v;{y}){x}^{\prime}({y}(u)){x}^{\prime}({y}(v))=\\
+\varphi({y}(u),{y}(v)){x}^{\prime}(y(u)){x}^{\prime}(y(v))-\varphi({z}(u),{z}(v))+\\
+\varphi(u,v){x}^{\prime}(y(u))y^{\prime}(u){x}^{\prime}(y(v))y^{\prime}(v)-\varphi({y}(u),{y}(v)){x}^{\prime}(y(u)){x}^{\prime}(y(v)).\endmultline$$
Compairing both sides we obtain an identity.
\par The condition (ii) is equivalent to ${\Omega}(u,v;{\Cal
X})=-{\Omega}(v,u;{\Cal X})$ which on the other hand is equivalent to the
skew-symmetry of the $\omega_{ij}$'s.
\par The condition (i) is needed since as noticed above $\Omega(u,v;x)$ is
divisible by $uv$. This requires that the r.h.s\. of (4.2) is divisible by
$uv$.
 From the definition of $x(u)$ it is clear that $x'(u)x'(v)$ is not divisible
by $uv$. It begins with a term  $x_1^2+2x_1x_2(u+v)+\ldots$. Suppose that
$\varphi(u,v)$ is not divisible by $uv$. Then $\varphi(x(u),x(v))$ is also not
divisible by $uv$, and so is the difference
$\varphi(u,v){z}^{\prime}(u){z}^{\prime}(v)-\varphi({z}(u),{z}(v))$,
as an easy analysis shows.\qed\enddemo
\par Next, we would like to find out for which classes of $\varphi$'s the
Jacobi identity
is satisfied. This will be an important step in the solution of the problem of
classifying all possible Lie-Poisson structures
on $G_{\infty}$. For this we use the following technical tool.
\par Let $\Cal U=\{u_1,u_2,\ldots\}$ be a countably infinite set of
indeterminates. Consider the ring of formal power series
$C^{\infty}(G_{\infty})\big[\big[\Cal U\big]\big]$ in $\Cal U$ over
the algebra $C^{\infty}(G_{\infty})$ defined as the inductive limit of
the
rings
$\left\{C^{\infty}(G_{\infty})\left[\left[u_1,\ldots,u_n\right]\right]\right\}_{n\in
\Bbb N}$. Then the map $\{\ ,\ \}\:C^{\infty}(G_{\infty})\times
C^{\infty}(G_{\infty})\to C^{\infty}(G_{\infty})$ induces a map $\{\
,\ \}\:C^{\infty}(G_{\infty})\big[\big[\Cal U\big]\big]\times
C^{\infty}(G_{\infty})\big[\big[\Cal U\big]\big]\to
C^{\infty}(G_{\infty})\big[\big[\Cal U\big]\big]$. In particular we have
$$\{ {\Cal X}(u),{\Cal X}(v)\}=\sum_{i,j=1}^{\infty}\{\Cal X_i,\Cal
X_j\}u^iv^j={\Omega}(u,v;{\Cal X}),$$
where $u=u_i$ and $v=u_j$ for some $u_i,u_j\in\Cal U$.
Then the Jacobi identities (1.1), in terms of generating series, can
be put together in a single equation
$$\{ {\Cal X}(w),\{ {\Cal X}(u),{\Cal X}(v)\}\}+\{ {\Cal X}(u),\{ {\Cal
X}(v),{\Cal X}(w)\}\}+\{ {\Cal X}(v),\{ {\Cal X}(w),{\Cal X}(u)\}\}=0, \tag
4.3$$
for any $u,v,w\in\Cal U$.
On the other hand we have
$$\aligned\{ {\Cal X}(w),\{ {\Cal X}(u),{\Cal X}(v)\}\}&=\{\Cal
X(w),\varphi(u,v)\Cal X'(u)\Cal X'(v)-\varphi(\Cal X(u),\Cal X(v)\}\\
&=\varphi(u,v)\biggl[\{ {\Cal X}(w),{\Cal X}^{\prime}(u)\}{\Cal
X}^{\prime}(v)+\{ {\Cal X}(w),{\Cal X}^{\prime}(v)\}{\Cal
X}^{\prime}(u)\biggr]-\\
&\qquad-{\partial}_1\varphi({\Cal X}(u) ,{\Cal X}(v))\{{\Cal X}(w),{\Cal
X}(u)\}-{\partial}_2\varphi({\Cal X}(u),{\Cal X}(v))\{ {\Cal X}(w),{\Cal
X}(v)\},\endaligned$$
 where $\partial_1$ denotes the derivative with respect to the first argument
and $\partial_2$ is the derivative with respect to the second argument. Also
$$\aligned\{ {\Cal X}(w),{\Cal X}^{\prime}(u)\}&={\partial}_u \{ {\Cal
X}(w),{\Cal X}(u)\}\\
&={\partial}_u\varphi(w,u){\Cal X}^{\prime}(w){\Cal
X}^{\prime}(u)+\varphi(w,u){\Cal X}^{\prime}(w){\Cal
X}^{\prime\prime}(u)+{\partial}_2\varphi({\Cal X}(w),{\Cal X}(u)){\Cal
X}^{\prime}(u),\endaligned$$
and we have similar formulae when considering the remaining two terms in (4.3)
with $w,u,v$ permuted.
\proclaim{Lemma 4.2} The solution {\rm (4.2)} satisfies the Jacobi
identity {\rm (4.3)} iff $\varphi(u,v)$ satisfies the
following functional (partial) differential equation
$$\varphi(u,v)\bigl[
{\partial}_u\varphi(w,u)+{\partial}_v\varphi(w,v)\bigr]+c.p.=0. \tag
4.4$$\endproclaim
\demo{Proof} After substituting (4.2) into (4.3), using the formulae derived
above, and collecting terms we obtain
$$\multline\biggl(\varphi (u,v)\left[\partial_{u}\varphi
(w,v)+\partial_{v}\varphi (w,v)\right]+c.p.\biggr)\Cal X'(w)\Cal X'(u)\Cal
X'(v)+\\
+\biggl(\varphi (\Cal X(v),\Cal X(u))[\partial_2\varphi(\Cal X(w),\Cal
X(v))+\partial_2\varphi(\Cal X(w),\Cal X(u))]+c.p.\biggr)=0.\endmultline\tag
*$$
Let us define $\Phi(w,u,v)$ by
$$\Phi(w,u,v):=\varphi(u,v)\bigl[
{\partial}_u\varphi(w,u)+{\partial}_v\varphi(w,v)\bigr]+c.p.$$
It is easily verified that $\Phi(w,u,v)$ is antisymmetric
with respect to each pair of its arguments. For example
$\Phi(w,u,v)=-\Phi(u,w,v)$. After evaluation at $x$ (*) becomes
$$\Phi(x(w),x(u),x(v))=x'(w)x'(u)x'(v)\Phi(w,u,v).$$
This equation is satisfied for every $x(u)$. In particular, it is true for
$x(u)=\lambda u$, where
 $\lambda\ne 0$. In this case the above equation is equivalent to
$\Phi(\lambda w,\lambda u,\lambda v)=\lambda^3\Phi(w,u,v).$
In other words $\Phi$ is homogeneous of degree 3.
But the only homogeneous function $\Phi(w,u,v)$ of degree 3 which is also
antisymmetric with respect to each pair of its arguments is $\Phi=0$. Therefore
the statement of the Lemma follows.\qed\enddemo
\proclaim{Theorem 4.3} The map $x \mapsto x^{-1}$ is an anti-Poisson
map.\endproclaim
\demo{Proof} Let $\overline{\Cal X}( u)$ denote the inverse of $\Cal X(u)$.
Then we have the identities
$$\overline{\Cal X}(\Cal X(u)) =u, \hskip 0.5in \text{and} \hskip
0.5in \Cal X(\overline{\Cal X}(u))=u ,$$
as well as (following from them)
$${\overline{\Cal X}}{}^{\prime}(\Cal X(u)){\Cal X}^{\prime}(u)=1,
\hskip 0.5in \text{and} \hskip 0.5in {\Cal X}^{\prime}(\overline{\Cal
X}(u)){\overline{\Cal X}}{}^{\prime}(u)=1.$$
On the other hand we have
$$\aligned 0&=\{ u,\Cal X( v)\}\\
&=\{ \overline{\Cal X}(\Cal X( u)) ,\Cal X(v )\}\\
&=\{ \overline{\Cal X}( w) ,\Cal X( v)\} {|}_{w=\Cal X( u)} +{\overline{\Cal X}
}{}^{\prime}( w) {|}_{w=\Cal X( u) }\{\Cal X( u) ,\Cal X( v)\}. \endaligned$$
Therefore,
$$\{\Cal X( v) ,\overline{\Cal X}( w)\}{|}_{w=\Cal X( u)}=\overline{\Cal
X}{}^{\prime}( w) {|}_{w=\Cal X( u)}\{\Cal X( u) ,\Cal X( v)\}. \tag 4.5$$
Also, we have the following chain of identities
$$\aligned 0&=\{ v,\overline{\Cal X}( w)\}{|}_{w=\Cal X( u)}\\
&=\{\overline{\Cal X}(\Cal X( v)) ,\overline{\Cal X}( w)\}{|}_{w=\Cal X( u)}\\
&=\{\overline{\Cal X}( s) ,\overline{\Cal X}( w)\}{|}_{s=\Cal X( v) ,w=\Cal X(
u)}+\overline{\Cal X}{}^{\prime}( s) {|}_{s=\Cal X( v)}\{\Cal X( v)
,\overline{\Cal X}( w)\} {|}_{w=\Cal X( u)}.\endaligned$$
Using (4.2) and (4.5), the last identity can be rewritten as
$$\aligned 0&=\varphi(\Cal X(v) ,\Cal X(u))\overline{\Cal X}{}^{\prime}(\Cal
X(v))\overline{\Cal X}{}^{\prime}(\Cal X(u))-\varphi(v,u)\\
&\qquad +\overline{\Cal X}{}^{\prime}(\Cal X(v))\overline{\Cal
X}{}^{\prime}(\Cal X(u))\left[\varphi(u,v){\Cal X}^{\prime}(u){\Cal
X}^{\prime}(v) -\varphi(\Cal X(u) ,\Cal X(v))\right]\\
&=\{\overline{\Cal X}(s) ,\overline{\Cal X}(w)\} {|}_{s=\Cal X(v) ,w=\Cal X(u)
} +\varphi(u,v)-\overline{\Cal X}{}^{\prime}(w)\overline{\Cal
X}{}^{\prime}(s)\varphi(w,s) .\endaligned$$
Thus,
$$\{\overline{\Cal X}(w) ,\overline{\Cal X}(s)\} =-\left[\overline{\Cal
X}{}^{\prime}(w)\overline{\Cal X}{}^{\prime}(s)\varphi(w,s)
-\varphi(\overline{\Cal X}(w) ,\overline{\Cal X}(s))\right],$$
and the proof is finished.\qed\enddemo
\proclaim{Theorem 4.4} For every $d\in {\zZ}$, the function
$\varphi(u,v)=uv(u^d-v^d)$ solves {\rm (4.4)}. In particular, for each
$d\in {\zZ_{+}}$, it gives rise to a family of Lie-Poisson structures on the
group $G_{\infty}$.\endproclaim
\demo{Proof} Substituting $\varphi(u,v)=uv(u^d-v^d)$ into the equation (4.4) we
obtain
$$[u^{d+1}v-uv^{d+1}][w^{d+1}-(d+1)wu^d+w^{d+1}-(d+1)wv^d]+c.p.\overset
?\to{=}0.$$
But this is an identity, since
$$2u^{d+1}w^{d+1}v-2uv^{d+1}w^{d+1}+c.p.=0,$$
$$(d+1)uwv^{2d+1}-(d+1)u^{2d+1}vw+c.p.=0,$$
and
$$(d+1)u^{d+1}v^{d+1}w-(d+1)wu^{d+1}v^{d+1}+c.p.=0.\qed$$
\enddemo
\par Writing formula (4.2) in components, with $\varphi(u,v)=uv(u^d-v^d)$, we
obtain
$${\omega}_{ij}(x)=(i-d)jx_jx_{i-d}-i(j-d)x_ix_{j-d}+x_i\sum_{\sum_{k=1}^{d+1}s_k=j}^{}x_{s_1}\ldots x_{s_{d+1}}-x_j\sum_{\sum_{k=1}^{d+1}s_k=i}^{}x_{s_1}\ldots x_{s_{d+1}}.$$
\par These formulae describe a countable family of Poisson-Lie structures on
$G_{\infty}$, thus answering the question of
existence of such.
\par In order to classify {\it all} Poisson-Lie structures on $G_{\infty}$ one
has, in particular, to classify all solutions of the functional
(partial) differential equation (4.4). The main result of this section
is formulated below.
\proclaim{Theorem 4.5} For each $d\in \zN$, and any $f_d(u) ,g_d(u)$ such that
$f_d^{\prime}(u) g_d(u) -f_d(u) g_d^{\prime}(u) =-d\lambda_{1,d+1} f_d(u)$,
where
$\lambda_{1,d+1}\ne 0$ is an arbitrary parameter, and $f_d$ has a zero
of order $d+1$ at $u=0$, there is a solution of {\rm (4.4)} given
by the series
$$\varphi_d(u,v) =\frac {1}{\lambda_{1,d+1}}\biggl[ f_d(u) g_d(v)
-f_d(v) g_d(u)\biggr].$$
The set of all solutions of {\rm (4.4)} is described
in this way.
\endproclaim\par
First, we prove a helpful lemma.
\proclaim{Lemma 4.6} For any two functions $f(u)$, $g(u)$ satisfying the
relation
$f^{\prime}(u)g(u)-f(u)g^{\prime}(u)=\alpha f(u)+\beta g(u)$, where $\alpha$,
$\beta$ are arbitrary
constants, {\rm (4.4)} has a solution in the form $\varphi
(u,v)=f(u)g(v)-f(v)g(u)$.
\endproclaim
\demo{Proof} After substituting $\varphi (u,v)=f(u)g(v)-f(v)g(u)$ into (4.4)
and collecting terms we obtain
$$\left[f(u)g^{\prime}(u)-f^{\prime}(u)g(u)\right]
f(w)g(v)-\left[f(u)g^{\prime}(u)-f^{\prime}(u)g(u)\right] f(v)g(w)+c.p.=0.$$
Using the relation
$$f^{\prime}(u)g(u)-f(u)g^{\prime}(u)=\alpha f(u)+\beta g(u)$$
the above equality transforms to
$$-\bigl[\alpha f(u)+\beta g(u)\bigr]f(w)g(v)+\bigl[\alpha f(u)+\beta
g(u)\bigr]f(v)g(w)+c.p.=0.$$
The last equality is equivalent to
$$\alpha\bigl[f(u)f(v)g(w)-f(u)f(w)g(v)+c.p.\bigr]+\beta\bigl[f(v)g(u)g(w)-f(w)g(u)g(v)+c.p.\bigr]=0.$$
But the expressions in the square brackets are identically zero as one can
easily check. Thus we obtain an identity.\qed\enddemo
\par In fact, we shall show that all solutions of the functional differential
equation (4.4) with the additional assumption that $\varphi(u,v)$ is divisible
by $uv$ are of the above form, with $\beta =0$.
\par We will seek the general solution of (4.4) as a formal power series
$\varphi(u,v) =\sum_{n,m=1}^{\infty}\lambda_{nm}u^nv^m$. Here, the antisymmetry
of $\varphi(u,v)$ implies the antisymmetry of $\lambda_{nm}$, namely
$\lambda_{nm}=-\lambda_{mn}$.
\par Substituting into (4.4) we obtain
$$\multline\sum\Sb k,n,r\endSb \sum_{s}
s[\lambda_{k-s+1,n}\lambda_{rs}+\lambda_{k,n-s+1}\lambda_{rs}+\lambda_{n-s+1,r}\lambda_{ks}+\lambda_{n,r-s+1}\lambda_{ks}+\lambda_{r-s+1,k}\lambda_{ns}+\\
+\lambda_{r,k-s+1}\lambda_{ns}]u^k v^n w^r=0,\endmultline$$
 or
$$\multline\sum\Sb s=1\endSb\Sp max(k,n,r)\endSp
s[(\lambda_{k-s+1,n}+\lambda_{k,n-s+1})\lambda_{rs}+(\lambda_{n-s+1,r}+\lambda_{n,r-s+1})\lambda_{ks}+(\lambda_{r-s+1,k}+\\
+\lambda_{r,k-s+1})\lambda_{ns}]=0.\endmultline\tag 4.6$$
 We may assume $k<n<r$, since if at least two of the indices $k,n,r$ are equal
then (4.6) is identically satisfied. Then $\max(k,n,r)=r$. Let $k=1$ and $n<r$,
then we have
$$\sum_{s=1}^{r}
s\left[\lambda_{2-s,n}\lambda_{rs}+\lambda_{1,n-s+1}\lambda_{rs}+\lambda_{n-s+1,r}\lambda_{1s}+\lambda_{n,r-s+1}\lambda_{1s}+\lambda_{r-s+1,1}\lambda_{ns}+\lambda_{r,2-s}\lambda_{ns}\right]=0,$$
which is equivalent to
$$\multline\lambda_{1n}\lambda_{r1}+\sum_{s=1}^{n}
s\lambda_{1,n-s+1}\lambda_{rs} +\sum_{s=1}^{n} s\lambda_{n-s+1,r}\lambda_{1s}
+\sum_{s=1}^{r} s\lambda_{n,r-s+1}\lambda_{1s} +\\
+\sum_{s=1}^{r}
s\lambda_{r-s+1,1}\lambda_{ns}+\lambda_{r1}\lambda_{n1}=0.\endmultline$$
  The first and the last terms in the above equation cancel each other. We make
the change of variables $s\mapsto n-s+1$, and $s\mapsto r-s+1$ in the
third and the forth terms respectively to obtain
$$\multline\sum_{s=1}^{n} s\lambda_{1,n-s+1}\lambda_{rs}
+\sum_{s=1}^{n}(n-s+1)\lambda_{sr}\lambda_{1,n-s+1}
+\sum_{s=1}^{r}(r-s+1)\lambda_{ns}\lambda_{1,r-s+1} +\\
+\sum_{s=1}^{r} s\lambda_{r-s+1,1}\lambda_{ns}=0,\endmultline$$
 which is equivalent to
$$\sum_{s=1}^{n}(n-2s+1)\lambda_{1,n-s+1}\lambda_{sr}=\sum_{s=1}^{r}(r-2s+1)\lambda_{1,r-s+1}\lambda_{sn}. \tag 4.7$$
\par A close look at the first several equations of (4.6)
$$\eqalign{&\lambda_{12}\lambda_{13}=0\cr
&\lambda_{12}(2\lambda_{14}+\lambda_{23})=0\cr
&\lambda_{13}(\lambda_{14}+\lambda_{23})=0\cr
&3\lambda_{14}\lambda_{23}-(\lambda_{23})^2-4\lambda_{13}\lambda_{24}+5\lambda_{12}\lambda_{34}=0\cr
&\lambda_{12}(3\lambda_{15}+2\lambda_{24})=0\cr
&\lambda_{13}\lambda_{15}+\lambda_{14}\lambda_{23}+\lambda_{12}\lambda_{34}=0\cr
&3\lambda_{15}\lambda_{23}-2\lambda_{23}\lambda_{24}-5\lambda_{13}\lambda_{25}+6\lambda_{12}\lambda_{35}=0\cr
&-\lambda_{14}\lambda_{15}-2\lambda_{14}\lambda_{24}+\lambda_{13}\lambda_{25}-\lambda_{12}\lambda_{35}=0\cr
&4\lambda_{15}\lambda_{24}-2(\lambda_{24})^2-5\lambda_{14}\lambda_{25}+\lambda_{23}\lambda_{25}+7\lambda_{12}\lambda_{45}=0\cr
&5\lambda_{15}\lambda_{34}-2\lambda_{24}\lambda_{34}-6\lambda_{14}\lambda_{35}+\lambda_{23}\lambda_{35}+7\lambda_{13}\lambda_{45}=0}$$
$$\vdots$$
shows that the solutions of this infinite system of quadrics fall
into three main classes. Namely those with $\lambda_{13}=0$, these with
$\lambda_{12}=0$, and these with $\lambda_{12}=0=\lambda_{13}$, the last one
being an intersection of the first two. We now proceed with describing these
classes.\newline
 (i) Let $\lambda_{12}\ne 0$ and $\lambda_{13} =0$. Then, from (4.7) with $n=2$
we obtain
$$-\lambda_{12}\lambda_{1r}
=\sum_{s=1}^{r}(r-2s+1)\lambda_{1,r-s+1}\lambda_{2s}.$$
 After multiplying both sides of the above equation by $u^r$, and summing
over $r$ we obtain
$$-\lambda_{12}\sum_{r=1}^{\infty}\lambda_{1r} u^r
=\sum_{r=1}^{\infty}\sum_{s=1}^{r}(r-2s+1)\lambda_{1,r-s+1}\lambda_{2s} u^r$$
$$\Updownarrow \hskip 1pt r\mapsto m+s-1$$
$$\eqalign{-\lambda_{12}\sum_{r=1}^{\infty}\lambda_{1r}
u^r&=\sum_{m=1}^{\infty}\sum_{s=1}^{\infty}(m-s)\lambda_{1m}\lambda_{2s}
u^{m+s-1}\cr\cr
&=\sum_{m=1}^{\infty}\sum_{s=1}^{\infty} m\lambda_{1m} u^{m-1}\lambda_{2s} u^s
-\sum_{m=1}^{\infty}\sum_{s=1}^{\infty} s\lambda_{2s} u^{s-1}\lambda_{1m}
u^m}.$$
\par Now, if we define $f_1(u)=\sum_{r=1}^{\infty}\lambda_{1r}u^r$, and
$g_1(u)=\sum_{s=1}^{\infty}\lambda_{2s}u^s$, the above equation becomes
$$-\lambda_{12}f_1(u)=f_1^{\prime}(u)g_1(u)-f_1(u)g_1^{\prime}(u).$$
 (ii) Let $\lambda_{12}=0$ and $\lambda_{13}\ne 0$. Then, from (4.7) with $n=3$
we have
$$-2\lambda_{13}\lambda_{1r}=\sum_{s=1}^{r}(r-2s+1)\lambda_{1,r-s+1}\lambda_{3s}.$$
 Define $f_2(u)$, and $g_2(u)$ to be
$f_2(u)=\sum_{r=1}^{\infty}\lambda_{1r}u^r$, and
$g_1(u)=\sum_{s=1}^{\infty}\lambda_{3s}u^s$.
Then, performing the same manipulations as in the case (i) above we obtain
$$-2\lambda_{13}f_2(u)=f_2^{\prime}(u)g_2(u)-f_2(u)g_2^{\prime}(u).$$
(iii) For $\lambda_{12}=0=\lambda_{13}$, let us assume temporarily that
$\lambda_{14}\ne 0$. Then from (4.7) with $n=4$ we have
$$-3\lambda_{14}\lambda_{1r}=\sum_{s=1}^{r}(r-2s+1)\lambda_{1,r-s+1}\lambda_{4s},$$
 and therefore
$$-3\lambda_{14}f_3(u)=f_3^{\prime}(u)g_3(u)-f_3(u)g_3^{\prime}(u).$$
 Where $f_3(u)=\sum_{s=1}^{\infty}\lambda_{1s}u^s$,
$g_3(u)=\sum_{s=1}^{\infty}\lambda_{4s}u^s$.
\par The above considerations suggest the following argument. If the first
non-zero element of the
set $\{\lambda_{1n}\}_{n\ge 2}$ is $\lambda_{1,d+1}$ $(d\ge 1)$, then from
(4.7) we deduce that
$$-d\lambda_{1,d+1}\lambda_{1r}=\sum_{s=1}^{r}(r-2s+1)\lambda_{1,r-s+1}\lambda_{d+1,s},$$
 and the above equation is equivalent to
$$-d\lambda_{1,d+1}f_d(u)=f_d^{\prime}(u)g_d(u)-f_d(u)g_d^{\prime}(u).$$
 Here $f_d$, $g_d$ are defined as $f_d(u)=\sum_{s=1}^{\infty}\lambda_{1s}u^s$,
and $g_d(u)=\sum_{s=1}^{\infty}\lambda_{d+1,s}u^s$.
(Obviously, the first $d$ terms in the definition of $f_d$ are zero.)
\par Thus, we will parametrize all classes of solutions of (4.4) by $d\in \zN$
such that
$\lambda_{12}=\ldots =\lambda_{1d} =0$, and $\lambda_{1,d+1}\ne 0$. In what
follows, we will show that
for each $d\in \zN$ with the above property $\varphi(u,v)$ is given by
$\varphi_d(u,v) =\frac {1}{\lambda_{1,d+1}}\left[ f_d(u) g_d(v) -f_d(v)
g_d(u)\right]$,
and therefore is a solution of (4.4), according to Lemma 4.6. This we will show
by proving that for each $d\in \zN$ we have
$$\lambda_{nm}=\frac{1}{\lambda_{1,d+1}}\biggl[\lambda_{1n}\lambda_{d+1,m}-\lambda_{1m}\lambda_{d+1,n}\biggr],
\qquad \text{for all}\quad {n,m\ge 1}. \tag 4.8$$
\proclaim{Lemma 4.7} For any fixed $d\in \zN$ such that $\lambda_{1n}=0$ for
$1\le n\le d$, $\lambda_{1,d+1}\ne 0$, it follows that $\lambda_{sn}=0$ for
$1\le s\le {d-1}$, $1\le n\le d$.\endproclaim
\demo{Proof} Since $1\le n\le d$, it follows from (4.7) that (the l.h.s\. is
zero)
$$\sum_{s=1}^{r-d}(r-2s+1)\lambda_{1,r-s+1}\lambda_{sn}=0. \tag 4.9$$
Since $n<r$, if\newline
(i) $r=d+1$, then  $d\lambda_{1,d+1}\lambda_{1n}=0$, which is an
identity;\newline
(ii) $r=d+2$, then
$(d+1)\lambda_{1,d+2}\lambda_{1n}+(d-1)\lambda_{1,d+1}\lambda_{2n}=0 \implies
(d-1)\lambda_{1,d+1}\lambda_{2n}=0 \implies \lambda_{2n}=0$.\newline
Assume now that $\lambda_{sn}=0$ for $1\le s\le m<d-1$, $1\le n\le d$. We would
like to show that this implies $\lambda_{m+1,n}=0$. But from (4.9) with
$r=d+m+1$ it follows that
$(d-m)\lambda_{1,d+1}\lambda_{m+1,n}=0 \implies \lambda_{m+1,n}=0.$
 Therefore, letting $r$ run in the interval $d+1\le r\le 2d-1$ finishes the
proof.\qed\enddemo
 \remark{Remark} For $r=2d$ and $r=2d+1$ we obtain identities. For
$r>2d+1$ we obtain relations which are particular cases of (4.8). For example,
when $r=2d+2$ we have
$$\lambda_{1,d+2}\lambda_{d+1,n}-\lambda_{1,d+1}\lambda_{d+2,n}=0 \implies
\lambda_{n,d+2}=\frac {\lambda_{1,d+2}}{\lambda_{1,d+1}} \lambda_{n,d+1}
.$$\endremark
\proclaim{Lemma 4.8} For each fixed $d\in \zN$, such that $\lambda_{1m}=0$ for
$1\le m\le d$, $\lambda_{1,d+1}\ne 0$, it follows that
$\lambda_{d+1,n}$ is a rational function of $\lambda_{1,d+1},\ldots
,\lambda_{1,d+n}$, $\forall\ {n\ge 1}$.\endproclaim
\demo{Proof} Let us consider (4.7) with $r=n+d$
$$\sum_{s=1}^{d+1}(d-2s+1)\lambda_{1,d+2-s}\lambda_{s,n+d}=\sum_{s=1}^{n}(n+d-2s+1)\lambda_{1,n+d-s+1}\lambda_{s,d+1} \hskip 10pt (n\ge 1), \tag 4.10$$
which is equivalent to
$$d\lambda_{1,d+1}\lambda_{1,n+d}=\sum_{s=1}^{n-1}(n+d-2s+1)\lambda_{1,n+d-s+1}\lambda_{s,d+1}+(d-n+1)\lambda_{1,d+1}\lambda_{n,d+1}$$
since only the first term on the l.h.s. is non-zero. Therefore, solving for
$\lambda_{d+1,n}$ we obtain
$$\lambda_{d+1,n}=-\frac
{1}{(d-n+1)\lambda_{1,d+1}}\left[d\lambda_{1,d+1}\lambda_{1,n+d}-\sum_{s=1}^{n-1}(n+d-2s+1)\lambda_{1,n+d-s+1}\lambda_{s,d+1}\right]. \tag 4.11$$
\par The above formula gives a recursive relation for $\lambda_{n,d+1}$'s,
whenever $n\ne d+1$. To finish the proof of
the lemma we write the first three relations. For $n=1$, (4.10) is an identity.
For $n=2$ we have
$$d\lambda_{1,d+1}\lambda_{1,2+d}=(d+1)\lambda_{1,2+d}\lambda_{1,d+1}+(d-1)\lambda_{1,d+1}\lambda_{2,d+1},$$
from which it follows that
$$\lambda_{2,d+1}=-\frac {1}{d-1}\lambda_{1,d+2}. \tag 4.12$$
For $n=3$ we have
$$\eqalign{\lambda_{d+1,3}&=-\frac
{1}{(d-2)\lambda_{1,d+1}}\biggl[d\lambda_{1,d+1}\lambda_{1,d+3}-(d+2)\lambda_{1,d+3}\lambda_{1,d+1}-d\lambda_{1,d+2}\lambda_{1,d+1}\biggr]\cr
&=-\frac
{1}{(d-2)\lambda_{1,d+1}}\left[-2\lambda_{1,d+1}\lambda_{1,d+3}+{d\over {d-1}}
(\lambda_{1,d+2})^2\right],}$$
 where we have used (4.11). Finally for $n=4$ we have
$$\multline\lambda_{d+1,4}=-\frac
{1}{(d-3)\lambda_{1,d+1}}\biggl[d\lambda_{1,d+1}\lambda_{1,d+4}-(d+3)\lambda_{1,d+4}\lambda_{1,d+1}-(d+1)\lambda_{1,d+3}\lambda_{2,d+1}-\\
\shoveright{-(d-1)\lambda_{1,d+2}\lambda_{3,d+1}\biggr]}\\
=-\frac {1}{(d-3)\lambda_{1,d+1}}\biggl[-3\lambda_{1,d+1}\lambda_{1,d+4}+\frac
{(3d-5)d}{(d-1)(d-2)}\lambda_{1,d+2}\lambda_{1,d+3}-\frac {d}{d-2}\frac
{(\lambda_{1,d+2})^3}{\lambda_{1,d+1}}\biggr].\endmultline$$
Therefore, repeating this process $n$ times we obtain that $\lambda_{d+1,n}$ is
a rational function of the type
as stated, since all $\lambda_{d+1,s}$ with $1\le s\le n-1$ are rational
functions. This completes the proof of the lemma.\qed\enddemo
\remark{Remark} Lemma 4.8 shows that for each $d\in \zN$ the series $f_d$ and
$g_d$ are defined only in terms of $\{\lambda_{1n}\}_{n\ge d+1}$.
That is, for each $d\in \zN$ we will have a solution of (4.4) determined by
this infinite set of parameters. It turns out
that these parameters are not completely independent. Namely, we have the
following lemma.\endremark
\proclaim{Lemma 4.9} For each $d\in \zN$ there exists the following single
relation between $\lambda_{1n}$'s (with $d+1\le n\le 2d+1$)
$$\lambda_{1,2d+1}=-\frac
{1}{d\lambda_{1,d+1}}\sum_{s=2}^{d}2(d+1-s)\lambda_{1,2d+2-s}\lambda_{s,d+1}.$$
Here $\lambda_{s,d+1}=\lambda_{s,d+1}(\lambda_{1,d+1},\ldots,\lambda_{1,d+s})$,
$2\le s\le d$,  are rational functions
according to the previous lemma.\endproclaim
\demo{Proof} From (4.11) we have
$$-(d-n+1)\lambda_{1,d+1}\lambda_{d+1,n}=d\lambda_{1,d+1}\lambda_{1,n+d}-\sum_{s=1}^{n-1}(n+d-2s+1)\lambda_{1,n+d-s+1}\lambda_{s,d+1}.$$
If $n=d+1$ the l.h.s\. is zero, and we obtain
$$\eqalign{d\lambda_{1,d+1}\lambda_{1,2d+1}&=\sum_{s=1}^{d}(2d+2-2s)\lambda_{1,2d+2-s}\lambda_{s,d+1}\cr
&=2d\lambda_{1,2d+1}\lambda_{1,d+1}+\sum_{s=2}^{d}2(d+1-s)\lambda_{1,2d+2-s}\lambda_{s,d+1}.}$$
 From this, the statement follows. This means that for each $d\in\zN$ we can
solve for $\lambda_{2d+1}$ in terms of $\lambda_{1,d+1},\ldots,\lambda_{1,2d}$,
and $\lambda_{1,2d+1}$ is a rational function of these variables. This is the
only relation between $\lambda_{1n}$'s. This is easily seen from (4.11).
Multiplying both sides of (4.10) by $(d+1-n)$ we see that the l.h.s\. of the
equality  so obtained vanishes if and only if $n=d+1$. From this we obtain
exactly one relation between
$\lambda_{1n}$'s for $d+1\le n\le 2d+1$. Thus for $d=1$ we obtain
$\lambda_{13}=0$, for $d=2$ we have $\lambda_{15}=\frac
{(\lambda_{14})^2}{\lambda_{13}}$, for $d=3$ we have
$\lambda_{17}=2\frac {\lambda_{15}\lambda_{16}}{\lambda_{14}}-\frac
{(\lambda_{15})^3}{(\lambda_{14})^2},$
and so on.\qed\enddemo
\par Summarizing, each natural number $d\in\zN$ labels a branch in the space of
solutions of (4.6) for which the set of parameters $\{\lambda_{1n}\}_{n\ge
d+1,n\ne 2d+1}$ forms a basis. Here $\lambda_{1,2d+1}$ is a rational function
of $\lambda_{1,d+1},\ldots,\lambda_{1,2d}$.
\proclaim{Lemma 4.10} For each fixed $d\in \zN$ and every $n,m\ge 1$ the
following formula is valid
$$\lambda_{mn}=\frac{1}{\lambda_{1,d+1}}\biggl[\lambda_{1m}\lambda_{d+1,n}-\lambda_{1n}\lambda_{d+1,m}\biggr].\tag 4.13$$\endproclaim
\demo{Proof} The plan of the proof is as follows\newline
(i) First we prove (4.13) for $1\le m\le d$ and $n\ge d+1$. For $n<d+1$ there
is nothing to prove, because all $\lambda_{mn}$'s are zero according to Lemma
4.7;\newline
(ii) Second, we prove (4.13) for $m\le 2d$ and $n>m$, using (i). Here, we prove
it
first for $\lambda_{d+2,d+3}$ (for $\lambda_{d+1,d+2}$ the statement is
trivial). Then,
we prove it for $\lambda_{d+2,n}$, and every $n\ge d+3$, using an inductive
argument. Next,
we prove that if the statement is true for $\lambda_{m-1,m}$, then it is true
for $\lambda_{m,m+1}$,
for some $m\le 2d$. Last, we fix $m\le 2d$, and use again an inductive argument
to prove it for $\lambda_{mn}$, and every $n>m$;\newline
(iii) Our third step is to apply induction to the argument $m$. Namely,
assuming that
the statement is true for $\lambda_{mn}$, where $m\le {(k-1)d}$, and every
$n>m$, we prove it for
$m\le kd$, and every $n>m$. \newline
\par The proof is technical, but not difficult, and uses extensively formula
(4.11), which
we now write as
$$\sum_{s=1}^{n-1}(n+d-2s+1)\lambda_{1,n+d-s+1}\lambda_{s,d+1}=(d-n+1)\lambda_{1,d+1}\lambda_{d+1,n}+d\lambda_{1,d+1}\lambda_{1,n+d}.\tag 4.14$$
 (i) Let $1\le m\le d$, and $N\ge d+1$. Fix $m$, and assume that the statement
is true for all $d+1\le n\le N-1$. We want to prove it for $\lambda_{mN}$. But
from (4.7) we have
$$\sum_{s=1}^{m-1}(m-2s+1)\lambda_{1,m-s+1}\lambda_{s,N+d}=-\sum_{s=d+1}^{N}(N+d-2s+1)\lambda_{1,N+d-s+1}\lambda_{ms}.$$
But the l.h.s\. of this equation is zero, since $m\le d $. Therefore
$$\align
0&=\sum_{s=d+1}^{N-1}(N+d-2s+1)\lambda_{1,N+d-s+1}\lambda_{ms}+(d-N+1)\lambda_{1,d+1}\lambda_{mN}\\
&=\frac
{1}{\lambda_{1,d+1}}\sum_{s=d+1}^{N-1}(N+d-2s+1)\lambda_{1,N+d-s+1}\left[\lambda_{1m}\lambda_{d+1,s}-\lambda_{1s}\lambda_{d+1,m}\right]+\\
&\qquad{ +(d-N+1)\lambda_{1,d+1}\lambda_{mN}}\\
&=-\frac
{\lambda_{d+1,m}}{\lambda_{1,d+1}}\sum_{s=d+1}^{N-1}(N+d-2s+1)\lambda_{1,N+d-s+1}\lambda_{1s}+(d-N+1)\lambda_{1,d+1}\lambda_{mN}\\
&=(d-N+1)\lambda_{d+1,m}\lambda_{1N}+(d-N+1)\lambda_{1,d+1}\lambda_{mN},\endalign$$
and we obtain
$$\lambda_{mN}=-\frac {\lambda_{1N}}{\lambda_{d+1,m}}\lambda_{1,d+1}.$$
This is exactly (4.13) with $\lambda_{1m}=0$, which is a consequence of the
assumption on $m$.\newline
 (ii.a) The statement is true for $\lambda_{d+1,d+2}$. Let us now assume that
$m=d+2$, $n=d+3$. Then from (4.7) we have
$$\sum_{s=1}^{d+1}(d-2s+3)\lambda_{1,d-s+3}\lambda_{s,2d+3}=\sum_{s=1}^{d+1}(2d-2s+3)\lambda_{1,2d-s+3}\lambda_{s,d+2}-2\lambda_{1,d+1}\lambda_{d+3,d+2}.$$
Now, we use (i) in both sides of the above equation to obtain
$$\multline\frac
{1}{\lambda_{1,d+1}}\sum_{s=1}^{d+1}(d-2s+3)\lambda_{1,d-s+3}\left[\lambda_{1s}\lambda_{d+1,2d+3}-\lambda_{1,2d+3}\lambda_{d+1,s}\right]=\\
=\frac
{1}{\lambda_{1,d+1}}\sum_{s=1}^{d+1}(2d-2s+3)\lambda_{1,2d-s+3}\left[\lambda_{1s}\lambda_{d+1,d+2}-\lambda_{1,d+2}\lambda_{d+1,s}\right]-2\lambda_{1,d+1}\lambda_{d+3,d+2},\endmultline$$
which is equivalent to
$$\multline\frac
{\lambda_{d+1,2d+3}}{\lambda_{1,d+1}}\sum_{s=1}^{d+1}(d-2s+3)\lambda_{1,d-s+3}\lambda_{1s}-\frac {\lambda_{1,2d+3}}{\lambda_{1,d+1}}\sum_{s=1}^{d+1}(d-2s+3)\lambda_{1,d-s+3}\lambda_{d+1,s}=\\
=-2\lambda_{1,d+1}\lambda_{d+3,d+2}+\frac
{\lambda_{1,2d+3}}{\lambda_{1,d+1}}\sum_{s=1}^{d+1}(2d-2s+4)\lambda_{1,2d-s+4}\lambda_{1s}-\\
-\frac
{\lambda_{1,d+2}}{\lambda_{1,d+1}}\sum_{s=1}^{d+1}(2d-2s+4)\lambda_{1,2d-s+4}\lambda_{d+1,s}.\endmultline$$
The first term on the l.h.s\. is zero, and the second term on the r.h.s\.
contains only one summand thus giving
$$\multline-\frac
{\lambda_{1,2d+3}}{\lambda_{1,d+1}}\sum_{s=1}^{d+1}(d-2s+3)\lambda_{1,d-s+3}\lambda_{d+1,s}=\\=-2\lambda_{1,d+1}\lambda_{d+3,d+2}+2\lambda_{1,d+2}\lambda_{d+1,d+2}
-\frac
{\lambda_{1,d+2}}{\lambda_{1,d+1}}\sum_{s=1}^{d+1}(2d-2s+3)\lambda_{1,2d-s+3}\lambda_{d+1,s}.\endmultline$$
The term on the l.h.s\. has only two summands. The third term on the r.h.s\. we
transform using (4.14). Thus,
$$\multline\frac {\lambda_{1,2d+3}}{\lambda_{1,d+1}}\left[
(d+1)\lambda_{1,d+2}\lambda_{1,d+1}+(d-1)\lambda_{1,d+1}\lambda_{2,d+1}\right]=\\=-2\lambda_{1,d+1}\lambda_{d+3,d+2}+2\lambda_{1,d+3}\lambda_{d+1,d+3}
+\frac
{\lambda_{1,d+2}}{\lambda_{1,d+1}}\left[-2\lambda_{1,d+1}\lambda_{d+1,d+3}+d\lambda_{1,d+1}\lambda_{2d+3}\right].\endmultline$$
Collecting all the terms we obtain
$$\lambda_{d+2,d+3}=\frac
{1}{\lambda_{1,d+1}}\biggl[\lambda_{1,d+2}\lambda_{d+1,d+3}-\lambda_{1,d+3}\lambda_{d+1,d+3}\biggr].$$
 (ii.b) Let us assume now that the statement is true for $\lambda_{d+2,k}$,
$d+3\le k\le n-1$. We would like to prove it for $k=n$. From (4.7) we have
$$\sum_{s=1}^{d+1}(d-2s+3)\lambda_{1,d-s+3}\lambda_{s,n+d}=\sum_{s=1}^{n-1}(n+d-2s+1)\lambda_{1,n+d-s+1}\lambda_{s,d+2}+(d-n+1)\lambda_{1,d+1}\lambda_{n,d+2}.$$
The l.h.s\. has only two summands. Also, applying the inductive hypothesis we
obtain
$${(d+1)\lambda_{1,d+2}\lambda_{1,n+d}+(d-1)\lambda_{1,d+1}\lambda_{2,n+d}=}$$
$$\align &=\frac
{1}{\lambda_{1,d+1}}\sum_{s=1}^{n-1}(n+d-2s+1)\lambda_{1,n+d-s+1}\biggl[\lambda_{1s}\lambda_{d+1,d+2}-\lambda_{1,d+2}\lambda_{d+1,s}\biggr]+
\\ &\qquad{ +(d-n+1)\lambda_{1,d+1}\lambda_{n,d+2}}\\
&=\frac
{\lambda_{d+1,d+2}}{\lambda_{1,d+1}}\sum_{d+1}^{n-1}(n+d-2s+1)\lambda_{1,n+d-s+1}\lambda_{1s}-\\
&\qquad{ -\frac
{\lambda_{1,d+2}}{\lambda_{1,d+1}}\sum_{s=1}^{n-1}(n+d-2s+1)\lambda_{1,n+d-s+1}\lambda_{d+1,s}+(d-n+1)\lambda_{1,d+1}\lambda_{n,d+2}}\\
&=(n-d-1)\lambda_{d+1,d+2}\lambda_{1n}+\frac
{\lambda_{1,d+2}}{\lambda_{1,d+1}}\biggl[
(d-n+1)\lambda_{1,d+1}\lambda_{d+1,n}+d\lambda_{1,d+1}\lambda_{1,n+d}\biggr]+\\
&\qquad{ +(d-n+1)\lambda_{1,d+1}\lambda_{n,d+2}.}\endalign$$
The second term in the above expression has been obtained by using formula
(4.14). Applying the inductive argument to $\lambda_{2,n+d}$, namely
$\lambda_{2,n+d}=-\frac
{1}{\lambda_{1,d+1}}\left[-\lambda_{1,n+d}\lambda_{d+1,2}\right]$, and
collecting terms we obtain
$$\lambda_{d+2,n}=\frac
{1}{\lambda_{1,d+1}}\biggl[\lambda_{1,d+2}\lambda_{d+1,n}-\lambda_{1n}\lambda_{d+1,d+2}\biggr].$$
 (ii.c) Suppose that the statement is true for $\lambda_{m-1,n}$, $\forall\
{n\ge m-1}$. We are going to prove it for $\lambda_{m,m+1}$. From (4.7) we have
$$\sum_{s=1}^{m-1}(m-2s+1)\lambda_{1,m-s+1}\lambda_{s,m+d+1}=\sum_{s=1}^{m}(m+d-2s+2)\lambda_{1,m+d-s+2}\lambda_{sm}+(d-m)\lambda_{1,d+1}\lambda_{m+1,m}.$$
Using the induction hypothesis the above equation transforms to
$$\multline\frac
{1}{\lambda_{1,d+1}}\sum_{s=1}^{m-1}(m-2s+1)\lambda_{1,m-s+1}\left[\lambda_{1s}\lambda_{d+1,m+d+1}-\lambda_{1,m+d+1}\lambda_{d+1,s}\right]=\\=(d-m)\lambda_{1,d+1}\lambda_{m+1,m}
+\frac
{1}{\lambda_{1,d+1}}\sum_{s=1}^{m-1}(m+d-2s+2)\lambda_{1,m+d-s+2}\left[\lambda_{1s}\lambda_{d+1,m}-\lambda_{1m}\lambda_{d+1,s}\right].\endmultline$$
Expanding, we have
$$\multline\frac
{\lambda_{d+1,m+d+1}}{\lambda_{1,d+1}}\sum_{s=d+1}^{m-d}(m-2s+1)\lambda_{1,m-s+1}\lambda_{1s}-\frac {\lambda_{1,m+d+1}}{\lambda_{1,d+1}}\sum_{s=1}^{m-d}(m-2s+1)\lambda_{1,m-s+1}\lambda_{d+1,s}=\\=(d-m)\lambda_{1,d+1}\lambda_{m+1,m}
+\frac
{\lambda_{d+1,m}}{\lambda_{1,d+1}}\sum_{s=d+1}^{m}(m+d-2s+2)\lambda_{1,m+d-s+2}\lambda_{1s}-\\
-\frac
{\lambda_{1m}}{\lambda_{1,d+1}}\sum_{s=1}^{(m+1)-1}(m+d-2s+2)\lambda_{1,m+d-s+2}\lambda_{d+1,s}.\endmultline$$
The first term on the l.h.s\. is zero. From the others we obtain
$$\multline-(2d-m+1)\lambda_{1,m+d+1}\lambda_{d+1,m-d}+\frac
{\lambda_{1,m+d+1}}{\lambda_{1,d+1}}\biggl[
(2d-m+1)\lambda_{1,d+1}\lambda_{d+1,m-d}+d\lambda_{1,d+1}\lambda_{1m}\biggr]=\\
\shoveleft{=(d-m)\lambda_{1,d+1}\lambda_{m+1,m}+(m-d)\lambda_{1,m+1}\lambda_{d+1,m}+}\\
+\frac {\lambda_{1m}}{\lambda_{1,d+1}}\biggl[
(d-m)\lambda_{1,d+1}\lambda_{d+1,m+1}+d\lambda_{1,d+1}\lambda_{1,m+d+1}\biggr].\endmultline$$
Collecting terms we arrive at
$$\lambda_{m,m+1}=\frac
{1}{\lambda_{1,d+1}}\biggl[\lambda_{1,m}\lambda_{d+1,m+1}-\lambda_{1,m+1}\lambda_{d+1,m}\biggr].$$
 (ii.d) Let us assume now, that the statement is true for $\lambda_{mk}$, where
$m\le 2d$, and all $k\ge m+1$ up to $k=n-1$. We will prove it for $k=n$.
Again, from (4.7) we have
$$\sum_{s=1}^{m-1}(m-2s+1)\lambda_{1,m-s+1}\lambda_{s,n+d}=\sum_{s=1}^{n+1}(n+d-2s+1)\lambda_{1,n+d-s+1}\lambda_{sm}+(d-n+1)\lambda_{1,d+1}\lambda_{nm}.$$
Next, we apply the induction hypothesis to obtain
$$\multline\frac
{1}{\lambda_{1,d+1}}\sum_{s=1}^{m-1}(m-2s+1)\lambda_{1,m-s+1}\left[\lambda_{1s}\lambda_{d+1,n+d}-\lambda_{1,n+d}\lambda_{d+1,s}\right]=\\
=\frac
{1}{\lambda_{1,d+1}}\sum_{s=1}^{n-1}(n+d-2s+1)\lambda_{1,n+d-s+1}\left[\lambda_{1s}\lambda_{d+1,m}-\lambda_{1m}\lambda_{d+1,s}\right]+(d-n+1)\lambda_{1,d+1}\lambda_{nm},\endmultline$$
which leads to
$$\multline\frac
{\lambda_{d+1,n+d}}{\lambda_{1,d+1}}\sum_{s=d+1}^{m-1}(m-2s+1)\lambda_{1,m-s+1}\lambda_{1s}-\frac {\lambda_{1,n+d}}{\lambda_{1,d+1}}\sum_{s=1}^{m-1}(m-2s+1)\lambda_{1,m-s+1}\lambda_{d+1,s}=\\
=(d-n+1)\lambda_{1,d+1}\lambda_{nm}+\frac
{\lambda_{d+1,m}}{\lambda_{1,d+1}}\sum_{s=d+1}^{n-1}(n+d-2s+1)\lambda_{1,n+d-s+1}\lambda_{1s}-\\
-\frac
{\lambda_{1m}}{\lambda_{1,d+1}}\sum_{s=1}^{n-1}(n+d-2s+1)\lambda_{1,n+d-s+1}\lambda_{d+1,s}.\endmultline\tag *$$
The first term on the l.h.s\. is zero. The second one we represent as (notice
that the sum is up to $m-d$)
$$-(2d-m+1)\lambda_{1,n+d}\lambda_{d+1,n-d}+\sum_{s=1}^{m-(d+1)}(m-2s+1)\lambda_{1,m-s+1}\lambda_{s,d+1}.$$
For the second term in the expression above we apply (4.14), namely
$$\sum_{s=1}^{m-(d+1)}(m-2s+1)\lambda_{1,m-s+1}\lambda_{s,d+1}=(2d-m+1)\lambda_{1,d+1}\lambda_{d+1,m-d}+d\lambda_{1,d+1}\lambda_{1m}.$$
We do the same for the last term on the r.h.s\. of (*) too. Thus finally, we
come to the equality
$$\multline-(2d-m+1)\lambda_{1,n+d}\lambda_{d+1,m-d}+(2d-m+1)\lambda_{1,n+d}\lambda_{d+1,m-d}+d\lambda_{1,n+d}\lambda_{1m}=\\
=(d-n+1)\lambda_{1,d+1}\lambda_{nm}+(n-d-1)\lambda_{d+1,m}\lambda_{1n}+(d-n+1)\lambda_{d+1,n}\lambda_{1m}+d\lambda_{1,n+d}\lambda_{1m},\endmultline$$
and after cancellation we get
$$\lambda_{mn}=\frac
{1}{\lambda_{1,d+1}}\biggl[\lambda_{1m}\lambda_{d+1,n}-\lambda_{1n}\lambda_{d+1,m}\biggr].$$
 (iii) In this case the arguments repeat vis-\`a-vis the arguments presented in
case (ii). Thus, we omit them.
\par This concludes the proof Lemma 4.10 and of Theorem 4.5.\qed\enddemo
\par Using the above results we are going to describe now a one-parameter
extension of the solution $\varphi (u,v) =uv(u^d-v^d)$ obtained in Theorem 4.5,
which
is of course a particular case of the general infinite-parameter solution.
\proclaim{Lemma 4.11} For each $d\in \zN\setminus \{1\}$, if
$\lambda_{1n}=\frac {(\lambda_{1,d+2})^{n-d-1}}{(\lambda_{1,d+1})^{n-d-2}}$,
for every ${n\ge d+1}$, it follows that
$\lambda_{n,d+1}=0$, for every ${n\ge d+1}$, and
$\lambda_{n,d+1}=-\frac {1}{d-1}\frac
{(\lambda_{1,d+2})^{n-1}}{(\lambda_{1,d+1})^{n-2}}$, for all $n$ such
that $2\le n\le d$.\endproclaim
\demo{Proof} From (4.10) we obtain
$$\lambda_{n,d+1}=\frac {1}{(d-n+1)\lambda_{1,d+1}}\left[d\frac
{(\lambda_{1,d+1})^{n-1}}{(\lambda_{1,d+1})^{n-3}}-\sum_{s=1}^{n-1}(n+d-2s+1)\frac {(\lambda_{1,d+2})^{n-s}}{(\lambda_{1,d+1})^{n-s-1}}\lambda_{s,d+1}\right]. \tag 4.15$$
If $n=2$ we have
$$\lambda_{2,d+1}=\frac
{1}{(d-n+1)\lambda_{1,d+1}}\biggl[d\lambda_{1,d+1}\lambda_{1,d+2}-(d+1)\lambda_{1,d+2}\lambda_{1,d+1}\biggr]=-\frac {1}{d-1}\lambda_{1,d+2}.$$
Let as assume now that $\lambda_{k,d+1}=-\frac {1}{d-1}\frac
{(\lambda_{1,d+2})^{k-1}}{(\lambda_{1,d+1})^{k-2}}$ for all $k$, such that
$2\le k\le n-1 <d$. We will now prove that this relation is true for $k=n\le
d$. From (4.14) we obtain
$$\allowdisplaybreaks\align\lambda_{n,d+1}&=\frac
{1}{(d-n+1)\lambda_{1,d+1}}\left[d\frac
{(\lambda_{1,d+2})^{n-1}}{(\lambda_{1,d+1})^{n-3}}+\frac
{1}{d-1}\sum_{s=1}^{n-1}(n+d-2s+1)\frac
{(\lambda_{1,d+2})^{n-s}}{(\lambda_{1,d+1})^{n-s-1}}\lambda_{s,d+1}\right]\\
&=\frac {1}{(d-n+1)\lambda_{1,d+1}}\left[d\frac
{(\lambda_{1,d+2})^{n-1}}{(\lambda_{1,d+1})^{n-3}}+\frac {1}{d-1}\frac
{(\lambda_{1,d+2})^{n-1}}{(\lambda_{1,d+1})^{n-3}}\sum_{s=2}^{n-1}(n+d-2s+1)\right]-\\
&\qquad{ -\frac {(\lambda_{1,d+2})^{n-1}}{(\lambda_{1,d+1})^{n-2}}} \\
&=\frac {1}{(d-n+1)\lambda_{1,d+1}}\biggl[-(n-1)+\frac {1}{d-1}(n-2)(n+d+1)-\\
&\qquad\qquad\qquad\qquad\qquad\qquad\qquad\qquad\qquad\qquad{ -\frac
{2}{d-1}\biggl(\frac {n(n-1)}{2}-1\biggr)\biggr]\frac
{(\lambda_{1,d+2})^{n-1}}{(\lambda_{1,d+1})^{n-3}}}\\
&=\frac {1}{(d-n+1)(d-1)}\biggl[-(n-1)(d-1)+(n-2)(n+d+1)-\\
&\qquad\qquad\qquad\qquad\qquad\qquad\qquad\qquad\qquad\qquad\qquad\qquad{
-n(n-1)+2\biggr]\frac {(\lambda_{1,d+2})^{n-1}}{(\lambda_{1,d+1})^{n-2}}}\\
&=\frac {1}{(d-n+1)(d-1)}(-d+n-1)\frac
{(\lambda_{1,d+2})^{n-1}}{(\lambda_{1,d+1})^{n-2}}\\\\
&=-\frac {1}{d-1}\frac
{(\lambda_{1,d+2})^{n-1}}{(\lambda_{1,d+1})^{n-2}}.\endalign$$
\par This proves the second part of the lemma. If $n=d+1$, then the
l.h.s\. of (4.15) is zero. We check for consistency of whether the r.h.s\. is
zero. The expression
in the brackets reads
$$-d+\frac {1}{d-1}\sum_{s=2}^{d}2(d-s+1)=\frac
{1}{d-1}\left[-d^2+d+2(d+1)(d-1)-2\biggl(\frac {d(d+1)}{2}-1\biggr)\right]=0.$$
 Now, we check whether the statement is true for $\lambda_{d+2,d+1}$.
 From (4.15) one has
$$\align\lambda_{d+2,d+1}&=-\left[-(d+1)+\frac
{1}{d-1}\sum_{s=2}^{d}(2d+3-2s)\right]\frac
{(\lambda_{1,d+2})^{d+1}}{(\lambda_{1,d+1})^d}\\
&=-\frac {1}{d-1}\left[-d^2+1+(d-1)(2d+3)-2\biggl(\frac
{d(d+1)}{2}-1\biggr)\right]\frac
{(\lambda_{1,d+2})^{d+1}}{(\lambda_{1,d+1})^d}=0.\endalign$$
Thus $\lambda_{d+2,d+1}=0$. We use again an inductive argument. We assume that
$\lambda_{k,d+1}=0$ for all $k$ such that $d+1\le k\le n-1$, and we would like
to show
that this implies $\lambda_{n,d+1}=0$. From (4.15) we have
$$\allowdisplaybreaks\align\lambda_{n,d+1}&=\frac {1}{(d-n+)}\left[ d\frac
{(\lambda_{1,d+2})^{n-1}}{(\lambda_{1,d+1})^{n-2}}-\sum_{s=1}^{d}(n+d-2s+1)\frac {(\lambda_{1,d+2})^{n-s}}{(\lambda_{1,d+1})^{n-2}}\lambda_{s,d+1}\right]\\
&=\frac {1}{(d-n+1)}\left[-(n-1)+\frac
{1}{d-1}\sum_{s=2}^{d}(n+d-2s+1)\right]\frac
{(\lambda_{1,d+2})^{n-1}}{(\lambda_{1,d+1})^{n-2}}\\
&=\frac {1}{(d-n+1)(d-1)}\biggl[-(n-1)(d-1)+(d-1)(n+d+1)-\\
&\qquad\qquad\qquad\qquad\qquad\qquad\qquad\qquad{ -2\biggl(\frac
{d(d+1)}{2}-1\biggr)\biggr]\frac
{(\lambda_{1,d+2})^{n-1}}{(\lambda_{1,d+1})^{n-2}}=0.}\endalign$$
This concludes the proof of the lemma.\qed\enddemo
\proclaim{Theorem 4.12} For every $d\ge 2$ and $\lambda := \frac
{\lambda_{1,d+2}}{\lambda_{1,d+1}}$ the function
$$\varphi_{d,\lambda }(u,v)=\frac {1}{(d-1)(1-\lambda u)(1-\lambda v)}\biggl\{
(d-1)uv(v^d-u^d)+\lambda du^2v^2(u^{d-1}-v^{d-1})\biggr\} \tag 4.16$$
is a solution of {\rm (4.4)}. It is a one-parameter extension of the
solution $\varphi (u,v)=uv(v^d-u^d)$, for every ${d\ge 2}$, which we obtain
back from {\rm (4.16)} by setting $\lambda =0$.\endproclaim
\remark{Remark} One can obtain the solution $\varphi (u,v)=uv(v-u)$, which
gives the Poisson-Lie structure for $d=1$ in Theorem 4.4, from
(4.16) in the following way. Rewriting (4.16) as
$$\varphi_{d,\lambda }(u,v)=\frac {uv(v^d-u^d)}{(1-\lambda u)(1-\lambda
v)}+\frac {\lambda du^2v^2}{(1-\lambda u)(1-\lambda v)}
\frac {u^{d-1}-v^{d-1}}{d-1},$$
we pass to the limit $d\to 1$ and then set $\lambda =0$.
\endremark
\demo{Proof} With the assumptions of Lemma 4.11 we have
$$f_d(u)=\sum_{n=d+1}^{\infty}\lambda_{1n}u^n =\sum_{n=d+1}^{\infty}\frac
{(\lambda_{1,d+2})^{n-d-1}}{(\lambda_{1,d+1})^{n-d-2}}u^n
=\lambda_{1,d+1}u^{d+1}\sum_{n=0}^{\infty}\biggl(\frac
{\lambda_{1,d+2}}{\lambda_{1,d+1}}u\biggr)^n
=\frac {\lambda_{1,d+1}}{1-\lambda u}u^{d+1},$$
and
$$\allowdisplaybreaks\align g_d(u,v)=-\sum_{n=1}^{d}\lambda_{n,d+1}u^n
&=-\lambda_{1,d+1}u+\frac {1}{d-1}\sum_{n=2}^{d}\frac
{(\lambda_{1,d+2})^{n-1}}{(\lambda_{1,d+1})^{n-2}}u^n\\
&=-\lambda_{1,d+1}u+\frac {\lambda_{1,d+2}u^2}{d-1}\sum_{n=2}^{d}\biggl(\frac
{\lambda_{1,d+2}}{\lambda_{1,d+1}}u\biggr)^{n-2}\\
&=-\lambda_{1,d+1}u+\frac {1}{d-1}\lambda_{1,d+2}u^2\left[\frac {1-(\lambda
u)^{d-1}}{1-\lambda u}\right]\\
&=\frac {\lambda_{1,d+1}}{(d-1)(1-\lambda u)}\left[-(d-1)(1-\lambda u)u+\lambda
u^2\left[1-(\lambda u)^{d-1}\right]\right]\\
&=\frac {\lambda_{1,d+1}}{(d-1)(1-\lambda u)}\biggl[d\lambda
u^2-u\bigl(d-1+(\lambda u)^d\bigr)\biggr].\endalign$$
By Theorem 4.5 we obtain
$$\aligned\varphi_{d,\lambda} (u,v)&=\frac
{1}{\lambda_{1,d+1}}\biggl[f_d(u)g_d(v)-f_d(v)g_d(u)\biggr]\\
&=\frac {1}{(d-1)(1-\lambda u)(1-\lambda v)}\biggl\{ v^{d+1}\left[ -\lambda
du^2+u(d-1+(\lambda u)^d)\right]-\\
&\qquad\qquad\qquad\qquad\qquad\qquad\qquad\qquad{ -u^{d+1}\left[ -\lambda
dv^2+v(d-1+(\lambda
v)^d)\right]\biggr\}.}\endaligned$$ After simplification we finally have
$$\varphi_{d,\lambda }(u,v)=\frac {1}{(d-1)(1-\lambda u)(1-\lambda v)}\biggl\{
(d-1)uv(v^d-u^d)+\lambda du^2v^2(u^{d-1}-v^{d-1})\biggr\}.\qed$$\enddemo
\remark{Remark} One can independently verify that (4.16) solves the functional
(partial) differential equation by directly substituting (4.16) into the
equation.\endremark
\par Summarizing, we showed existence of an infinite parameter family
of Poisson-Lie structures on the group $G_{\infty}$. In the next section we
shall show that
this family exhausts all such possible structures on $G_{\infty}$.
}

\head 5. The group $G_{\infty}$ and the $r$-matrix\endhead

{
\par In this section we describe the correspondence between the solution
$${\Omega}(u,v;{x})=\varphi(u,v){x}^{\prime}(u){x}^{\prime}(v)-\varphi({x}(u),{x}(v)) \tag 5.1$$
of the cocycle equation and the classical $r$-matrix on $\Cal
G_{\infty}$ \cite{STS2}.
It turns out that there is a one-to-one correspondence between the Poisson-Lie
structures on $G_{\infty}$ given by (5.1) and $r$-matrices on $\Cal
G_{\infty}$.
Namely, if we write $\varphi
(u,v)=\sum_{m,n=1}^{\infty}\lambda_{mn}u^mv^n$, where
$\lambda_{mn}=-\lambda_{nm}$, and $r=r_{ij}e_i\wedge e_j$ be the
classical $r$-matrix, where $r_{ij}=-r_{ji}$ and
$\{e_i\}_{i\ge 0}$ form a basis of the Lie algebra $\Cal G_{\infty}$, one shows
that $\lambda_{i+1,j+1}=r_{ij}$, $\forall_{i,j\ge 0}$. This we will prove by
demonstrating that the $\lambda_{ij}$'s and the $r_{ij}$'s
satisfy the same infinite system of algebraic equations. In what follows,
$\{x_i\}_{i\ge 1}$ will be again a set of local coordinates of a point of the
group $G_{\infty}$. We recall also
that a $1$-cochain $\alpha\:\Cal G_{\infty}\to \Cal G_{\infty}\wedge\Cal
G_{\infty}$ acting on a basis element $e_n\in \Cal G_{\infty}$ is written as
$\alpha (e_n)=\alpha_{ij}^{n}e_i\wedge e_j$, where summation is understood over
the repeated indices. If $\alpha$ is $1$-cocycle, then
$$(\delta\alpha)(e_l,e_m)=e_l.\alpha (e_m)-e_m.\alpha (e_l)-\alpha
(\left[e_l,e_m\right])=0.\tag 5.2$$
where $\delta$ is the coboundary operator in the Chevalley-Eilenberg
cohomology of Lie algebras.
\par Let us assume for a momment that $\Cal G$ is a finite-dimensional Lie
algebra. Let $r=r_{ij}e_i\wedge e_j\in \wedge^2\Cal G$ be a 0-cochain, and let
$\alpha\:\Cal G\to \Cal G\wedge\Cal G$ be defined as $\alpha=\delta r$. Let us
also define $<r,r>\in \otimes^3\Cal G$ as
$$<r,r>:=\left[r^{12},r^{13}\right]+\left[r^{12},r^{23}\right]+\left[r^{13},r^{23}\right], \tag 5.3$$
where we have
$$\left[r^{12},r^{13}\right]:= r_{ij}r_{kl}\left[e_i,e_k\right]\wedge e_j\wedge
e_l=r_{ij}r_{kl}C_{n}^{ik}e_n\wedge e_j\wedge e_l,$$
$$\left[r^{12},r^{23}\right]:= r_{ni}r_{kl}e_n\wedge \left[e_i,e_k\right]\wedge
e_l=r_{ni}r_{kl}C_{j}^{ik}e_n\wedge e_j\wedge e_l,$$
$$\left[r^{13},r^{23}\right]:= r_{ni}r_{jk}e_n\wedge e_j\wedge
\left[e_i,e_k\right]=r_{ni}r_{jk}C_{l}^{ik}e_n\wedge e_j\wedge e_l.$$
\flushpar In the above expressions we used
$\left[e_i,e_k\right]=C_{n}^{ik}e_n$, where $C_{n}^{ki}$ are the structure
constants of the Lie algebra $\Cal G$. Now, we can rewrite (5.3) in tensor
notation as
$$<r,r>=[C_{n}^{ik}r_{ij}r_{kl}+C_{j}^{ik}r_{il}r_{kn}+C_{l}^{ik}r_{in}r_{kj}]e_n\wedge e_j\wedge e_l. \tag 5.4$$
\par We have the following:
\proclaim{Lemma 5.1 [D2]} The coboundary $\alpha$ satisfies the co-Jacobi
identity, i.e\., $\alpha$ defines a Lie-bialgebra structure on $\Cal G$, if and
only if
$<r,r>$ is $\Cal G$-invariant with respect to the adjoint action of $\Cal G$ on
itself.\endproclaim
\demo{Proof} Let $\alpha(e_n)=\alpha_{ij}^ne_i\wedge e_j$. Then
$\alpha$ satisfies the co-Jacobi identity (cf. (iii), Def\. 3.1, Sec\. 1) if
and only if
$$\alpha_{ij}^n\alpha_{sp}^j+\alpha_{pj}^n\alpha_{is}^j+\alpha_{sj}^n\alpha_{pi}^j=0.\tag 5.5$$
Using the fact that $\alpha$ is a coboundary $\alpha(e_n)=\delta
r(e_n)=(r_{is}C_j^{ns}+r_{sj}C_i^{ns})e_i\wedge e_j$ we rewrite (5.5) as
$$\multline(r_{is}C_j^{ns}+r_{sj}C_i^{ns})(r_{kp}C_l^{jp}+r_{pl}C_k^{jp})+\\
+(r_{ks}C_j^{ns}+r_{sj}C_k^{ns})(r_{lp}C_i^{jp}+r_{pi}C_l^{jp})+\\
+(r_{ls}C_j^{ns}+r_{sj}C_l^{ns})(r_{ip}C_k^{jp}+r_{pk}C_i^{jp})=0.\endmultline$$
This system of equations is equivalent to
$$\align
&r_{is}C_j^{ns}r_{kp}C_l^{jp}+r_{sj}C_i^{ns}r_{kp}C_l^{jp}+r_{is}C_j^{ns}r_{pl}C_k^{jp}+r_{sj}C_i^{ns}r_{pl}C_k^{jp}\\
&+r_{ks}C_j^{ns}r_{lp}C_i^{jp}+r_{sj}C_k^{ns}r_{lp}C_i^{jp}+r_{ks}C_j^{ns}r_{pi}C_l^{jp}+r_{sj}C_k^{ns}r_{pi}C_l^{jp}\\
&+r_{ls}C_j^{ns}r_{ip}C_k^{jp}+r_{sj}C_l^{ns}r_{ip}C_k^{jp}+r_{ls}C_j^{ns}r_{pk}C_i^{jp}+r_{sj}C_l^{ns}r_{pk}C_i^{jp}=0.\endalign$$
Next we perform some algebraic manipulations using the Jacobi identity for the
Lie algebra structure constants of $\Cal G$
$$C_m^{ij}C_n^{mk}+C_m^{jk}C_n^{mi}+C_m^{ki}C_n^{mj}=0.$$
We observe that
$$\align
r_{is}C_j^{ns}r_{kp}C_l^{jp}+r_{ks}C_j^{ns}r_{pi}C_l^{jp}&=r_{is}C_j^{ns}r_{kp}C_l^{jp}+r_{kp}C_j^{np}r_{si}C_l^{js}\\
&=r_{is}r_{kp}(C_j^{ns}C_l^{jp}-C_j^{np}C_l^{js})\\
&=-r_{is}r_{kp}(C_j^{sn}C_l^{jp}+C_j^{np}C_l^{js})\\
&=r_{is}r_{kp}C_j^{ps}C_l^{jn}.\endalign$$
Similarly we have
$$r_{is}C_j^{ns}r_{pl}C_k^{jp}+r_{ls}C_j^{ns}r_{ip}C_k^{jp}=r_{is}r_{pl}C_j^{ps}C_k^{jn},$$
and
$$r_{ks}C_j^{ns}r_{lp}C_i^{jp}+r_{ls}C_j^{ns}r_{pk}C_i^{jp}=r_{ks}r_{lp}C_j^{ps}C_i^{jn}.$$
Therefore the co-Jacobi identity is equivalent to
$$\aligned
&r_{is}r_{kp}C_j^{ps}C_l^{jn}+r_{sj}C_i^{ns}r_{kp}C_l^{jp}+r_{sj}C_i^{ns}r_{pl}C_k^{jp}\\
+&r_{sj}C_k^{ns}r_{lp}C_i^{jp}+r_{sj}C_k^{ns}r_{pi}C_l^{jp}+r_{is}r_{pl}C_j^{ps}C_k^{jn}\\
+&r_{ks}r_{lp}C_j^{ps}C_i^{jn}+r_{sj}C_l^{ns}r_{ip}C_k^{jp}+r_{sj}C_l^{ns}r_{pk}C_i^{jp}=0.\endaligned\tag 5.6$$
\par Let us now consider the system of equations $e_m.<r,r>=0$ for every $m\ge
0$. That is,
$$[C_{n}^{ik}r_{ij}r_{kl}+C_{j}^{ik}r_{ni}r_{kl}+C_{l}^{ik}r_{ni}r_{jk}]e_m.(e_n\wedge e_j\wedge e_l)=0. \tag 5.7$$
Calculating
$$\align e_m.(e_n\wedge e_j\wedge e_l)&=[e_m,e_n]\wedge e_j\wedge e_l+e_n\wedge
[e_m,e_j]\wedge e_l+e_n\wedge e_j\wedge [e_m,e_l]\\
&=C_s^{mn}e_s\wedge e_j\wedge e_l+C_s^{mj}e_n\wedge e_s\wedge
e_l+C_s^{ml}e_n\wedge e_j\wedge e_s,\endalign$$
and renaming indices when necessary we obtain that (5.7) is equivalent to
$$\aligned
&r_{ij}r_{kl}C_s^{ik}C_n^{ms}+r_{is}C_n^{ik}r_{kl}C_j^{ms}+r_{ij}C_n^{ik}r_{ks}C_l^{ms}\\
+&r_{si}r_{kl}C_j^{ik}C_n^{ms}+r_{ni}C_s^{ik}r_{kl}C_j^{ms}+r_{ni}C_j^{ik}r_{ks}C_l^{ms}\\
+&r_{si}r_{jk}C_l^{ik}C_n^{ms}+r_{ni}C_l^{ik}r_{sk}C_j^{ms}+r_{ni}C_s^{ik}r_{jk}C_l^{ms}=0.\endaligned\tag 5.8$$
Again after renaming indices we conclude that the above system of equations is
identical to (5.6). This concludes the proof.\qed\enddemo
\par A subclass of coboundary Lie-bialgebra structures is obtained when
$<r,r>=0$; written explicitly, this condition has the form
$$C_{n}^{ik}r_{ij}r_{kl}+C_{j}^{ik}r_{il}r_{kn}+C_{l}^{ik}r_{in}r_{kj}=0. \tag
5.9$$
This is the so-called Classical Yang-Baxter Equation (CYBE) [STS2].
\par For the Lie algebra $\Cal G_{\infty}$ the structure constants are
$C_{k}^{ij}=(i-j)\delta_{k}^{i+j}$, where $i,j,k\ge 1$. One easily sees that
the arguments given above apply in this case,
since the presence of the Kronecker symbol in the formula for the
structure constants as well as the fact that $r_{ij}=0$ whenever $i<0,\ j<0$
make all sums {\it finite}.
 Therefore (5.9) becomes
$$(i-k)\delta_{n}^{i+k}r_{ij}r_{kl}+(i-k)\delta_{j}^{i+k}r_{il}r_{kn}+(i-k)\delta_{l}^{i+k}r_{in}r_{kj}=0,$$
or
$$\sum_{k=0}^{\max
{(n,j,l)}}\left[(n-2k)r_{n-k,j}r_{kl}+(j-2k)r_{j-k,l}r_{kn}+(l-2k)r_{l-k,n}r_{kj}\right]=0. \tag 5.10$$
Using the fact that
$\sum_{k=0}^{n}(n-k)r_{n-k,j}r_{kl}=\sum_{k=0}^{n}kr_{kj}r_{n-k,l},$
we have
$$\align\sum_{k=0}^{n}(n-2k)r_{n-k,j}r_{kl}&=\sum_{k=0}^{n}\left[
(n-k)r_{n-k,j}r_{kl}-kr_{n-k,j}r_{kl}\right]\\
&=\sum_{k=0}^{n}k\left[r_{kj}r_{n-k,l}-r_{n-k,j}r_{kl}\right]\\
&=\sum_{k=0}^{n}k\left[r_{kj}r_{n-k,l}+r_{j,n-k}r_{kl}\right],\endalign$$
and similarly for the other two terms in (5.10).
Thus, (5.10) assumes the form
$$\sum_{k=0}^{\max
{(n,j,l)}}k\left[(r_{n-k,l}+r_{n,l-k})r_{kj}+(r_{j,n-k}+r_{j-k,n})r_{kl}+(r_{l,j-k}+r_{l-k,j})r_{kn}\right]=0. \tag 5.11$$
(In the above formulae we implicitly assume that $r_{ij}=0$ whenever $i<0$ or
$j<0$.)
This is the Classical Yang-Baxter Equation for $\Cal G_{\infty}$.
\par We now proceed with the proof of the following important result.
\proclaim{Lemma 5.2} For the Lie algebra $\Cal G_{\infty}$ the system of
equations $e_n.<r,r>=0$, $n\in\zZ_+$, implies $<r,r>=0$.\endproclaim
\remark{Remark} This means that all coboundary Lie-bialgebra
structures on $\Cal G_{\infty}$ are given by the solutions of the
CYBE.\endremark
\demo{Proof} For $\Cal G_{\infty}$ the equation $e_m.<r,r>=0$, for any
$m\in\zZ_+$, becomes
$A+B+C=0,$
where
$$\align
A=&(i-k)(m-s)\delta_{s}^{i+k}\delta_{n}^{m+s}r_{ij}r_{kl}+(i-k)(m-s)\delta_{n}^{i+k}\delta_{j}^{m+s}r_{is}r_{kl}+(i-k)(m-s)\delta_{n}^{i+k}\delta_{l}^{m+s}r_{ij}r_{ks}\\
=&(i-k)(m-i-k)\delta_{n}^{m+i+k}r_{ij}r_{kl}+(2i-n)(2m-j)r_{i,j-m}r_{n-i,l}+(2i-n)(2m-j)r_{ij}r_{n-i,l-m}\\
=&(2i+m-n)(2m-n)r_{ij}r_{n-m-i,l}+(2i-n)(2m-j)r_{i,j-m}r_{n-i,l}+(2i-n)(2m-j)r_{ij}r_{n-i,l-m},\endalign$$
$$B=(2i-j)(2m-n)r_{n-m,i}r_{j-i,l}+(2i-j+m)(2m-j)r_{ni}r_{j-m-i,l}+(2i-j)(2m-l)r_{ni}r_{j-i,l-m},$$
and
$$C=(2i-l)(2m-n)r_{n-m,i}r_{j,l-i}+(2i-l)(2m-j)r_{ni}r_{j-m,l-i}+(2i-l+m)(2m-l)r_{ni}r_{j,l-m-i}.$$
Therefore the system of equations $e_m.<r,r>=0$ is equivalent to the following
system of equations
$$\align&\sum_{i=0}^{n-m}(2i+m-n)(2m-n)r_{ij}r_{n-m-i,l}+\sum_{i=0}^{j-m}(2i-j+m)(2m-j)r_{ni}r_{j-m-i,l}+\\
&+\sum_{i=0}^{l-m}(2i-l+m)(2m-l)r_{ni}r_{j,l-m-i}\\
&+\sum_{i=0}^{n}[(2i-n)(2m-j)r_{i,j-m}r_{n-i,l}+(2i-n)(2m-j)r_{ij}r_{n-i,l-m}]\\
&+\sum_{i=0}^{j}[(2i-j)(2m-n)r_{n-m,i}r_{j-i,l}+(2i-j)(2m-l)r_{ni}r_{j-i,l-m}]\\
&+\sum_{i=0}^{l}[(2i-l)(2m-n)r_{n-m,i}r_{j,l-i}+(2i-l)(2m-j)r_{ni}r_{j-m,l-i}]=0.\endalign$$
\par In order to prove that $e_m.<r,r>=0$ implies $<r,r>=0$ it is enough to
prove this implication for $m=0$, i.e\., that
$e_0.<r,r>=0$ implies $<r,r>=0$. Let $m=0$ in the above system of equations.
Then we have
$$\align
-&\sum_{i=0}^{n}(2i-n)nr_{ij}r_{n-i,l}-\sum_{i=0}^{j}(2i-j)jr_{ni}r_{j-i,l}-\sum_{i=0}^{l}(2i-l)lr_{ni}r_{j,l-i}\\
+&\sum_{i=0}^{n}[-(2i-n)jr_{ij}r_{n-i,l}-(2i-n)lr_{ij}r_{n-i,l}]\\
+&\sum_{i=0}^{j}[-(2i-j)nr_{ni}r_{j-i,l}-(2i-j)lr_{ni}r_{j-i,l}]\\
+&\sum_{i=0}^{l}[-(2i-l)nr_{ni}r_{j,l-i}-(2i-l)jr_{ni}r_{j,l-i}]=0,\endalign$$
which is finally equivalent to
$$(n+j+l)\biggl\{\sum_{i=0}^{n}(2i-n)r_{ij}r_{n-i,l}+\sum_{i=0}^{j}(2i-j)r_{ni}r_{j-i,l}+\sum_{i=0}^{l}(2i-l)r_{ni}r_{j,l-i}\biggr\}=0.$$
Now, using the identities
$$\align\sum_{i=0}^{n}(2i-n)r_{ij}r_{n-i,l}&=\sum_{i=0}^{n}ir_{ij}r_{n-i,l}+\sum_{i=0}^{n}(i-n)r_{ij}r_{n-i,l}\\
&=\sum_{i=0}^{n}i(r_{ij}r_{n-i,l}-r_{n-i,j}r_{il})\\
&=\sum_{i=0}^{n}i(r_{ij}r_{n-i,l}+r_{j,n-i}r_{il}),\endalign$$
and similarly for
$$\sum_{i=0}^{j}(2i-j)r_{ni}r_{j-i,l}=\sum_{i=0}^{j}i(r_{ni}r_{j-i,l}+r_{j-i,n}r_{il}),$$
and for
$$\sum_{i=0}^{l}(2i-l)r_{ni}r_{j,l-i}=\sum_{i=0}^{l}i(r_{ni}r_{j,l-i}+r_{n,l-i}r_{ij}),$$
we conclude that $e_0.<r,r>=0$ is equivalent to
$$(n+j+l)\Bigl\{\sum_{i=0}^{\max
{(n,j,l)}}i\left[(r_{n-i,l}+r_{n,l-i})r_{ij}+(r_{j,n-i}+r_{j-i,n})r_{il}+(r_{l,j-i}+r_{l-i,j})r_{in}\right]\Bigr\}=0.$$
\par Since $n+j+l=0$ if and only if $n=j=l=0$, and in this case $e_m.<r,r>$ is
identically zero for any $m\in\zZ_+$, the above system of equations implies
that the expression in the curly brackets vanishes, and by formula (5.11) this
is exactly $<r,r>=0$.
\qed\enddemo
\par As we have seen, if $\alpha$ is a coboundary of $r$ then $\alpha$ has the
form
$$\alpha (e_n)=(\delta r)(e_n)=e_n.r=r_{ij}(\left[e_n ,e_i\right]\wedge
e_j+e_i\wedge\left[ e_n
,e_j\right])=(r_{lj}C_{i}^{nl}+r_{il}C_{j}^{nl})e_i\wedge e_j.$$
Therefore,
$$\alpha_{ij}^{n}=r_{lj}C_{i}^{nl}+r_{il}C_{j}^{nl}=(2n-i)r_{i-n,j}+(2n-j)r_{i,j-n}.\tag 5.12$$
\par Now, we turn our attention to the Poisson bracket on the group, written in
local coordinates, by writing ${\Omega}(u,v;{x})$ in components. In order to do
this we shall need the formula
$$
x(u)^n=\sum_{s_1=1}^{\infty}x_{s_1}u^{s_1}\ldots\sum_{s_n=1}^{\infty}x_{s_n}u^{s_n}=\sum_{s_1=1}^{\infty}\ldots\sum_{s_n=1}^{\infty}x_{s_1}\ldots x_{s_n}u^{s_1+\ldots +s_n}=\sum_{i=n}^{\infty}\left(\sum_{\bigl(\sum_{k=1}^{n}s_k\bigr)=i}x_{s_1}\ldots x_{s_n}\right)u^i.$$
Then we have
$${\Omega}(u,v;{x})=\varphi(u,v){x}^{\prime}(u){x}^{\prime}(v)-\varphi({x}(u),{x}(v))=$$
$$\multline
=\sum_{p,q=1}^{\infty}\lambda_{pq}u^pv^q\sum_{i=1}^{\infty}ix_iu^{i-1}\sum_{j=1}^{\infty}jx_jv^{j-1}-\\
-\sum_{p,q=1}^{\infty}\lambda_{pq}\sum_{i=p}^{\infty}\left(\sum_{\bigl(\sum_{k=1}^{p}r_k\bigr)=i}x_{r_1}\ldots x_{r_p}\right)u^i\sum_{j=q}^{\infty}\left(\sum_{\bigl(\sum_{l=1}^{q}s_l\bigr)=j}x_{s_1}\ldots x_{s_q}\right)v^j\endmultline$$
$$\multline
=\sum_{p,q=1}^{\infty}\sum_{i=1}^{\infty}\sum_{j=1}^{\infty}\lambda_{pq}ix_ijx_ju^{p+i-1}v^{q+j-1}-\\
-\sum_{p,q=1}^{\infty}\sum_{i=p}^{\infty}\sum_{j=q}^{\infty}\lambda_{pq}\left(\sum_{\bigl(\sum_{k=1}^{p}r_k\bigr)=i}x_{r_1}\ldots x_{r_p}\sum_{\bigl(\sum_{l=1}^{q}s_l\bigr)=j}x_{s_1}\ldots x_{s_q}\right)u^iv^j\endmultline$$
$$=\sum_{p,q=1}^{\infty}\sum_{i=p}^{\infty}\sum_{j=q}^{\infty}\left(\lambda_{i-p+1,j-q+1}px_pqx_q-\lambda_{pq}\sum_{\bigl(\sum_{k=1}^{p}r_k\bigr)=i}x_{r_1}\ldots x_{r_p}\sum_{\bigl(\sum_{l=1}^{q}s_l\bigr)=j}x_{s_1}\ldots x_{s_q}\right)u^iv^j$$
$$=\sum_{i,j=1}^{\infty}\left[\sum_{p=1}^{i}\sum_{q=1}^{j}\left(\lambda_{i-p+1,j-q+1}px_pqx_q-\lambda_{pq}\sum_{\bigl(\sum_{k=1}^{p}r_k\bigr)=i}x_{r_1}\ldots x_{r_p}\sum_{\bigl(\sum_{l=1}^{q}s_l\bigr)=j}x_{s_1}\ldots x_{s_q}\right)\right]u^iv^j.$$
Therefore for $\{x_i,x_j\}=\omega_{ij}(x)$ we obtain
$$\multline\omega_{ij}(x)=\\
=\sum_{p=1}^{i}\sum_{q=1}^{j}px_pqx_q\lambda_{i-p+1,j-q+1}-\sum_{p=1}^{i}\sum_{q=1}^{j}\lambda_{pq}\left(\sum_{\bigl(\sum_{k=1}^{p}r_k\bigr)=i}x_{r_1}\ldots x_{r_p}\sum_{\bigl(\sum_{l=1}^{q}s_l\bigr)=j}x_{s_1}\ldots x_{s_q}\right).\endmultline\tag 5.13$$
\par Before we continue further, let us deduce the following useful formulae
$$\eqalign{\frac {\partial}{\partial
x_n}\sum_{\bigl(\sum_{k=1}^{p}r_k\bigr)=i}x_{r_1}\ldots
x_{r_p}&=\sum_{\bigl(\sum_{k=1}^{p}r_k\bigr)=i}\biggl(\sum_{l=1}^{p}x_{r_1}\ldots \delta_{r_l}^{n}\ldots x_{r_p}\biggr)\cr
&=\sum_{l=1}^{p}\left(\sum_{\bigl(\sum_{k=1}^{p}r_k\bigr)=i}x_{r_1}\ldots\delta_{r_l}^{n}\ldots x_{r_p}\right)\cr
&=p\sum_{\bigl(\sum_{k=1}^{p-1}r_k\bigr)=i-n}x_{r_1}\ldots x_{r_{p-1}},}$$
as well as
$$\sum_{\bigl(\sum_{k=1}^{p-1}r_k\bigr)=i-n}x_{r_1}\ldots
x_{r_{p-1}}{|}_{e}=\sum_{\bigl(\sum_{k=1}^{p-1}r_k\bigr)=i-n}\delta_{r_1}^1\ldots\delta_{r_{p-1}}^1=\delta_{i-n}^{p-1},$$
$$\biggl(\frac {\partial}{\partial
x_n}\sum_{\bigl(\sum_{k=1}^{p}r_k\bigr)=i}x_{r_1}\ldots
x_{r_p}\biggr){\biggm|}_{e}=p\delta_{i-n}^{p-1}.$$
\flushpar Differentiating (5.13) with respect to $x_n$ we obtain
$$\eqalign{\frac {\partial\omega_{ij}}{\partial
x_n}{\biggm|}_{x}&=\sum_{p=1}^{i}\sum_{q=1}^{j}(
p\delta_{p}^{n}qx_q\lambda_{i-p+1,j-q+1}+px_p\delta_{q}^{n}q\lambda_{i-p+1,j-q+1})\cr
&-\sum_{p=1}^{i}\sum_{q=1}^{j}\left[
p\sum_{\bigl(\sum_{k=1}^{p-1}r_k\bigr)=i-n}x_{r_1}\ldots
x_{r_{p-1}}\sum_{\bigl(\sum_{l=1}^{q}s_l\bigr)=j}x_{s_1}\ldots
x_{s_{q}}\right]\lambda_{pq}\cr
&-\sum_{p=1}^{i}\sum_{q=1}^{j}\left[
q\sum_{\bigl(\sum_{k=1}^{p}r_k\bigr)=i}x_{r_1}\ldots
x_{r_{p}}\sum_{\bigl(\sum_{l=1}^{q-1}s_k\bigr)=j-n}x_{s_1}\ldots
x_{s_{q-1}}\right]\lambda_{pq}.}$$
 From the above formula we have (keeping in mind that
$x_p{|}_{e}=\delta_{p}^{1}$)
$$\eqalign{\beta_{ij}^{n}:=\frac {\partial\omega_{ij}}{\partial
x_n}{\biggm|}_{e}&=\sum_{p=1}^{i}p\delta_{p}^{n}\lambda_{i-p+1,j}+\sum_{q=1}^{j}q\delta_{q}^{n}\lambda_{i,j-q+1}-\sum_{p=1}^{i}\sum_{q=1}^{j}\lambda_{pq}\left[p\delta_{i-n+1}^{p}\delta_{j}^{q}+q\delta_{i}^{p}\delta_{j-n+1}^{q}\right]\cr
&=\sum_{p=1}^{i}(p\delta_{p}^{n}\lambda_{i-p+1,j}-p\delta_{i-n+1}^{p}\lambda_{pj})+\sum_{q=1}^{j}(q\delta_{q}^{n}\lambda_{i,j-q+1}-q\delta_{j-n+1}^{q}\lambda_{iq})\cr
&=n\lambda_{i-n+1,j}-(i-n+1)\lambda_{i-n+1,j}+n\lambda_{i,j-n+1}-(j-n+1)\lambda_{i,j-n+1}\cr
&=(2n-i-1)\lambda_{i-n+1,j}+(2n-j-1)\lambda_{i,j-n+1}.}$$
Thus, we finally obtain that
$$\beta_{ij}^{n}=(2n-i-1)\lambda_{i-n+1,j}+(2n-j-1)\lambda_{i,j-n+1},\tag
5.14$$
which is the same as (5.12) after we make the identification
$r_{i,j-1}=\lambda_{i+1,j}$, for every $i\ge 0,j\ge 1$, which is
equivalent to $r_{ij}=\lambda_{i+1,j+1}$, for every $i,j\ge 0$, since $r_{ij}$
and $\lambda_{ij}$ are antisymmetric.
\par On the other hand the system of equations (5.11) is exactly the
same as the system of equations (4.6) for the entries of the infinite
matrix $\Lambda=(\lambda_{ij})$. Thus $\Lambda$ and $r$ differ only
by a non-zero multiple.
\par The results of the calculations made in this section may be
summarized in the following theorem.
\proclaim{Theorem 5.3} There is a one-to-one correspondence between the
coboundary Lie bialgebra
structures on $\Cal G_{\infty}$ given by $r$ and the Poisson-Lie
structures of the type {\rm (5.1)} on $G_{\infty}$.
Since all Lie bialgebra structures on $\Cal G_{\infty}$ are given by
$r$ (cf. Theorem {\rm 3.1}), Theorem {\rm 4.5} gives
a classification of all solutions of the classical Yang-Baxter equation for
$\Cal G_{\infty}$.\endproclaim
\demo{Proof} Recall that $\omega_{mn}$ satisfy the infinite system of
functional equations (1.6):
$$\omega _{mn}({z})=\omega _{kl}(x)\pd {z_m}{x_k} \pd {z_n}{x_l} +
\omega_{kl}{(y)}\pd {z_m}{y_k} \pd {z_n}{y_l} , \ \ \ \text{where}\ \ x,y\in
G_{\infty},\tag 5.15$$
and $z_n=z_n(x,y)$ is given by formula (2.1)
$$z _n =\sum_{i=1}^n
x_i\sum_{\bigl(\sum_{\alpha=1}^{i}{j_{\alpha}}\bigr)=n}{y_{j_1}\ldots
y_{j_i}}.\tag 5.16$$
 From (5.16) it follows that
$$\pd {z_i}{y_k}\biggm|_{y=e}=(i-k+1)x_{i-k+1}\ \ \ \ \ \text{and}\ \ \ \ \ \pd
{z_i}{x_k}\biggm|_{y=e}=\delta^k_i.\tag 5.17$$
\par Let us fix $n\in\zN$ and consider a subsystem of the system of equations
(5.15) for all $\omega_{ij}$ with $1\le i<j\le n$.
After differentiating (5.15) with respect to $y_j$, for each $j$ such that
$1\le j\le n$, and setting $y=e$ we deduce that $\omega_{mn}$ satisfy the
following inhomogeneous system of linear partial differential equations
$$\multline\sum_{i=j}^n(i+1-j)x_{i+1-j}\pd
{\omega_{mn}}{x_i}=\omega_{m+1-j,n}(x)(m+1-j)+\omega_{m,n+1-j}(x)(n+1-j)+\\
+\sum_{k=1}^m\sum_{l=1}^n\beta^j_{kl}(m+1-k)(n+1-l)x_{m+1-k}x_{n+1-l},\endmultline\tag 5.18$$
for $1\le j\le n$, and where $\beta^j_{kl}=\pd {\omega_{kl}}{y_j}\Bigm|_{y=e}$.
\par The idea of the proof is as follows. Let $\beta^j_{kl}$ be given by
(5.14). For each $n\in\zN$ the general solution of (5.18) is a
linear combination of the general solution of the homogeneous system of
equations
$$\sum_{i=j}^n(i+1-j)x_{i+1-j}\pd
{\omega_{mn}}{x_i}=\omega_{m+1-j,n}(x)(m+1-j)+\omega_{m,n+1-j}(x)(n+1-j)\tag
5.19$$
and a particular solution of the inhomogeneous system (5.18). We now show that
for each $n\in\zN$ and $1\le m< n$ the system (5.18) has a
unique solution by demonstrating that the only solution of the homogeneous
system is the zero solution. Therefore, since every solution
of the system of functional equations (5.15) is a solution of the
system of partial differential equations (5.18) it follows that the class of
solutions of (5.15) found in Sec\. 4 exhausts all possible solutions of (5.15).
We will prove that the only solution of (5.19) is the zero solution by
induction applied in several steps. Recall that
 $$\omega_{mn}(e)=0,\ \ \ \ \text{for every}\ \ n,m\in\zN.\tag 5.20$$
In the following arguments we implicitly assume that $\omega_{mn}=0$ whenever
$n<1$ or $m<1$.
\par (i) If $n=1$ there is nothing to prove. Let $n=2$. Then from (5.19) we
obtain
$$x_1\pd {\omega_{12}}{x_1}=3\omega_{12}.\tag 5.21$$
The most general solution of this equation is
$\omega_{12}(x)=Cx_1^3,$
where $C$ is an arbitrary constant.
 From (5.20) it follows that $C=0$. Therefore $\omega_{12}(x)=0$ is the only
solution of (5.21). Let $n=3$. Then from (5.19) we obtain
$$\align x_1\pd {\omega_{13}}{x_1} +2x_2\pd
{\omega_{13}}{x_2}&=4\omega_{13}\tag 5.22\\
x_1\pd {\omega_{13}}{x_2}&=0.\tag 5.23\endalign$$
 From (5.23) it follows that $\omega_{13}(x)=\omega_{13}(x_1)$. From (5.22) we
deduce that $\omega_{13}(x_1)$ satisfies the
equation
$$x_1\pd {\omega_{13}}{x_1}=4\omega_{13}.\tag 5.24$$
Therefore $\omega_{13}(x)=Cx_1^5$, and from (5.20) it follows that $C=0$, and
$\omega_{13}(x)=0$. Let us assume now that $\omega_{1k}(x)=0$
for $2\le k\le n-1$. From (5.19) we have
$$\sum_{i=j}^n(i+1-j)x_{i+1-j}\pd
{\omega_{1n}}{x_i}=\omega_{2-j,n}(x)(2-j)+\omega_{1,n+1-j}(x)(n+1-j),\ \ \
\text{for}\ \ 1\le j\le n, \tag 5.25$$
which implies
$$\align\sum_{i=1}^nix_{i}\pd {\omega_{1n}}{x_i}&=(n+1)\omega_{1n}(x),\tag
5.26\\
\sum_{i=j}^n(i+1-j)x_{i+1-j}\pd {\omega_{1n}}{x_i}&=0,\ \ \ \text{for} \ \ 2\le
j\le n.\tag 5.27\endalign$$
 We used above the induction hypothesis: $\omega_{1k}(x)=0$
for $2\le k\le n-1$ from which follows that the r.h.s of (5.27) is zero. From
(5.27) it follows that $\omega_{1n}(x)=\omega_{1n}(x_1)$. Then from (5.26) it
follows that $\omega_{1n}(x_1)$ satisfies
$$x_1\pd {\omega_{1n}}{x_1}=(n+1)\omega_{1n}.\tag 5.28$$
 From (5.28) we have that $\omega_{1n}(x)=Cx_1^{n+1}$ for an arbitrary constant
$C$. Applying again (5.20) we conclude that $\omega_{1n}(x)=0$.
Therefore $\omega_{1n}(x)=0$ for every $n\in\zN$.
\par (ii) Let $m=2$ and $n=3$. Then from (5.19) we have the following
homogeneous system of partial differential equations for $\omega_{23}$:
$$\align x_1\pd {\omega_{23}}{x_1} +2x_2\pd {\omega_{23}}{x_2}+3x_3\pd
{\omega_{23}}{x_3}&=5\omega_{23}\\
x_1\pd {\omega_{23}}{x_2} +2x_2\pd {\omega_{23}}{x_3}&=\omega_{13}=0\\
x_1\pd {\omega_{23}}{x_3}&=-\omega_{12}=0.\endalign$$
Arguing in a similar manner as above we obtain that $\omega_{23}(x)=0$. Let us
assume that $\omega_{2k}(x)=0$ for all $k$ such that
$3\le k\le n-1$. We now prove that $\omega_{2n}(x)=0$. From (5.19) we have
$$\sum_{i=j}^n(i+1-j)x_{i+1-j}\pd
{\omega_{2n}}{x_i}=\omega_{3-j,n}(x)(3-j)+\omega_{2,n+1-j}(x)(n+1-j),\ \ \
\text{for}\ \ 1\le j\le n. \tag 5.29$$
After using the induction hypothesis and the already proved fact that
$\omega_{1n}=0$, for every $n\in\zN$, (5.29) yields
$$\align\sum_{i=1}^nix_{i}\pd {\omega_{2n}}{x_i}&=(n+2)\omega_{2n}(x),\tag
5.30\\
\sum_{i=j}^n(i+1-j)x_{i+1-j}\pd {\omega_{2n}}{x_i}&=0,\ \ \ \text{for} \ \ 2\le
j\le n.\tag 5.31\endalign$$
Therefore from (5.31) and (5.30) it follows that
$\omega_{2n}(x)=\omega_{2n}(x_1)=Cx_1^{n+2}$, and imposing (5.20) again we
obtain that
$\omega_{2n}(x_1)=0$. Thus $\omega_{2n}(x)=0$ for every $n\in\zN$.
\par (iii) Let us assume that $\omega_{sn}=0$ for all $s$ such that $1\le s\le
m-1$, for some $m\ge 2$ and all $n>s$. We will prove that
$\omega_{mn}=0$ for all $n\ge m$. Let $n=m+1$. From (5.19) we have
$$\multline\sum_{i=j}^{m+1}(i+1-j)x_{i+1-j}\pd
{\omega_{m,m+1}}{x_i}=\omega_{m+1-j,n}(x)(m+1-j)+\omega_{m,m+2-j}(x)(m+2-j),\\
\text{for}\ \ 1\le j\le m+1.\endmultline \tag 5.32$$
We apply now the induction hypothesis and deduce from (5.32) the following
system of equations
$$\align\sum_{i=1}^{m+1}ix_{i}\pd
{\omega_{m,m+1}}{x_i}&=(2m+1)\omega_{m,m+1}(x),\tag 5.33\\
\sum_{i=j}^{m+1}(i+1-j)x_{i+1-j}\pd {\omega_{m,m+1}}{x_i}&=0,\ \ \ \text{for} \
\ 2\le j\le m+1.\tag 5.34\endalign$$
 From (5.34) it follows that $\omega_{m,m+1}(x)=\omega_{m,m+1}(x_1)$, and from
(5.33) we deduce that $\omega_{m,m+1}(x)$ must satisfy
$$x_1\pd {\omega_{m,m+1}}{x_1}=(2m+1)\omega_{m,m+1}.\tag 5.35$$
The solution of the above equation is $\omega_{m,m+1}(x)=Cx_1^{2m+1}$, where
$C$ is an arbitrary constant. Then from $\omega_{m,m+1}(e)=C=0$ we obtain that
$\omega_{m,m+1}(x)=0$. Finally, we assume that $\omega_{mk}=0$ for all $k$ such
that $m+1\le k\le n-1$, and we prove it for $k=n$.
Indeed,  from (5.19), after applying the induction hypothesis, we obtain
$$\align\sum_{i=1}^nix_{i}\pd {\omega_{mn}}{x_i}&=(m+n)\omega_{mn}(x),\tag
5.36\\
\sum_{i=j}^{n}(i+1-j)x_{i+1-j}\pd {\omega_{mn}}{x_i}&=0,\ \ \ \text{for} \ \
2\le j\le n.\tag 5.37\endalign$$
Again, from (5.37) it follows that $\omega_{mn}(x)=\omega_{mn}(x_1)$, and that
$\omega_{mn}(x)$ must satisfy
$$x_1\pd {\omega_{mn}}{x_1}=(m+n)\omega_{mn}.\tag 5.38$$
 From here we conclude that $\omega_{mn}(x)=Cx_1^{m+n}$ for an arbitrary
constant $C$. But the requirement $\omega_{mn}(e)=0$ fixes
the value of this constant to be $C=0$. Therefore $\omega_{mn}(x)=0$.
\par Thus, we showed that for every $m,n\in\zN$ the only solution of (5.19) is
the zero solution. Therefore the system of partial
differential equations (5.18) has a unique solution. The existence of the
solution follows from the existence of the solution
of the system of functional equations (5.15) of which (5.18) is a consequence.
Thus, the structure constants $\beta^j_{kl}$ of the Lie-bialgebra $\Cal
G_{\infty}$, as given by (5.14), determine uniquely all Poisson-Lie structures
on the group $G_{\infty}$.
The proof of Theorem 5.3 is completed.
\qed\enddemo
\par We conclude this section by writing an explicit formula for the family of
Lie-bialgebra structures arising from the family of Poisson-Lie structures
obtained in Theorem 4.4, as well as the more general one-parameter family, of
which it is a particular case for $d\ge 1$.
An elegant way to do this is by deriving a global formula for the Lie-bialgebra
structures on $\Cal G_{\infty}$ in terms of
generating series and solutions of
$$\varphi(u,v)\bigl[
{\partial}_u\varphi(w,u)+{\partial}_v\varphi(w,v)\bigr]+c.p.=0.$$
Let us define $\Cal A_n(u,v)$ as
$$\Cal A_n(u,v):=\frac {\partial}{\partial {x_n}}{\Omega}(u,v;{\Cal
X})\biggm|_e=\sum_{i,j=1}^{\infty}\alpha_{ij}^{n}u^iv^j.$$
Then we have the following lemma.
\proclaim{Lemma 5.4} The generating series $\Cal A_n(u,v)$ is given by
$$\Cal
A_n(u,v)=n\varphi(u,v)\bigl(u^{n-1}+v^{n-1}\bigr)-\bigl[u^n{\partial}_u\varphi(u,v)+v^n{\partial}_v\varphi(u,v)\bigr].$$
\endproclaim
\demo{Proof} We use formula (5.1) and the following facts.  If
$x(u)=\sum_{i=1}^{\infty}x_iu^i$ then
$\Cal X'(u)\bigm|_x=\sum_{i=1}^{\infty}ix_iu^{i-1} \implies \Cal
X'(u)\bigm|_e=1,$
and also
$\frac {\partial}{\partial {x_n}}\Cal
X'(u)\bigm|_x=\sum_{i=1}^{\infty}i\delta^n_iu^{i-1}=nu^{n-1}.$
 From this it follows that
$$\frac {\partial}{\partial {x_n}}\biggl[\Cal X'(u)\Cal
X'(v)\biggr]\biggm|_x=nu^{n-1}x'(v)+nv^{n-1}x'(u)\implies \frac
{\partial}{\partial {x_n}}\biggl[\Cal X'(u)\Cal
X'(v)\biggr]\biggm|_e=n(u^{n-1}+v^{n-1}).$$
Finally we have
$$\eqalign{\frac {\partial}{\partial {x_n}}\varphi\bigl(\Cal X(u),\Cal
X(v)\bigr)\biggm|_e&=
\frac {\partial\Cal X(u)}{\partial {x_n}}\partial_1\varphi\bigl(\Cal X(u),\Cal
X(v)\bigr)\biggm|_e+
\frac {\partial\Cal X(v)}{\partial {x_n}}\partial_2\varphi\bigl(\Cal X(u),\Cal
X(v)\bigr)\biggm|_e\cr\cr
&=u^n{\partial}_u\varphi(u,v)+v^n{\partial}_v\varphi(u,v).}$$
In the above equality we used that $\frac {\partial\Cal X(u)}{\partial
{x_n}}\bigm|_e=u^n$.
\qed
\enddemo
\proclaim{Proposition 5.5 [Mi,Ta]} For each $d\in \zN$ the family of
Poisson-Lie structures (5.1) given by
$\varphi_d (u,v)=uv(v^d-u^d)$ gives rise
to the following family of Lie-bialgebra structures on $\Cal G_{\infty}$:
$$\alpha(e_n)=2ne_{d}\wedge e_n-2(n-d)e_0\wedge e_{d+n}, \qquad (n\ge 0), \tag
5.39$$
where $\{e_n\}_{n\in \zZ_{+}}$ is a basis for $\Cal G_{\infty}$.\endproclaim
\demo{Proof} The generating series $\Cal A_{n,d}$ in this case is
$$\eqalign{\Cal
A_{n,d}(u,v)=&n\bigl[uv^{d+1}-vu^{d+1}\bigr]\bigl(u^{n-1}+v^{n-1}\bigr)-\cr
&\qquad
-\biggl\{u^n\bigl[v^{d+1}-(d+1)vu^d\bigr]+v^n\bigl[(d+1)uv^d-u^{d+1}\bigr]\biggr\}\cr
=&(n-1)u^nv^{d+1}-(n-1)v^nu^{d+1}+(n-d-1)uv^{n+d}-(n-d-1)vu^{n+d}\cr
=&\biggl\{(n-1)\bigl[\delta^n_i\delta^{d+1}_j-\delta^n_j\delta^{d+1}_i\bigr]+
(n-d-1)\bigl[\delta^1_i\delta^{d+n}_j-\delta^1_j\delta^{d+n}_i\bigr]\biggr\}u^iv^j\cr
=&\sum_{i,j=1}^{\infty}\alpha_{ij|d}^{n}u^iv^j,}$$
where
$$\alpha_{ij|d}^{n}=\biggl\{(n-1)\bigl[\delta^n_i\delta^{d+1}_j-\delta^n_j\delta^{d+1}_i\bigr]+
(n-d-1)\bigl[\delta^1_i\delta^{d+n}_j-\delta^1_j\delta^{d+n}_i\bigr]\biggr\}.$$
Therefore
$$\eqalign{\alpha_d(e_n)=\alpha_{ij|d}^{n}e_i\wedge e_j&=(n-1)\bigl[e_n\wedge
e_{d+1}-e_{d+1}\wedge e_{n}\bigr]+(n-d-1)\bigl[e_1\wedge e_{d+n}-e_{d+n}\wedge
e_1\bigr]\cr
&=2(n-1)e_n\wedge e_{d+1}+2(n-d-1)e_1\wedge e_{n+d},}$$
and after shifting indices by 1 we obtain
$$\alpha_d(e_n)=-2ne_d\wedge e_{n}+2(n-d)e_0\wedge e_{n+d}\ \ \ \ \text{for
every}\ \ n\in\zZ_+.$$
A second way to derive the above formula is by using Theorem 5.3.
The $r$-matrix is given in this case by
$r_{ij}=\delta_{i+1}^{1}\delta_{j+1}^{d+1}-\delta_{i+1}^{d+1}\delta_{j+1}^{1}=\lambda_{i+1,j+1}$.
 Therefore, using (5.12), we have
$$\eqalign{\alpha (e_n)=\alpha_{ij}^{n}e_i\wedge
e_j=&\left[(2n-i)r_{i-n,j}+(2n-j)r_{i,j-n}\right] e_i\wedge e_j\cr
=&-(2n-i)\delta_{i-n+1}^{d+1}\delta_{j+1}^{1}e_i\wedge
e_j+(2n-i)\delta_{i-n+1}^{1}\delta_{j+1}^{d+1}e_i\wedge e_j\cr
&\qquad-(2n-j)\delta_{i+1}^{d+1}\delta_{j-n+1}^{1}e_i\wedge
e_j+(2n-j)\delta_{i+1}^{1}\delta_{j-n+1}^{d+1}e_i\wedge e_j\cr
=&-(n-d)e_{d+n}\wedge e_0+ne_n\wedge e_{d+1}-ne_{d+1}\wedge e_n+(n-d)e_0\wedge
e_{d+n}\cr
=&-2ne_{d}\wedge e_{n}+2(n-d)e_0\wedge e_{d+n}.}$$
This concludes the proof.\qed \enddemo
\proclaim{Proposition 5.6} For every $d\ge 2$ the family of
Poisson-Lie structures described by {\rm (4.15)} gives rise to the following
family of Lie-bialgebra structures on $\Cal G_{\infty}$:
$$\align\alpha_{d,\lambda}(e_n)=&2\sum_{i=d+n}^{\infty}(2n-i){\lambda}^{i-(n+d)}e_0\wedge e_i-2n\sum_{i=d}^{\infty}{\lambda}^{i-d}e_i\wedge e_n+\\
&\qquad +\frac
{2}{d-1}\sum_{i=d+n}^{\infty}\sum_{j=1}^{d-1}(2n-i)\lambda^{i+j-(n+d)}e_i\wedge
e_j\tag 5.40\\
 &\qquad +\frac
{2}{d-1}\sum_{i=d}^{\infty}\sum_{j=n+1}^{d+n-1}(2n-j)\lambda^{i+j-(n+d)}e_i\wedge e_j,\endalign$$
for every $n\in\zZ_+$.
\endproclaim
\demo{Proof} Let $\alpha_{d,\lambda}(e_n)=\alpha^n_{ij|d}e_i\wedge e_j$, where
$\{e_n\}_{n\ge 1}$ is a basis of $\Cal G_{\infty}$, and
$$\alpha^n_{ij|d}=(2n-i-1)\lambda_{i-n+1,j}+(2n-j-1)\lambda_{i,j-n+1},\
\ \ \ \text{for every}\ \ \ n,i,j\in\zN.$$
With the assumptions of Lemma 4.11 and Theorem 4.12 we have
$$\lambda_{ij}=\frac
{1}{\lambda_{1,d+1}}\biggl[\lambda_{1i}\lambda_{d+1,j}-\lambda_{1j}\lambda_{d+1,i}\biggr],$$
where
$$\align \lambda_{n,d+1}&=0,\ \ \ \ \text{for every}\ \ \ n\ge d+1,\\
\lambda_{n,d+1}&=-\frac {1}{d-1}\frac
{(\lambda_{1,d+2})^{n-1}}{(\lambda_{1,d+1})^{n-2}}=-\frac
{1}{d-1}\lambda_{1,d+1}\lambda^{n-1},\ \ \
\text{for}\ \ \ 2\le n\le d,\tag 5.41\\
\lambda_{1n}&=\frac
{(\lambda_{1,d+2})^{n-d-1}}{(\lambda_{1,d+1})^{n-d-2}}=\lambda_{1,d+1}\lambda^{n-d-1},\ \ \ \ \text{for every}\ \ \ n\ge d+1,\endalign$$
and where we have introduced $\lambda:= \frac
{\lambda_{1,d+2}}{\lambda_{1,d+1}}$. Then
$$\allowdisplaybreaks\align
\alpha_{d,\lambda}(e_n)=&\Bigl[(2n-i-1)\lambda_{i-n+1,j}+(2n-j-1)\lambda_{i,j-n+1}\Bigr]e_i\wedge e_j\\
=&\frac
{1}{\lambda_{1,d+1}}\biggl\{(2n-i-1)\Bigl[\lambda_{1,i-n+1}\lambda_{d+1,j}-\lambda_{1j}\lambda_{d+1,i-n+1}\Bigr]+\\
&\qquad\qquad\qquad\qquad+(2n-j-1)\Bigl[\lambda_{1i}\lambda_{d+1,j-n+1}-\lambda_{1,j-n+1}\lambda_{d+1,i}\Bigr]\biggr\}e_i\wedge e_j\\
=&\frac
{1}{\lambda_{1,d+1}}\Biggl\{\sum_{i=d+n}^{\infty}\sum_{j=1}^{d}(2n-i-1)\lambda_{1,i-n+1}\lambda_{d+1,j}e_i\wedge e_j-\\
&\qquad-\sum_{i=n}^{d+n-1}\sum_{j=d+1}^{\infty}(2n-i-1)\lambda_{1j}\lambda_{d+1,i-n+1}e_i\wedge e_j\\
&\qquad+\sum_{i=d+1}^{\infty}\sum_{j=n}^{d+n-1}(2n-j-1)\lambda_{1i}\lambda_{d+1,j-n+1}e_i\wedge e_j-\\
&\qquad-\sum_{i=1}^{d}\sum_{j=d+n}^{\infty}(2n-j-1)\lambda_{1,j-n+1}\lambda_{d+1,i}e_i\wedge e_j\Biggr\}\\
=&\frac
{1}{\lambda_{1,d+1}}\Biggl\{\sum_{i=d+n}^{\infty}\sum_{j=2}^{d}(2n-i-1)\lambda_{1,i-n+1}\lambda_{d+1,j}e_i\wedge e_j+\\
&\qquad+\lambda_{d+1,1}\sum_{i=d+n}^{\infty}(2n-i-1)\lambda_{1,i-n+1}e_i\wedge
e_1\\
&\qquad-\sum_{i=n+1}^{d+n-1}\sum_{j=d+1}^{\infty}(2n-i-1)\lambda_{1j}\lambda_{d+1,i-n+1}e_i\wedge e_j-
\lambda_{d+1,1}\sum_{j=d+1}^{\infty}(n-1)\lambda_{1j}e_n\wedge e_j\\
&\qquad+\sum_{i=d+1}^{\infty}\sum_{j=n+1}^{d+n-1}(2n-j-1)\lambda_{1i}\lambda_{d+1,j-n+1}e_i\wedge e_j+
\lambda_{d+1,1}\sum_{i=d+1}^{\infty}(n-1)\lambda_{1i}e_i\wedge e_n\\
&\qquad-\sum_{i=2}^{d}\sum_{j=d+n}^{\infty}(2n-j-1)\lambda_{1,j-n+1}\lambda_{d+1,i}e_i\wedge e_j-\\
&\qquad-\lambda_{d+1,1}\sum_{j=d+n}^{\infty}(2n-j-1)\lambda_{1,j-n+1}e_1\wedge
e_j\Biggr\}\\
=&\frac
{1}{\lambda_{1,d+1}}\Biggl\{2\sum_{i=d+n}^{\infty}\sum_{j=2}^{d}(2n-i-1)\lambda_{1,i-n+1}\lambda_{d+1,j}e_i\wedge e_j+\\
&\qquad+2\sum_{i=d+1}^{\infty}\sum_{j=n+1}^{d+n-1}(2n-j-1)\lambda_{1i}\lambda_{d+1,j-n+1}e_i\wedge e_j\\
&\qquad+2\lambda_{1,d+1}\sum_{i=d+n}^{\infty}(2n-i-1)\lambda_{1,i-n+1}e_1\wedge
e_i-
2\lambda_{1,d+1}\sum_{i=d+1}^{\infty}(n-1)\lambda_{1i}e_i\wedge e_n\Biggr\}\\
=&\lambda_{1,d+1}\Biggl\{\frac
{2}{d-1}\sum_{i=d+n}^{\infty}\sum_{j=2}^{d}(2n-i-1)\lambda^{i+j-(n+d+1)}e_i\wedge e_j+\\
&\qquad+\frac
{2}{d-1}\sum_{i=d+1}^{\infty}\sum_{j=n+1}^{d+n-1}(2n-j-1)\lambda^{i+j-(n+d+1)}e_i\wedge e_j\\
&\qquad+2\sum_{i=d+n}^{\infty}(2n-i-1)\lambda^{i-(n+d)}e_1\wedge e_i-
2\sum_{i=d+1}^{\infty}(n-1)\lambda^{i-(d+1)}e_i\wedge e_n\Biggr\},\endalign$$
where we used formulae (5.41) to obtain the last equality. Hence, after
normalizing by the factor $\lambda_{1,d+1}\ne 0$ and shifting
indices by 1 we obtain (5.40).\qed\enddemo

\remark{Remark} One can show directly that $\alpha_{d,\lambda}$ satisfies the
co-Jacobi identity. The r.h.s\. of (5.40) is understood
as an element of the completed tensor product $\Cal G_{\infty}\widehat
{\otimes}\Cal G_{\infty}$ [Di].\endremark
}
\head 6. The group $G_{0\infty}$ and Poisson-Lie structures on it\endhead
\par In this section we study the group $G_{0\infty}$ of which $G_{\infty}$ is
a subgroup. We classify all Poisson-Lie
structures on $G_{0\infty}$ corresponding to coboundary Lie-bialgebra
structures on the Lie algebra $\Cal G_{0\infty}$ of  $G_{0\infty}$.
\par Let $X=\{x_i\}_{i\in\zZ_+}$ be a countable set of indeterminates. Let
$k\bigl[\bigl[ X\bigr]\bigr]$ be the ring of formal power series over $X$
whitout constant term with the
standard multiplication. Here $k$ is a commutative field assumed to be of
characteristic zero.
Let $Y=\{y_i\}_{i\in\zZ_+}$ be a second set of indeterminates, and $k[[Y]]$ be
the corresponding ring of formal power series over $Y$.
Consider the formal group $G_{0\infty}$ defined by a formal group law
$F=(F_i)_{i\in\zZ_+}$ [Se,Di] in countably infinite number of variables, where
$F_i\in k[[X,Y]]$ for every $i\in\zZ_+$, induced by a substitution of formal
power series in one variable.
Let $x(u)=\sum_{i=0}^{\infty}x_iu^i\in k[[X]][[u]]$ and
$y(u)=\sum_{i=0}^{\infty}y_iu^i\in k[[Y]][[u]]$ be elements in the rings of
formal power series with a constant term in the variable $u$ over the rings
$k[[X]]$ and $k[[Y]]$ respectively.
The multiplication of formal power series in the variable $u$ is defined again
as the substitution:
$$\eqalign{(xy)(u)=x(y(u))&=\sum_{i=0}^{\infty}x_i(y(u))^i\cr
&=\sum_{i=0}^{\infty}x_i\left[y_0^i+\sum_{j=1}^{\infty}(\sum_{\bigl(\sum_{\alpha =1}^{i}s_{\alpha}\bigr)=j}y_{s_1}\ldots y_{s_{i}})u^j\right]\cr
&=\sum_{i=0}^{\infty}x_iy_0^i+\sum_{i=1}^{\infty}x_i\sum_{j=1}^{\infty}\left(\sum_{\bigl(\sum_{\alpha =1}^{i}s_{\alpha}\bigr)=j}y_{s_1}\ldots y_{s_{i}}\right)u^j\cr
&=\sum_{i=0}^{\infty}x_iy_0^i+\sum_{j=1}^{\infty}\left(\sum_{i=1}^{\infty}x_i\sum_{\bigl(\sum_{\alpha =1}^{i}s_{\alpha}\bigr)=j}y_{s_1}\ldots y_{s_{i}}\right)u^j.}\tag 6.1$$
Therefore from (6.1) we obtain
$$\aligned F_0(X,Y)=&\sum_{i=0}^{\infty}x_iy_0^i,\\
F_j(X,Y)=&\sum_{i=1}^{\infty}x_i\sum_{\bigl(\sum_{\alpha
=1}^{i}s_{\alpha}\bigr)=j}y_{s_1}\ldots y_{s_{i}},\ \ \ \text{for every}\ \
j\ge 1.\endaligned\tag 6.1.1$$
This is a model of the group of diffeomorphisms of $\zR^1$ not necessarily
leaving the point $u=0$ fixed.
The identity here is $e=(0,1,0,0,\ldots )$. Formulae (6.1.1) have the
following interpretation. The ring $k[[X]]$ is naturally graded.
Namely, let us introduce a degree $|\  |\: X\to \Bbb Z_{+}$ defined on
the generators by $|x_i|:=i$. We extend it to monomials as
$|x_{i_1}\ldots x_{i_n}|=i_1+\cdots +i_n$. The grading on $k[[X]]$ and
$k[[Y]]$ induces a grading on $k[[X,Y]]$ in an obvious way. Then
$F_i(X,Y)=\sum_n f_n(X,Y)$, where each $f_n(X,Y)$ is a finite linear
combination of monomials of degree $n$. Clearly $G_{\infty}$, if viewed as a
formal group, will be a subgroup of $G_{0\infty}$. We define a Poisson
structure $\omega_x=\omega_{ij}(x)\pd {}{x_i}\wedge\pd {}{x_j}$ on the group
$G_{0\infty}$ as a bi-derivation $\omega_x\:k\bigl[\bigl[ X\bigr]\bigr]\otimes
k\bigl[\bigl[ X\bigr]\bigr]\to k\bigl[\bigl[ X\bigr]\bigr]$, where
$\omega_{ij}(x)=-\omega_{ji}(x)\in k\bigl[\bigl[ X\bigr]\bigr]$, satisfying the
Jacobi identity. The methods developed in analyzing the Poisson-Lie structures
on $G_{\infty}$ apply without major
changes to the case of $G_{0\infty}$, but with two important
differences. Namely, Theorem 4.1 still holds with ${\Omega}(u,v;{\Cal X})$
defined as
${\Omega}(u,v;{\Cal X}):=\sum_{i,j=0}^{\infty}\omega_{ij}u^iv^j,$
but in the solution of the cocycle equation
$${\Omega}(u,v;{x})=\varphi(u,v){x}^{\prime}(u){x}^{\prime}(v)-\varphi({x}(u),{x}(v)) \tag 6.2$$
$\varphi (u,v)$ does not have to be divisible by $uv$. Thus, this condition is
dropped. This change affects the analysis of the equation
$$\varphi (u,v)\left[\partial_{u}\varphi (w,u)+\partial_{v}\varphi
(w,v)\right]+c.p.=0, \tag 6.3$$
as well as the structure of its solutions.
\par Still, we shall show that the Poisson-Lie structures on $G_{0\infty}$ fall
into two main classes.
Also, as we found in the previous section, the solution (6.2) corresponds to a
cocycle that is a coboundary in the Lie algebra $\Cal G_{\infty}$. This result
carries over to the case of $\Cal G_{0\infty}$ without change. We
showed that all cocycles on $\Cal G_{0\infty}$ are coboundaries. This
allows us to completely classify them. That is, we shall give a
classification of all
$r$-matrices on $\Cal G_{0\infty}$.
\par  We start with two special solutions of (6.3). These two solutions are
only two ``points'' in the otherwise infinite-parameter space of solutions of
(6.3).
\proclaim{\bf Theorem 6.1} The functions\newline
(i)  $\varphi (u,v)=u-v$, and\newline
(ii) $\varphi (u,v)=e^{\lambda u}-e^{\lambda v}$, where $\lambda$ is an
arbitrary parameter,\newline
are solutons of {\rm (6.3)}, thus giving rise to two Poisson-Lie structures
on $G_{0\infty}$.\newline
These are the only solutions of {\rm (6.3)} of the form $\varphi
(u,v)=a(u)-a(v)$.\endproclaim
\par The proof is straightforward and we omit it. Here, we only write the
Poisson brackets in coordinates. For the case (i) we have
$$\omega_{ij}(x)=i(j+1)x_ix_{j+1}-(i+1)jx_{i+1}x_{j}-x_i\delta_{j}^{0}+x_j\delta_{i}^{0}, \qquad i,j\in \zZ_+.$$
Notice, that there are no terms higher than quadratic in the right hand side.
For the case (ii) we obtain
$$\eqalign{\omega_{ij}(x)&=(j+1)x_{j+1}\sum_{p=0}^{i+1}\frac
{px_p}{(i-p+1)!}-(i+1)x_{i+1}\sum_{q=0}^{j+1}\frac {qx_q}{(j-q+1)!}\cr
&\qquad-\delta_{j}^{0}\sum_{p=0}^{i}\frac {1}{p!}\sum_{r_1+\ldots
+r_p=i}x_{r_1}\ldots x_{r_p}+\delta_{i}^{0}\sum_{q=0}^{j}\frac
{1}{q!}\sum_{r_1+\ldots +r_q=j}x_{r_1}\ldots x_{r_q}.}$$
\par Now, we proceed with the main result of this section.
\proclaim{Theorem 6.2} All solutions of {\rm (6.3)} fall into the following two
classes\newline
(a) The first class is given by Theorem {\rm 4.4.}\newline
(b) The second class is given by
$$\varphi (u,v)=\frac {1}{\lambda_{01}}\biggl[f(u)g(v)-f(v)g(u)\biggr]$$
where $\lambda_{01}\ne 0$, for any formal power series $f(u)$ and $g(u)$
satisfying the relation
$$f'(u)g(u)-f(u)g'(u)=\lambda_{01}g(u)-2\lambda_{02}f(u).$$
Here, $\lambda_{01}$ and $\lambda_{02}$ are arbitrary parameters with
$\lambda_{01}$ being subject to the above restriction: $\lambda_{01}\ne
0$.\endproclaim
\demo{Proof} We look again for a solution of (6.3) in a form of a formal power
series $\varphi (u,v)=\sum_{m,n=0}^{\infty}\lambda_{mn}u^mv^n$, where
$\lambda_{mn}=-\lambda_{nm}$. Substituting into (6.3) we obtain (in a similar
way as in the proof of Theorem 4.5)
$$\sum_{s=0}^{k+1}(k-2s+1)\lambda_{n,k-s+1}\lambda_{rs}-\sum_{s=0}^{n+1}(n-2s+1)\lambda_{k,n-s+1}\lambda_{rs} -\sum_{s=0}^{r+1}(r-2s+1)\lambda_{n,r-s+1}\lambda_{ks}=0 \tag 6.4$$
where $k<n<r$.
\par Let $k=0$. Then (6.4) becomes
$$\lambda_{n1}\lambda_{r0}-\lambda_{n0}\lambda_{r1}-\sum_{s=0}^{n+1}(n-2s+1)\lambda_{0,n-s+1}\lambda_{rs}-\sum_{s=0}^{r+1}(r-2s+1)\lambda_{n,r-s+1}\lambda_{0s}=0. \tag 6.4*$$
Notice that the summation in the sums above is to $n$ and $r$ respectively. If
we now let $n=1$, we finally obtain
$$-\lambda_{01}\lambda_{1r}+2\lambda_{02}\lambda_{0r}=\sum_{s=0}^{r}(r-2s+1)\lambda_{1,r-s+1}\lambda_{0s}=-\sum_{s=0}^{r}(r-2s+1)\lambda_{0,r-s+1}\lambda_{1s}.$$
Therefore,
$$\sum_{s=0}^{r}(r-s+1)\lambda_{0,r-s+1}\lambda_{1s}-\sum_{s=0}^{r}s\lambda_{0,r-s+1}\lambda_{1s}=\lambda_{01}\lambda_{1r}-2\lambda_{02}\lambda_{0r}. \tag 6.5$$
Let us define the functions $f(u)=\sum_{n=0}^{\infty}\lambda_{0n}u^n$ and
$g(u)=\sum_{n=0}^{\infty}\lambda_{1n}u^n$. Then, after multiplying both sides
of (6.5) by $u^r$, and summing over $r$ we obtain
$$\sum_{r=0}^{\infty}\sum_{s=0}^{r}(r-s+1)\lambda_{0,r-s+1}\lambda_{1s}u^r-\sum_{r=0}^{\infty}\sum_{s=0}^{r}s\lambda_{0,r-s+1}\lambda_{1s}u^r=\lambda_{01}g(u)-2\lambda_{02}f(u),$$
which, after the change of variables $r=q+s-1$, is equivalent to
$$\sum_{s=0}^{\infty}\sum_{q=0}^{\infty}q\lambda_{0q}\lambda_{1s}u^{q+s-1}-\sum_{s=0}^{\infty}\sum_{q=0}^{\infty}s\lambda_{0q}\lambda_{1s}u^{q+s-1}=\lambda_{01}g(u)-2\lambda_{02}f(u).$$
Thus, the functions $f(u)$, and $g(u)$ satisfy
$f'(u)g(u)-f(u)g'(u)=\lambda_{01}g(u)-2\lambda_{02}f(u).$\qed\enddemo
\par Next, we split the remaining part of the proof into five lemmas.
\proclaim{Lemma 6.3} For every $n\ge 1$ the parameters $\lambda_{1n}$ are
rational functions of $\lambda_{0n}$, provided that $\lambda_{01}\ne
0$.\endproclaim
\remark{Remark} The case of $\lambda_{01}=0$ will be treated in Lemmas 6.5,
6.6, and 6.7.\endremark
\demo{Proof} From (6.5) we have
$$\sum_{s=0}^{r-1}(r-2s+1)\lambda_{0,r-s+1}\lambda_{1s}-(r-1)\lambda_{01}\lambda_{1r}=\lambda_{01}\lambda_{1r}-2\lambda_{02}\lambda_{0r},$$
thus
$$\lambda_{1r}=\frac
{1}{r\lambda_{01}}\left[\sum_{s=0}^{r-1}(r-2s+1)\lambda_{0,r-s+1}\lambda_{1s}+2\lambda_{02}\lambda_{0r}\right].$$
This gives us a recursive relation for $\lambda_{1r}$. Here are the first
several $\lambda_{1r}$'s
$$\eqalign{&\lambda_{12}=\frac
{1}{\lambda_{01}}\biggl[3\lambda_{03}\lambda_{10}+2(\lambda_{02})^2\biggr]=\frac {1}{2}\left[2\frac {(\lambda_{02})^2}{\lambda_{01}}-3\lambda_{03}\right]\cr
&\lambda_{13}=\frac {1}{3}\left[2\frac
{\lambda_{02}\lambda_{03}}{\lambda_{01}}-4\lambda_{04}\right]\cr
&\lambda_{14}=\frac {1}{24}\left[2\frac
{(\lambda_{02})^2\lambda_{03}}{(\lambda_{01})^2}-9\frac
{(\lambda_{03})^2}{\lambda_{01}}+20\frac
{\lambda_{02}\lambda_{04}}{\lambda_{01}}-30\lambda_{05}\right]\cr
&\vdots}$$
Therefore, the claim is proved by induction.\qed\enddemo
\proclaim{Lemma 6.4} If $\lambda_{01}\ne 0$, then we have the following
formula:
$$\lambda_{nr}=\frac
{1}{\lambda_{01}}\biggl[\lambda_{0n}\lambda_{1r}-\lambda_{1n}\lambda_{0r}\biggr],\qquad \text{for all}\quad {n,r\ge 0}.$$\endproclaim
\demo{Proof} From (6.4*) with $n=2$ we obtain
$$\lambda_{12}\lambda_{0r}-\lambda_{01}\lambda_{2r}+3\lambda_{03}\lambda_{0r}-\sum_{s=0}^{r+1}(r-2s+1)\lambda_{2,r-s+1}\lambda_{0s}=0.\tag 6.6$$
If we now let $r=3$, we obtain a formula for $\lambda_{23}$:
$$\lambda_{12}\lambda_{03}-\lambda_{01}\lambda_{23}+3(\lambda_{03})^2-2\lambda_{23}\lambda_{01}-2\lambda_{12}\lambda_{03}-4\lambda_{02}\lambda_{04}=0,$$
from which
$$\eqalign{\lambda_{23}&=\frac {1}{\lambda_{01}}\left[-\frac {1}{3}\frac
{(\lambda_{02})^2}{\lambda_{01}}\lambda_{03}+\frac {2}{3}(\lambda_{03})^2-\frac
{4}{3}\lambda_{02}\lambda_{04}\right]\cr
&=\frac {1}{\lambda_{01}}\biggl[\frac {1}{3}\biggl(2\frac
{\lambda_{02}\lambda_{03}}{\lambda_{01}}-4\lambda_{04}\biggr)\lambda_{02}-\frac
{1}{2}\biggl(2\frac
{(\lambda_{02})^2}{\lambda_{01}}-3\lambda_{03}\biggr)\lambda_{03}\biggr]\cr
&=\frac
{1}{\lambda_{01}}\biggl[\lambda_{02}\lambda_{13}-\lambda_{12}\lambda_{03}\biggr].}$$
\par Now, let us assume that $\lambda_{2k}=\frac
{1}{\lambda_{01}}\left[\lambda_{02}\lambda_{1k}-\lambda_{12}\lambda_{0k}\right]$ for $1\le k\le r-1$. We need to prove that it is true for $k=r$.
 From (6.6) we have
$$\lambda_{12}\lambda_{0r}-\lambda_{01}\lambda_{2r}+3\lambda_{03}\lambda_{0r}+(1-r)\lambda_{01}\lambda_{2r}+\sum_{s=0}^{r-1}(r-2s+1)\lambda_{0,r-s+1}\lambda_{2s}=0.$$
Using the induction hypothesis we transform the above equation into
$$\lambda_{12}\lambda_{0r}+3\lambda_{03}\lambda_{0r}+\frac
{1}{\lambda_{01}}\sum_{s=0}^{r-1}(r-2s+1)\lambda_{0,r-s+1}\left[\lambda_{02}\lambda_{1s}-\lambda_{12}\lambda_{0s}\right]=r\lambda_{01}\lambda_{2r}.$$
Next, using the fact that
$\sum_{s=2}^{r-1}(r-2s+1)\lambda_{0,r-s+1}\lambda_{0s}=0$ we transform the
above equality further into
$$(2-r)\lambda_{12}\lambda_{0r}+3\lambda_{03}\lambda_{0r}+\frac
{\lambda_{02}}{\lambda_{01}}\sum_{s=0}^{r-1}(r-2s+1)\lambda_{0,r-s+1}\lambda_{1s}=r\lambda_{01}\lambda_{2r}.$$
The final step is to use (6.5) and the formula for $\lambda_{12}$ obtained in
the proof of Lemma 6.3. This leads to
$$\biggl(2\frac
{(\lambda_{02})^2}{\lambda_{01}}-3\lambda_{03}\biggr)\lambda_{0r}-r\lambda_{12}\lambda_{0r}+3\lambda_{03}\lambda_{0r}+\frac {\lambda_{02}}{\lambda_{01}}\biggl[r\lambda_{01}\lambda_{1r}-2\lambda_{02}\lambda_{0r}\biggr]=r\lambda_{01}\lambda_{2r}.$$
After collecting all terms we obtain
$$\lambda_{2r}=\frac
{1}{\lambda_{01}}\biggl[\lambda_{02}\lambda_{1r}-\lambda_{12}\lambda_{0r}\biggr].$$
\par Finally, we assume that $\lambda_{nk}=\frac
{1}{\lambda_{01}}\left[\lambda_{0n}\lambda_{1k}-\lambda_{1n}\lambda_{0k}\right]$ for $1\le k\le r-1$ and each $n\ge 1$. Then, we show that it is true for $k=r$. The steps are essentially the same as in the previous calculation. From (6.4*) we have
$$\lambda_{1n}\lambda_{0r}-\lambda_{0n}\lambda_{1r}+\sum_{s=0}^{n-1}(n-2s+1)\lambda_{0,n-s+1}\lambda_{sr}+\sum_{s=0}^{r-1}(r-2s+1)\lambda_{0,r-s+1}\lambda_{ns}=(r+n-2)\lambda_{01}\lambda_{nr}.$$
Applying the induction hypothesis, the above equation transforms to
$$(2-r)\lambda_{1n}\lambda_{0r}+(n-2)\lambda_{1r}\lambda_{0n}-\frac
{\lambda_{0r}}{\lambda_{01}}\sum_{s=0}^{n-1}(n-2s+1)\lambda_{0,n-s+1}\lambda_{1s}+\frac {\lambda_{0n}}{\lambda_{01}}\sum_{s=0}^{r-1}(r-2s+1)\lambda_{0,r-s+1}\lambda_{1s}=$$
$$=(r+n-2)\lambda_{01}\lambda_{nr}.$$
Now, we use (6.5) to obtain
$$(2-r)\lambda_{1n}\lambda_{0r}+(n-2)\lambda_{1r}\lambda_{0n}-\frac
{\lambda_{0r}}{\lambda_{01}}\biggl[n\lambda_{01}\lambda_{1n}-2\lambda_{02}\lambda_{0n}\biggr]+\frac {\lambda_{0n}}{\lambda_{01}}\biggl[r\lambda_{01}\lambda_{1r}-2\lambda_{02}\lambda_{0r}\biggr]=$$
$$=(r+n-2)\lambda_{01}\lambda_{nr},$$
and collecting terms we have
$$\lambda_{nr}=\frac
{1}{\lambda_{01}}\biggl[\lambda_{0n}\lambda_{1r}-\lambda_{1n}\lambda_{0r}\biggr].$$
This completes the proof of the lemma.\qed\enddemo
\par As a consequence,
$$\varphi (u,v)=\sum_{n,m=0}^{\infty}\lambda_{nm}u^nv^m=\frac
{1}{\lambda_{01}}\biggl[f(u)g(v)-f(v)g(u)\biggr],$$
where $f(u)=\sum_{s=0}^{\infty}\lambda_{0s}u^s$ and
$g(u)=\sum_{s=0}^{\infty}\lambda_{1s}u^s.$
\proclaim{Lemma 6.5} If $\lambda_{01}=0$, then $\lambda_{02}=0$. If in addition
we assume that $\lambda_{12}\ne 0$, then it follows that $\lambda_{0n}=0$,
$\forall\ {n\ge 2}$.\endproclaim
\demo{Proof} Again, from (6.5) with $\lambda_{01}=0$ we obtain
$$\sum_{s=1}^{r-1}(r-2s+1)\lambda_{0,r-s+1}\lambda_{1s}+2\lambda_{02}\lambda_{0r}=0,\qquad\text{for
all}\quad{r\ge 2}.$$
For $r=2$ we have $(\lambda_{02})^2=0$. Using this fact, the above equation
reduces to
$$\sum_{s=1}^{r-1}(r-2s+1)\lambda_{0,r-s+1}\lambda_{1s}=0,\qquad\text{for
all}\quad {r\ge 3}. \tag 6.7$$
For $r=3$ we have an identity. For $r=4$ we obtain
$\lambda_{03}\lambda_{12}-\lambda_{02}\lambda_{13}=\lambda_{03}\lambda_{12}=0\implies \lambda_{03}=0.$
Therefore, if we assume that $0=\lambda_{01}=\lambda_{02}=\ldots
=\lambda_{0,r-2}$, equation (6.7) gives
$(r-3)\lambda_{0,r-1}\lambda_{12}=0\implies \lambda_{0,r-1}=0$,for all ${r>3}$.
As a result the system of equations (6.4) reduces to the system (4.6). But for
$\lambda_{12}\ne 0$ the solution
of that system of equations is the one corresponding to $d=1$ as described by
Theorem 4.5.\qed\enddemo
\proclaim{Lemma 6.6} Suppose that $\lambda_{01}=0=\lambda_{02}$. If
$\lambda_{12}=\ldots =\lambda_{1n}=0$ for some $n\ge 2$, then
$\lambda_{03}=\ldots =\lambda_{0,n+1}=0$.\endproclaim
\demo{Proof} From (6.4) with $k=0$ and $r=n+1$ we obtain
$$\multline\sum_{s=1}^{n}(n-2s+1)\lambda_{n+1,n-s+1}\lambda_{0s}+\sum_{s=1}^{n+1}(n-2s+2)\lambda_{0,n-s+2}\lambda_{ns}=\\
-\lambda_{0,n+1}\lambda_{1n}+\lambda_{0n}\lambda_{1,n+1}-(n+1)(\lambda_{0,n+1})^2+(n+2)\lambda_{0,n+2}\lambda_{0n}\endmultline \tag 6.8$$
For $n=2$ the above equation gives
$3(\lambda_{03})^2=\lambda_{03}\lambda_{12}$. Thus if $\lambda_{12}=0$ then
$\lambda_{03}=0$.
\par Assume that $\lambda_{12}=\ldots =\lambda_{1,n-1}=0$ and
$\lambda_{03}=\ldots =\lambda_{0n}=0$. Then from (6.8) we obtain
$(n+1)(\lambda_{0,n+1})^2=-(n-1)\lambda_{0,n+1}\lambda_{1n}.$
Therefore if $\lambda_{1n}=0$, then $\lambda_{0,n+1}=0$.\qed\enddemo
\proclaim{Lemma 6.7} Suppose that $\lambda_{01}=0$, and $\lambda_{12}=\ldots
=\lambda_{1d}=0$, and $\lambda_{1,d+1}\ne 0$ for some $d\in \zN$. Then
$\lambda_{0n}=0$ for every $n\ge 2$.\endproclaim
\demo{Proof} From Lemma 6.6 it follows that $\lambda_{02}=\ldots
=\lambda_{0,d+1}=0$. Further, from (6.7) we have
$$\sum_{s=d+1}^{n}(n-2s+1)\lambda_{0,n-s+1}\lambda_{1s}=0. \tag 6.9$$
Let $n=2d+2$. Then
$$\sum_{s=d+1}^{2d+2}(2d-2s+3)\lambda_{0,2d-s+3}\lambda_{1s}=\lambda_{0,d+2}\lambda_{1,d+1}-\lambda_{0,d+1}\lambda_{1,d+2}=0.$$
But $\lambda_{0,d+1}=0$ according to Lemma 6.5. Therefore
$\lambda_{0,d+2}\lambda_{1,d+1}=0$. Since $\lambda_{1,d+1}\ne 0$ then
$\lambda_{0,d+2}=0$. Now, from (6.9) with $n=2d+3$ we obtain
$\sum_{s=d+1}^{2d+3}(2d-2s+4)\lambda_{0,2d-s+4}\lambda_{1s}=2\lambda_{0,d+3}\lambda_{1,d+1}=0,$
therefore $\lambda_{0,d+3}=0$. By induction the statetment follows.\qed\enddemo
\par The above arguments show that if $\lambda_{01}=0$ the classification of
the solutions of (6.3) reduces to the classification of all those solutions
$\varphi (u,v)$ of (6.3) which are divisible by $uv$. But all such solutions
have been given by Theorem 4.5, and the proof of Theorem 6.2 is finished.\qed
\par Next, we discuss an interesting particular case which follows from the
general formulae.
First we write the Lie bialgebra structures on $\Cal G_{0\infty}$ that
correspond to the Poisson-Lie structures (i) and (ii) of Theorem 6.1. Namely,
if $e_n\in \Cal G_{0\infty}$ is a basis element we have
$$\eqalign{\alpha(e_n)&=\alpha_{ij}^{n}e_i\wedge e_j\cr
&=\left[(\delta_{i-n+1}^{1}\delta_{j+1}^{0}-\delta_{i-n+1}^{0}\delta_{j+1}^{0})(2n-i)+(\delta_{i+1}^{1}\delta_{j-n+1}^{0}-\delta_{i+1}^{0}\delta_{j-n+1}^{1})(2n-j)\right]e_i\wedge e_j\cr
&=-2ne_{-1}\wedge e_{n}+2(n+1)e_{0}\wedge e_{n-1}}.$$
Thus, for the case (i) we obtain
$$\alpha(e_n)=-2ne_{-1}\wedge e_{n}+2(n+1)e_{0}\wedge e_{n-1}.\tag 6.10$$
Similarly, for the case (ii)
$$\eqalign{\alpha(e_n)=&\biggl[\frac {1}{(i-n+1)!}\delta_{j+1}^{0}-\frac
{1}{(j+1)!}\delta_{i-n+1}^{0}\biggr](2n-i)e_i\wedge e_j\cr
&\qquad+\biggl[\frac {1}{(i+1)!}\delta _{j-n+1}^{0}-\frac
{1}{(j-n+1)!}\delta_{i+1}^{0}\biggr](2n-j)e_i\wedge e_j\cr
=&-2\sum_{j=n-1}^{\infty}\frac {(2n-j)}{(j-n+1)!}e_{-1}\wedge
e_{j}+2(n+1)\sum_{i=-1}^{\infty}\frac {1}{(i+1)!}e_i\wedge e_{n-1}.}$$
\par The commutator for $\Cal G_{0\infty}$ has the standard form
$\left[e_n,e_m\right]=(n-m)e_{n+m}, \hskip 15pt (n,m\ge -1).$
That is, $\Cal G_{0\infty}$ is the Witt algebra.
\par It is well known that the Witt algebra contains $sl_2$ as a Lie
subalgebra. In our notation the defining relations are given by
$$\left[e_1,e_{-1}\right]=2e_0,\qquad\left[e_1,e_0\right]=e_1,\qquad\left[e_0,e_{-1}\right]=e_{-1}.$$
If we now turn to the formula (6.10)
$$\alpha(e_n)=-2ne_{-1}\wedge e_n+2(n+1)e_0\wedge e_{n-1},$$
it gives us the following Lie bialgebra structure on $sl_2$:
$$\alpha(e_{-1})=0,\qquad\alpha(e_0)=2e_0\wedge
e_{-1},\qquad\alpha(e_1)=-2e_{-1}\wedge e_1.$$
On the other hand if we take formula (5.15)
$$\alpha(e_n)=-2ne_d\wedge e_n+2(n-d)e_0\wedge e_{d+n},$$
and consider it for $d=1$, that is
$$\alpha(e_n)=-2ne_1\wedge e_n+2(n-1)e_0\wedge e_{n+1}.$$
It gives a second Lie bialgebra structure on $sl_2$, namely
$$\alpha(e_{-1})=2e_1\wedge e_{-1},\qquad\alpha(e_0)=-2e_0\wedge
e_1,\qquad\alpha(e_1)=0.$$
However, it turns out that these two structures are isomorphic
under the action of the group of automorphisms of the Witt algebra.
\head 7. Elements of Representation Theory\endhead
\par Let $G$ be a Poisson-Lie group and $\overline{\omega}$ be a Poisson-Lie
structure on $G$.
Let $V$ be a space on which $G$ acts, i.e., there is a map $G\times V\to V$.
Such a space is called a $G$-space. Assume that $V$ is equipped with
a Poisson structure $\omega$. Recall the following definition [STS1].
\definition{Definition 7.1} The action of $G$ on $V$ is called Poisson if the
map $G\times V\to V$ is Poisson.
Here $G\times V$ is equipped with the product Poisson structure.\enddefinition
\par In this section we study the following problem. Suppose that we are given
the Poisson-Lie group $G_{\infty}$.
Consider the space $V_{\lambda}=\{x(u)(du)^{\lambda}\mid
x(u)=\sum_{i=0}^{\infty}x_iu^i\}$, $\lambda\in \zR$.
The space $V_{\lambda}$ is sometimes referred to as the space of
$\lambda$-densities (Jacobians) over the real line.
The group $G_{\infty}$ acts naturally on $V_{\lambda}$. Let $y\in G_{\infty}$
and $x(u)(du)^{\lambda}\in V_{\lambda}$.
 Then the action of
$G_{\infty}$ on $V_{\lambda}$ is defined by
$$x(u)(du)^{\lambda}\mapsto
x\left(y(u)\right)\left(y'(u)\right)^{\lambda}(du)^{\lambda},$$
 where $y(u)=\sum_{i=1}^{\infty}y_iu^{i}$, and
$$\align
\left(y'(u)\right)^{\lambda}&=\left(\sum_{i=1}^{\infty}iy_iu^{i-1}\right)^{\lambda}\\
&=\left(y_1+\sum_{i=2}^{\infty}iy_iu^{i-1}\right)^{\lambda}\\
&=y_1^{\lambda}\left(1+\sum_{i=2}^{\infty}i\frac
{y_i}{y_1}u^{i-1}\right)^{\lambda}\\
&=y_1^{\lambda}\left[1+\frac {\lambda}{1\text{!}}\sum_{i=2}^{\infty}i\frac
{y_i}{y_1}u^{i-1}+\frac
{\lambda(\lambda-1)}{2\text{!}}\left(\sum_{i=2}^{\infty}i\frac
{y_i}{y_1}u^{i-1}\right)^2+\ldots\right].\endalign$$
\par A question arises. Are there Poisson structures on the space $V_{\lambda}$
such that the above action of
$G_{\infty}$ on $V_{\lambda}$ is a Poisson action?
In other words, is there a Poisson strucuture $\omega$ on $V_{\lambda}$ such
that the map $G_{\infty}\times V_{\lambda}\to V_{\lambda}$ is
 a Poisson map?
Here again $G_{\infty}\times V_{\lambda}$ is equipped with the product Poisson
structure.
\par Let $y(u)=\sum_{i=1}^{\infty}y_iu^i\in G_{\infty}$, and
$x(u)(du)^{\lambda}\in V_{\lambda}$.
Let us define
$z_{\lambda}(u):=x(y(u))\bigl[y'(u)\bigr]^{\lambda}=\sum_{i=0}^{\infty}z_iu^i,$
where $z_i=z_i(x,y;\lambda)$ are the coordinates of $z_{\lambda}$.  If we also
introduce the notation
$J(u):= y'(u)=\sum_{i=1}^{\infty}iy_iu^{i-1}$, we have
$$\eqalign{x(u)(du)^{\lambda}\mapsto x(y(u))(y'(u))^{\lambda}(du)^{\lambda}
&=x(y(u))(J(u))^{\lambda}(du)^{\lambda}\cr
&=z_{\lambda}(u)(du)^{\lambda}\cr
&=\sum_{i=0}^{\infty}x_i(y(u))^i\left(\sum_{i=1}^{\infty}iy_iu^{i-1}\right)^{\lambda}(du)^{\lambda}.}$$
 \par Defining $z(u):=x(y(u))$ and using the definition
$z_{\lambda}(u)=x(y(u))\bigl[J(u)\bigr]^{\lambda}$
we deduce that
$z'_{\lambda}(u)=z'(u)(J(u))^{\lambda}+z(u)\left[(J(u))^{\lambda}\right]',$
where $\prime$ stands for the derivative with respect to $u$.
 \par An argument analogous to the argument given in Sec\. 1 implies that
the map $G_{\infty}\times V_{\lambda}\to V_{\lambda}$ is Poisson if and only if
$$\omega _{ij}(z)=\omega _{kl}(x)\pd {z_i}{x_k} \pd {z_j}{x_l} +
\overline{\omega}_{kl}{(y)}\pd {z_i}{y_k} \pd {z_j}{y_l}.
\tag 7.2$$
Here  $\omega _{ij}(x)=\{x_i,x_j\}$ and $\overline{\omega}
_{ij}(y)=\{y_i,y_j\}$ where $\{x_i\}_{i\in\zZ_+}$ and $\{y_i\}_{i\in\zZ_+}$
are the coordinates on $V_{\lambda}$ and $G_{\infty}$ respectively.
 Also, let us introduce in a manner similar to the
one used in Sec\. 4 a generating series for the Poisson strucutures on
$V_{\lambda}$ as $\Omega(u,v;\Cal X):= \sum_{i,j=0}^{\infty}\omega_{ij}u^iv^j.$
\proclaim{Lemma 7.2} The multiplicativity condition {\rm (7.2)} is equivalent
to the following functional equation
$$\eqalign{\Omega(u,v;z_{\lambda})=&\Omega(y(u),y(v);x)(J(u))^{\lambda}(J(v))^{\lambda}+
\overline{\Omega}(u,v;y)z'(u)(J(u))^{\lambda-1}z'(v)(J(v))^{\lambda-1}+\cr
&\qquad+\lambda\partial_u\overline{\Omega}(u,v;y)z(u)(J(u))^{\lambda-1}z'(v)(J(v))^{\lambda-1}+\cr
&\qquad+\lambda\partial_v\overline{\Omega}(u,v;y)z'(u)(J(u))^{\lambda-1}z(v)(J(v))^{\lambda-1}+\cr
&\qquad+\lambda^2\partial^2_{u,v}\overline{\Omega}(u,v;y)z(u)(J(u))^{\lambda-1}z(v)(J(v))^{\lambda-1}.}\tag 7.3$$
Here $\overline{\Omega}$ stands for the generating series of the Poisson-Lie
strucutures on $G_{\infty}$.\endproclaim
\demo{Proof} Multiplying both sides of equation (7.2) by $u^iv^j$, summing over
$i$ and $j$, and using the definition of
$\Omega$ we obtain
$$\Omega(u,v;z_{\lambda})=\omega _{kl}(x)\pd {z_{\lambda}}{x_k} \pd
{z_{\lambda}}{x_l}+
\overline{\omega}_{kl}(y)\pd {z_{\lambda}}{y_k} \pd {z_{\lambda}}{y_l}.\tag
7.4$$
\par On the other hand, the following formulae are valid
$$\align\pd {z_{\lambda}}{x_k}=&
\left(\sum_{j=1}^{\infty}y_ju^j\right)^k\left(\sum_{i=1}^{\infty}iy_iu^{i-1}\right)^{\lambda}=
\bigl[y(u)\bigr]^k\bigl[J(u)\bigr]^{\lambda},\\\\\allowdisplaybreak
\pd {z_{\lambda}}{y_k}=&
\left[\sum_{i=0}^{\infty}ix_i\left(\sum_{j=1}^{\infty}y_ju^j\right)^{i-1}u^k\right]\left(\sum_{i=1}^{\infty}iy_iu^{i-1}\right)^{\lambda}+\left[\sum_{i=0}^{\infty}x_i\left(\sum_{j=1}^{\infty}y_ju^j\right)^{i}\right]
\lambda\left(\sum_{i=1}^{\infty}iy_iu^{i-1}\right)^{\lambda-1}ku^{k-1}\\\\
=&x'(y(u))(J(u))^{\lambda}u^k+\lambda x(y(u))[J(u)]^{\lambda -1}ku^{k-1}\\
=&z'(u)[J(u)]^{\lambda-1}u^k+\lambda
z_{\lambda}(u)[J(u)]^{-1}ku^{k-1}.\endalign$$
Therefore equation (7.4) takes the form
$$\align\Omega(u,v;z_{\lambda})=&\Omega(y(u),y(v);x)(J(u))^{\lambda}(J(v))^{\lambda}+\\
&\qquad+\overline{\omega}_{kl}(y)\biggl[z'(u)(J(u))^{\lambda-1}u^k+\lambda
z(u)(J(u))^{\lambda-1}ku^{k-1}\biggr]\times\\
&\qquad\qquad\qquad\qquad\qquad\times\biggl[z'(v)(J(v))^{\lambda-1}v^l+\lambda
z(v)(J(v))^{\lambda-1}lv^{l-1}\biggr]\\\allowdisplaybreak
=&\Omega(y(u),y(v);x)(J(u))^{\lambda}(J(v))^{\lambda}+
\overline{\omega}_{kl}(y)z'(u)(J(u))^{\lambda-1}u^kz'(v)(J(v))^{\lambda-1}v^l+\\
&\qquad+\overline{\omega}_{kl}(y)\lambda
z(u)(J(u))^{\lambda-1}ku^{k-1}z'(v)(J(v))^{\lambda-1}v^l+\\
&\qquad+\overline{\omega}_{kl}(y)z'(u)(J(u))^{\lambda-1}u^k\lambda
z(v)(J(v))^{\lambda-1}lv^{l-1}+\\
&\qquad+\lambda^2\overline{\omega}_{kl}(y)z(u)(J(u))^{\lambda-1}ku^{k-1}z(v)(J(v))^{\lambda-1}lv^{l-1}\\\allowdisplaybreak
=&\Omega(y(u),y(v);x)(J(u))^{\lambda}(J(v))^{\lambda}+
\overline{\Omega}(u,v;y)z'(u)(J(u))^{\lambda-1}z'(v)(J(v))^{\lambda-1}+\\
&\qquad+\lambda\partial_u\overline{\Omega}(u,v;y)z(u)(J(u))^{\lambda-1}z'(v)(J(v))^{\lambda-1}+\\
&\qquad+\lambda\partial_v\overline{\Omega}(u,v;y)z'(u)(J(u))^{\lambda-1}z(v)(J(v))^{\lambda-1}+\\
&\qquad+\lambda^2\partial^2_{u,v}\overline{\Omega}(u,v;y)z(u)(J(u))^{\lambda-1}z(v)(J(v))^{\lambda-1}.\endalign$$
This concludes the proof of the Lemma.\qed\enddemo
\par In the above formulae $\overline{\Omega}$ is given by (cf\. Sec\. 4)
$$\overline{\Omega}(u,v;{y})=\varphi(u,v){y}^{\prime}(u){y}^{\prime}(v)
-\varphi({y}(u),{y}(v)), \tag 7.5$$
where the function $\varphi (u,v)$ satisfies the equation
$$\varphi (u,v)\left[\partial_{u}\varphi (w,u)+\partial_{v}\varphi
(w,v)\right]+c.p.=0. \tag 7.6$$
In other words $\varphi (u,v)$ is given by $\varphi
(u,v)=f(u)g(v)-f(v)g(u),$ (cf. Sec. 4.),
where the functions $f$ and $g$ satisfy the relation
$$f'(u)g(u)-f(u)g'(u)=\alpha f'(u)+\beta g'(u)\hskip 0.2in \biggl(\implies
f''(u)g(u)-f(u)g''(u)=\alpha f''(u)+\beta g''(u)\biggr).
\tag 7.7$$
Here, $\alpha$ and $\beta$ are arbitrary constants.
\par Now we will describe a class of solutions of equation (7.3).
\proclaim{Theorem 7.3} If $\varphi (u,v)$ is defined by the equation
{\rm (7.6)} and $\overline{\Omega}$ is defined by {\rm (7.5)} one has
the following solution of {\rm (7.3)}:
$$\eqalign{{\Omega}(u,v;x)=\varphi (u,v)&x'(u)x'(v)+\lambda\partial_u\varphi
(u,v)x(u)x'(v)+\cr
&+\lambda\partial_v\varphi (u,v)x'(u)x(v)+\lambda^2\partial^2_{u,v}\varphi
(u,v)x(u)x(v).} \tag 7.8$$
\endproclaim
\demo{Proof} Using formula (7.5) in the r.h.s\. of equation (7.4) we obtain
$$\allowdisplaybreaks\align
{r.h.s.}&=\Omega(y(u),y(v);x)(J(u))^{\lambda}(J(v))^{\lambda}+
\varphi (u,v)z'(u)z'(v)(J(u))^{\lambda}(J(v))^{\lambda}-\\
&\qquad-\varphi({y}(u),{y}(v))z'(u)z'(v)(J(u))^{\lambda-1}(J(v))^{\lambda-1}
+\lambda\partial_u\varphi (u,v)z_{\lambda}(u)
z'(v)(J(v))^{\lambda}+\\
&\qquad+\varphi (u,v)z(u)z'(v)\bigl[(J(u))^{\lambda}\bigr]'(J(v))^{\lambda}-
\lambda\partial_{1}\varphi({y}(u),{y}(v))z_{\lambda}(u)z'(v)(J(v))^{\lambda-1}+\\
&\qquad+\lambda\partial_v\varphi (u,v)z_{\lambda}(v)z'(u)(J(u))^{\lambda}+
\varphi (u,v)z'(u)z(v)\bigl[(J(v))^{\lambda}
\bigr]'(J(u))^{\lambda}-\\
&\qquad-\lambda\partial_{2}\varphi({y}(u),{y}(v))z_{\lambda}(v)z'(u)(J(u))^{\lambda-1}+
\lambda^2\partial^2_{u,v}\varphi (u,v)z_{\lambda}(u)z_{\lambda}(v)\\
&\qquad+\lambda\partial_u\varphi
(u,v)z_{\lambda}(u)z(v)\bigl[(J(v))^{\lambda}\bigr]'+
\lambda\partial_v\varphi
(u,v)z_{\lambda}(v)z(u)\bigl[(J(u))^{\lambda}\bigr]'+\\
&\qquad+\varphi
(u,v)z(u)z(v)\bigl[(J(u))^{\lambda}\bigr]'\bigl[(J(v))^{\lambda}\bigr]'-
\lambda^2\partial^2_{1,2}\varphi({y}(u),{y}(v))z_{\lambda}(u)z_{\lambda}(v).\endalign$$
We used above the formulae
$$\allowdisplaybreaks\align\partial_u\overline{\Omega}(u,v;{y})=&\partial_u\varphi (u,v)J(u)J(v)+\varphi (u,v)J'(u)J(v)-
\partial_{1}\varphi({y}(u),{y}(v))J(u),\quad \text{and}\\
\partial_v\overline{\Omega}(u,v;{y})=&\partial_v\varphi (u,v)J(u)J(v)+\varphi
(u,v)J(u)J'(v)-
\partial_{2}\varphi({y}(u),{y}(v))J(v).\endalign$$
 For the l.h.s\. of equation (7.4) we have
$$\allowdisplaybreaks\align l.h.s.=&\varphi
(u,v)z'_{\lambda}(u)z'_{\lambda}(v)+
\lambda\partial_u\varphi (u,v)z_{\lambda}(u)z'_{\lambda}(v)+\\
&\qquad+\lambda\partial_v\varphi (u,v)z'_{\lambda}(u)z_{\lambda}(v)+
\lambda^2\partial^2_{u,v}\varphi (u,v)z_{\lambda}(u)z_{\lambda}(v)\\
=&\varphi
(u,v)\biggl\{z'(u)(J(u))^{\lambda}+z(u)\left[(J(u))^{\lambda}\right]'\biggr\}
\biggl\{z'(v)(J(v))^{\lambda}+z(v)\left[(J(v))^{\lambda}\right]'\biggr\}+\\
&\qquad+\text{remaining terms}\\\allowdisplaybreak
=&\varphi (u,v)z'(u)z'(v)(J(u))^{\lambda}(J(v))^{\lambda}+
\varphi (u,v)z'(u)z(v)(J(u))^{\lambda}\left[(J(v))^{\lambda}\right]'+\\
&\qquad+\varphi (u,v)z(u)z'(v)(J(v))^{\lambda}\left[(J(u))^{\lambda}\right]'+
\varphi
(u,v)z(u)z(v)\left[(J(u))^{\lambda}\right]'\left[(J(v))^{\lambda}\right]'+\\
&\qquad+\lambda\partial_u\varphi (u,v)z_{\lambda}(u)z'(v)(J(v))^{\lambda}+
\lambda\partial_u\varphi
(u,v)z_{\lambda}(u)z(v)\left[(J(v))^{\lambda}\right]'+\\
&\qquad+\lambda\partial_v\varphi (u,v)z_{\lambda}(v)z'(u)(J(u))^{\lambda}+
\lambda\partial_v\varphi
(u,v)z_{\lambda}(v)z(u)\left[(J(u))^{\lambda}\right]'+\\
&\qquad+\lambda^2\partial^2_{u,v}\varphi
(u,v)z_{\lambda}(u)z_{\lambda}(v).\endalign$$
On the other hand we also have
$$\eqalign{\Omega(y(u),y(v);x)(J(u))^{\lambda}(J(v))^{\lambda}&=
\varphi({y}(u),{y}(v))z'(u)z'(v)(J(u))^{\lambda-1}(J(v))^{\lambda-1}+\cr
&\qquad+\lambda\partial_{1}\varphi({y}(u),{y}(v))z_{\lambda}(u)z'(v)(J(v))^{\lambda-1}+\cr
&\qquad+\lambda\partial_{2}\varphi({y}(u),{y}(v))z_{\lambda}(v)z'(u)(J(u))^{\lambda-1}+\cr
&\qquad+\lambda^2\partial^2_{1,2}\varphi({y}(u),{y}(v))z_{\lambda}(u)z_{\lambda}(v).}$$
\par After comparing the terms on the l.h.s\. and the r.h.s\. of equation (7.4)
we obtain an identity. This concludes the
proof of the Theorem.\qed\enddemo
\remark{Remark} Notice that for $\lambda\ne 0$ we can not have an inhomogeneous
term of the form
$\overline{\varphi}({x}(u),{x}(v))$ in (7.8). Had it been the case it would
impose on $\overline{\varphi}(u,v)$
the condition of being a homogeneous function of degree 1 in both arguments. In
order to satisfy the equation (7.4)
$\overline{\varphi}(u,v)$ must have the property
$\overline{\varphi}(z_{\lambda}(u),z_{\lambda}(v))=
\overline{\varphi}(z(u)(J(u))^{\lambda},z(v)(J(v))^{\lambda})=
\overline{\varphi}(z(u),z(v))(J(u))^{\lambda}(J(v))^{\lambda})$. But since
$\overline{\varphi}(u,v)$ must
be also antisymmetric it follows that the only function with these properties
is $\overline{\varphi}=0$.
On the contrary, for $\lambda=0$ we have a solution of (7.4) of the form
$$\Omega(u,v;x)={\varphi}(u,v)x'(u)x'(v)-\overline{\varphi}(x(u),x(v)),$$
where the function $\overline{\varphi}$ has the property
$\overline{\varphi}(u,v)=-\overline{\varphi}(v,u)$. The Jacobi
identity for $\Omega$ then implies that $\overline{\varphi}$ must satisfy
(7.6).
\endremark
\par Next we come to the following remarkable fact.
\proclaim{Theorem 7.4} If $\varphi (u,v)$ is defined by {\rm (7.6)}
then the solution {\rm (7.8)} satisfies the Jacobi identity, thus
defining a class of Poisson structures on $V_{\lambda}$ for which the action of
$G_{\infty}$ on $V_{\lambda}$ is a
Poisson action.\endproclaim
\demo{Proof} We shall use again the bracket $\bigl\{\Cal X(u),\Cal
X(v)\bigr\}:=\sum_{i,j=0}^{\infty}\{\Cal X_i,\Cal X_j\}u^iv^j$
introduced in Sec\. 4.
Then we have
$$\eqalign{\bigl\{\Cal X(u),\Cal X(v)\bigr\}=&\varphi (u,v)\Cal X'(u)\Cal
X'(v)+\lambda\partial_u\varphi (u,v)\Cal X(u)\Cal X'(v)+\cr
&+\lambda\partial_v\varphi (u,v)\Cal X'(u)\Cal
X(v)+\lambda^2\partial^2_{u,v}\varphi (u,v)\Cal X(u)\Cal X(v),} \tag 7.9$$
as well as
$$\multline\partial_u\bigl\{\Cal X(u),\Cal X(v)\bigr\}=\partial_v\varphi
(u,v)\Cal X'(u)\Cal X'(v)+
\varphi (u,v)\Cal X'(u)\Cal X''(v)+\\
+\lambda^2\partial^2_{u,v}\varphi (u,v)\Cal X(u)\Cal
X'(v)+\lambda\partial_u\varphi (u,v)\Cal X(u)\Cal
X''(v)+\lambda\partial^2_{v}\varphi (u,v)\Cal X'(u)\Cal
X(v)+\lambda\partial_v\varphi (u,v)\Cal X'(u)\Cal X'(v)+\cr
+\lambda^2\partial^3_{u,v^2}\varphi (u,v)\Cal X(u)\Cal
X(v)+\lambda^2\partial^2_{u,v}\varphi (u,v)\Cal X(u)\Cal X'(v).\endmultline$$
\par The Jacobi identity
$\bigl\{\Cal X(w),\bigl\{\Cal X(u),\Cal X(v)\bigr\}\bigr\}+c.p.=0$
is equivalent to the following equation:
$\Cal A+\Cal B+\Cal C+\Cal D=0.$
Here,
$$\align\Cal A=&\varphi (u,v)\bigl\{\Cal X(w),\Cal X'(u)\Cal
X'(v)\bigr\}+c.p.\\
=&\varphi (u,v)\biggl[\partial_u\bigl\{\Cal X(w),\Cal X(u)\bigr\}\Cal
X'(v)+\partial_v\bigl\{\Cal X(w),\Cal X(v)\bigr\}\Cal X'(u)\biggr]+
c.p.\\\allowdisplaybreak
\Cal B=&\lambda\partial_{u}\varphi (u,v)\bigl\{\Cal X(w),\Cal X(u)\Cal
X'(v)\bigr\}+c.p.\\
=&\lambda\partial_{u}\varphi (u,v)\biggl[\bigl\{\Cal X(w),\Cal X(u)\bigr\}\Cal
X(v)+
\partial_v\bigl\{\Cal X(w),\Cal X(v)\bigr\}\Cal X(u)\biggr]+
c.p.\\\allowdisplaybreak
\Cal C=&\lambda\partial_{v}\varphi (u,v)\bigl\{\Cal X(w),\Cal X'(u)\Cal
X(v)\bigr\}+c.p.\\
=&\lambda\partial_{v}\varphi (u,v)\biggl[\partial_u\bigl\{\Cal X(w),\Cal
X(u)\bigr\}\Cal X(v)+
\bigl\{\Cal X(w),\Cal X(v)\bigr\}\Cal X'(u)\biggr]+
c.p.\\\allowdisplaybreak
\Cal D=&\lambda^2\partial^2_{u,v}\varphi (u,v)\bigl\{\Cal X(w),\Cal X(u)\Cal
X(v)\bigr\}+c.p.\\
=&\lambda^2\partial^2_{u,v}\varphi (u,v)\biggl[\bigl\{\Cal X(w),\Cal
X(u)\bigr\}\Cal X(v)+
\bigl\{\Cal X(w),\Cal X(v)\bigr\}\Cal X(u)\biggr]+
c.p.\endalign$$
For the expressions in the square brackets for each term we therefore obtain
$$\allowdisplaybreaks\align\Cal A'=&\partial_{u}\varphi (w,u)\Cal X'(w)\Cal
X'(u)\Cal X'(v)+\varphi (w,u)\Cal X'(w)\Cal X''(u)\Cal X'(v)+\\
&\qquad+\lambda\partial^2_{w,u}\varphi (w,u)\Cal X(w)\Cal X'(u)\Cal
X'(v)+\lambda\partial_{w}\varphi (w,u)\Cal X(w)\Cal X''(u)\Cal X'(v)+\\
&\qquad+\lambda\partial^2_{u}\varphi (w,u)\Cal X'(w)\Cal X(u)\Cal
X'(v)+\lambda\partial_{u}\varphi (w,u)\Cal X'(w)\Cal X'(u)\Cal X'(v)+\\
&\qquad+\lambda^2\partial^3_{w,u^2}\varphi (w,u)\Cal X(w)\Cal X(u)\Cal X'(v)+
\lambda^2\partial^2_{w,u}\varphi (w,u)\Cal X(w)\Cal X'(u)\Cal X'(v)+\\
&\qquad+\partial_{v}\varphi (w,v)\Cal X'(w)\Cal X'(u)\Cal X'(v)+\varphi
(w,v)\Cal X'(w)\Cal X'(u)\Cal X''(v)+\\
&\qquad+\lambda\partial^2_{w,v}\varphi (w,v)\Cal X(w)\Cal X'(u)\Cal
X'(v)+\lambda\partial_{w}\varphi (w,v)\Cal X(w)\Cal X'(u)\Cal X''(v)+\\
&\qquad+\lambda\partial^2_{v}\varphi (w,v)\Cal X'(w)\Cal X'(u)\Cal
X(v)+\lambda\partial_{v}\varphi (w,v)\Cal X'(w)\Cal X'(u)\Cal X'(v)+\\
&\qquad+\lambda^2\partial^3_{w,v^2}\varphi (w,v)\Cal X(w)\Cal X'(u)\Cal X(v)+
\lambda^2\partial^2_{w,v}\varphi (w,v)\Cal X(w)\Cal X'(u)\Cal X'(v),\endalign$$
$$\allowdisplaybreaks\align\Cal B'=&\varphi (w,u)\Cal X'(w)\Cal X'(u)\Cal
X'(v)+\lambda\partial_{w}\varphi (w,u)\Cal X(w)\Cal X'(u)\Cal X'(v)+\\
&\qquad+\lambda\partial_{u}\varphi (w,u)\Cal X'(w)\Cal X(u)\Cal
X'(v)+\lambda^2\partial^2_{w,u}\varphi (w,u)\Cal X(w)\Cal X(u)\Cal X'(v)+\\
&\qquad+\partial_{v}\varphi (w,v)\Cal X'(w)\Cal X(u)\Cal X'(v)+\varphi
(w,v)\Cal X'(w)\Cal X(u)\Cal X''(v)+\\
&\qquad+\lambda\partial^2_{w,v}\varphi (w,v)\Cal X(w)\Cal X(u)\Cal
X'(v)+\lambda\partial_{w}\varphi (w,v)\Cal X(w)\Cal X(u)\Cal X''(v)+\\
&\qquad+\lambda\partial^2_{v}\varphi (w,v)\Cal X'(w)\Cal X(u)\Cal
X(v)+\lambda\partial_{v}\varphi (w,v)\Cal X'(w)\Cal X(u)\Cal X'(v)+\\
&\qquad+\lambda^2\partial^3_{w,v^2}\varphi (w,v)\Cal X(w)\Cal X(u)\Cal X(v)+
\lambda^2\partial^2_{w,v}\varphi (w,v)\Cal X(w)\Cal X(u)\Cal X'(v),\endalign$$
$$\allowdisplaybreaks\align\Cal C'=&\partial_{u}\varphi (w,u)\Cal X'(w)\Cal
X'(u)\Cal X(v)+\varphi (w,u)\Cal X'(w)\Cal X''(u)\Cal X(v)+\\
&\qquad+\lambda\partial^2_{w,u}\varphi (w,u)\Cal X(w)\Cal X'(u)\Cal
X(v)+\lambda\partial_{w}\varphi (w,u)\Cal X(w)\Cal X''(u)\Cal X(v)+\\
&\qquad+\lambda\partial^2_{u}\varphi (w,u)\Cal X'(w)\Cal X(u)\Cal
X(v)+\lambda\partial_{u}\varphi (w,u)\Cal X'(w)\Cal X'(u)\Cal X(v)+\\
&\qquad+\lambda^2\partial^3_{w,u^2}\varphi (w,u)\Cal X(w)\Cal X(u)\Cal X(v)+
\lambda^2\partial^2_{w,u}\varphi (w,u)\Cal X(w)\Cal X'(u)\Cal X(v)+\\
&\qquad+\varphi (w,v)\Cal X'(w)\Cal X'(u)\Cal X'(v)+\lambda\partial_{w}\varphi
(w,v)\Cal X(w)\Cal X'(u)\Cal X'(v)+\\
&\qquad+\lambda\partial^2_{w,v}\varphi (w,v)\Cal X(w)\Cal X'(u)\Cal
X(v)+\lambda\partial_{v}\varphi (w,v)\Cal X'(w)\Cal X'(u)\Cal X(v),\endalign$$
$$\allowdisplaybreaks\align\Cal D'=&\varphi (w,u)\Cal X'(w)\Cal X'(u)\Cal
X(v)+\lambda\partial_{w}\varphi (w,u)\Cal X(w)\Cal X'(u)\Cal X(v)+\\
&\qquad+\lambda\partial_{u}\varphi (w,u)\Cal X'(w)\Cal X(u)\Cal
X(v)+\lambda^2\partial^2_{w,u}\varphi (w,u)\Cal X(w)\Cal X(u)\Cal X(v)+\\
&\qquad+\varphi (w,v)\Cal X'(w)\Cal X(u)\Cal X'(v)+\lambda\partial_{w}\varphi
(w,v)\Cal X(w)\Cal X(u)\Cal X'(v)+\\
&\qquad+\lambda\partial_{v}\varphi (w,v)\Cal X'(w)\Cal X(u)\Cal
X(v)+\lambda^2\partial^2_{w,v}\varphi (w,u)\Cal X(w)\Cal X(u)\Cal
X(v).\endalign$$
\par As a result we will study the four terms in the Jacobi identity written as
$$\xalignat3 \Cal A&=\varphi (u,v)\Cal A'+c.p.&\qquad\Cal
B&=\lambda\partial_{u}\varphi (u,v)\Cal B'+c.p.\\
\Cal C&=\lambda\partial_{v}\varphi (u,v)\Cal C'+c.p.&\qquad\Cal
D&=\lambda^2\partial^2_{u,v}\varphi (u,v)\Cal D'+c.p.\endxalignat$$
\par We split the analysis of the Jacobi identity into seven steps (A)-(G).
\newline
 (A) Terms of the form $\Cal X'\Cal X'\Cal X'$. From them we obtain (after
cyclicly permuting the arguments of some of them):
$$\biggl\{\varphi (u,v)\bigl[\partial_{u}\varphi (w,u)+\partial_{v}\varphi
(w,v)\bigr]+c.p.\biggr\}\Cal X'(w)\Cal X'(u)\Cal X'(v).$$
But from (7.6) it follows that the above term is zero.\newline
 (B) Terms of the form $\Cal X\Cal X''\Cal X'$ cancel each other out after
cyclic permutation.\newline
 (C) Terms of the form $\Cal X\Cal X''\Cal X$ cancel each other out after
cyclic permutation.\newline
 (D) Terms of the form $\Cal X\Cal X\Cal X$ give
$$\multline\lambda^3\biggl[\partial_{u}\varphi (u,v)\partial^3_{w,v^2}\varphi
(w,v)+
\partial_{v}\varphi (u,v)\partial^3_{w,u^2}\varphi (w,u)+c.p.\biggr]\Cal
X(w)\Cal X(u)\Cal X(v)+\\
+\lambda^4\biggl[\partial^2_{u,v}\varphi (u,v)\biggl(\partial^2_{w,u}\varphi
(w,u)+\partial^2_{w,v}\varphi (w,v)\biggr)+c.p.\biggr]
\Cal X(w)\Cal X(u)\Cal X(v).\endmultline$$
Since $\varphi (u,v)$ is a solution of (7.6) the results obtained in Ch. VI
showed that  $\varphi (u,v)=f(u)g(v)-f(v)g(u)$, where
$f$ and $g$ satisfy (7.7). Therefore the $\lambda^3$-term becomes
$$\multline\lambda^3\biggl\{\biggl[f'(u)g(v)-f(v)g'(u)\biggr]\biggl[f'(w)g''(v)-f''(v)g'(w)\biggr]+\\
+\biggl[f(u)g'(v)-f'(v)g(u)\biggr]\biggl[f'(w)g''(u)-f''(u)g'(w)\biggr]+c.p.\biggr\}\Cal X(w)\Cal X(u)\Cal X(v),\endmultline$$
and the $\lambda^4$-term assumes the form
$$\multline \lambda^4\biggl\{\biggl[f'(u)g'(v)-f'(v)g'(u)\biggr]\times\\
\times\biggl[f'(w)g'(u)-f'(u)g'(w)+f'(w)g'(v)-f'(v)g'(w)\biggr]
+c.p.\biggr\}\Cal X(w)\Cal X(u)\Cal X(v).\endmultline$$
Using the identities
$f'(u)g(u)-f(u)g'(u)=\alpha f'(u)+\beta g'(u)\hskip 0.2in \Bigl(\implies
f''(u)g(u)-f(u)g''(u)=\alpha f''(u)+\beta g''(u)\Bigr)$
one easily shows that both terms proportional to $\lambda^3$ and $\lambda^4$
respectively are identically zero.\newline
 (E) Terms of the form $\Cal X\Cal X\Cal X'$ give
$$\multline\lambda^2\biggl[\varphi (u,v)\partial^3_{w,u^2}\varphi (w,u)+\varphi
(v,w)\partial^3_{u,w^2}\varphi (u,w)+\\
+\partial_{u}\varphi (u,v)\partial^2_{w,v}\varphi (w,v)+\partial_{w}\varphi
(w,u)\partial^2_{u}\varphi (v,u)\\
+\partial_{u}\varphi (w,u)\partial^2_{w}\varphi (v,w)+\partial_{w}\varphi
(v,w)\partial^2_{u,v}\varphi (u,v)\biggr]
\Cal X(w)\Cal X(u)\Cal X'(v)+c.p.\endmultline$$
The expression in the square brackets becomes identically zero after using
(7.7) in a similar way as in (D). Also we have a term
proportional to $\lambda^3$ which reads
$$\multline\lambda^3\biggl[\partial_{u}\varphi (u,v)\partial^2_{w,u}\varphi
(w,u)+\partial_{u}\varphi (u,v)\partial^2_{w,v}\varphi (w,v)+\\
+\partial_{w}\varphi (v,w)\partial^2_{u,v}\varphi (u,v)+\partial_{w}\varphi
(v,w)\partial^2_{u,w}\varphi (u,w)+\\
+\partial_{w}\varphi (v,w)\partial^2_{w,u}\varphi (w,u)+\partial_{u}\varphi
(u,v)\partial^2_{v,w}\varphi (v,w)+\\
+\partial_{w}\varphi (w,v)\partial^2_{u,v}\varphi (u,v)+\partial_{u}\varphi
(v,u)\partial^2_{w,u}\varphi (w,u)\biggr]
\Cal X(w)\Cal X(u)\Cal X'(v)+c.p.,\endmultline$$
which is identically zero.\newline
 (F) Terms of the form $\Cal X\Cal X'\Cal X'$ give
$$\multline\lambda\biggl[\varphi (u,v)\partial^2_{w,u}\varphi (w,u)+\varphi
(w,u)\partial^2_{w}\varphi (v,w)+\\
+\varphi (u,v)\partial^2_{w,v}\varphi (w,v)+\varphi (v,w)\partial^2_{w}\varphi
(u,w)+\\
+\partial_{w}\varphi (w,u)\partial_{u}\varphi (v,u)+\partial_{w}\varphi
(v,w)\partial_{v}\varphi (u,v)\biggr]
\Cal X(w)\Cal X'(u)\Cal X'(v)+c.p.\endmultline$$
The expression in the square brackets can be shown to be identically zero after
using (7.7) in a similar way as in (D) and (E).
\par The other two terms of the same form are
$$\multline\lambda^2\biggl[\varphi (u,v)\partial^2_{w,u}\varphi (w,u)+\varphi
(u,v)\partial^2_{w,v}\varphi (w,v)+\\
+\varphi (v,u)\partial^2_{w,u}\varphi (w,u)+\varphi
(u,v)\partial^2_{w,v}\varphi (v,w)\biggr]
\Cal X(w)\Cal X'(u)\Cal X'(v)+c.p.,\endmultline$$
which is identically zero, and
$$\multline\lambda^2\biggl[\partial_{w}\varphi (w,u)\partial_{u}\varphi
(u,v)+\partial_{w}\varphi (w,u)\partial_{w}\varphi (v,w)+\\
+\partial_{w}\varphi (w,u)\partial_{u}\varphi (v,u)+\partial_{w}\varphi
(v,w)\partial_{v}\varphi (u,v)+\\
+\partial_{v}\varphi (u,v)\partial_{w}\varphi (w,v)+\partial_{w}\varphi
(v,w)\partial_{w}\varphi (u,w)\biggr]
\Cal X(w)\Cal X'(u)\Cal X'(v)+c.p.,\endmultline$$
which is again identically zero.\newline
 (G) Terms of the form $\Cal X'\Cal X''\Cal X'$ cancel each other.
\par Thus all terms have been covered and the proof of Theorem  7.4 is
completed.\qed\enddemo
\par The consequence of Theorem 7.3 and Theorem 7.4 is that for each
Poisson-Lie structure on $G_{\infty}$ defined by a function
$\varphi$ satisfying the equation (7.6) there exists a Poisson structure on
$V_{\lambda}$ for which the action of $G_{\infty}$
is Poisson. Thus we obtain a series of representations
$\Cal V_{\varphi,\lambda}$ of the Poisson-Lie group $G_{\infty}$ in the
Poisson
spaces $V_{\lambda}$.
\head 8. Quantization\endhead
\par This section is devoted to the quantization of some of the Poisson-Lie
structures on the group $G_{\infty}$ restricted to the finite dimensional
factor-groups $G_n$.
We shall construct explicitly families of finitely generated quantum
(semi)groups.
Their quasi-classical limits are the finite-dimensional Poisson-Lie groups
endowed with Poisson-Lie structures, which are
restrictions of the Poisson-Lie structures on the group $G_{\infty}$, belonging
to the countable family obtained by Theorem 4.4.
This means that we shall consider factor-groups $G_n=G_{\infty}\text{mod}
\{u^{n+1}\}$, for $n\ge 5$ (cf\. Sec\. 2). In this aproach to quantization we
shall start from the quasi-classical limits, i.e\., the corresponding
Poisson-Lie groups, and
reconstruct from this data their quantum counterparts. This
quantization procedure is a deformation
quantization of the Poisson algebra $\Cal A$ of $C^{\infty}$ functions on the
corresponding finite-dimensional Poisson-Lie groups to a
non-commutative non-cocommutative bialgebra $\Cal A_h$ such that $\Cal
A_h/h\Cal A_h\cong \Cal A$. The second postulate of quantization requires
that the deformation be `flat', i.e., the dimension of
$\Cal A_h$ as a $k[[h]]$-module, for a field $k$ of characteristic zero, be the
same as the dimension of $\Cal
A$ as a $k$-module. For the general philosophy underlying the method
we refer the reader to [D2].
\par (i) Let $X=\{x_i\}_{i\in\zN}$ be a set of indeterminates. Let us introduce
a grading on the algebra $k[X]$, where $k$ denotes the ground field (assumed to
be of characteristic zero), by assigning a degree (denoted $|\ |$) to each of
the generators
$x_i$ of $k[X]$ by the following definition:
$$|x_i|=i-1, \qquad \text{for every  } i\in\zN. \tag 8.1$$
We extend it on monomials by $|AB|=|A|+|B|$, for every two monomials $ A,B.$
\par As we have mentioned above, in this section we shall be concerned with the
quantization problem for the countable family of Poisson-Lie
structures on the group $G_{\infty}$, found in Sec\. 4, the formulae for which
we now recall:
$$\eqalign{\{x_i,x_j\}=&(i-d)jx_jx_{i-d}-i(j-d)x_ix_{j-d}\cr &\qquad
+x_i\sum_{\bigl(\sum_{k=1}^{d+1}s_k\bigr)=j}^{}x_{s_1}\ldots
x_{s_{d+1}}-x_j\sum_{\bigl(\sum_{k=1}^{d+1}s_k\bigr)=i}^{}x_{s_1}\ldots
x_{s_{d+1}},\qquad\text{for all}\quad {d\in\zN}.} \tag 8.2$$
It is clear from the right hand side of (8.2) that for each $d\in\zN$ the
degree of the bracket $\{x_i,x_j\}$ is given by
$$|\{x_i,x_j\}|=|x_i|+|x_j|-d=i+j-d-2. \tag 8.3$$
\par (ii) Let $X=\{x_i\}_{i\in\zN}$ be a set of indeterminates and let
$\langle X\rangle$ be a free associative semigroup with identity on
$X$. Let $k[[h]]\langle X\rangle$ be the semigroup algebra of $\langle
X\rangle$ over the ring of formal power series $k[[h]]$ in the
parameter $h$. Here $k$ is assumed to be a field of characteristic
zero. Consider the
set of relations
 $$\Cal R_h=\{x_ix_j-x_jx_i=f_{ij|h}(x)\mid i,j\in\zN\}\tag 8.4$$
where $f_{ij|h}(x)$ are polynomials
in $x_i$ with coefficients in $k[[h]]$, such that $f_{ij|0}=0$. Let $\Cal I_h$
be
the ideal generated by $\Cal R_h$. Define $\Cal A_h:=k[[h]]\langle
X\rangle/\Cal I_h$, and define a grading on $k\langle X\rangle$ as
was done in (i).
\par The following is the first postulate of quantization. As
explained earlier we require that $\Cal
A_h/h\Cal A_h\cong k\langle X\rangle$, i.e.,
$$\bigl[x_i,x_j\bigr]=h\{x_i,x_j\}+{ O}(h^2).\tag 8.5$$
Here $[x_i,x_j]=x_ix_j-x_jx_i$, and the product is the product in $\Cal A_h$.
In other words we would like to recover the Poisson-Lie bracket on $G_{\infty}$
(or the factor groups $G_n=G_{\infty}\text{mod }u^{n+1}$)
in the quasi-classical limit $h\to 0$. This also means that for $h\to 0$ we
should have
$$f_{ij|h}(x)=h\{x_i,x_j\}+{ O}(h^2).\tag 8.6$$
After computing the degree of the right hand side of the above equality
$|h\{x_i,x_j\}|=|h|+|\{x_i,x_j\}|=|h|+i+j-d-2,$
we deduce that for each $d\in\zN$ the parameter $h$ must have degree $|h|=d$,
since
$|[x_i,x_j]|=i+j-2.$
\par (iii) Consider the semigroup algebra $k[[h]]\langle X\rangle$ of
$\langle X\rangle$ over the ring $k[[h]]$. For each $d\in\zN$ consider the set
of relations
$$
\Cal R^d_h=\{x_ix_j=x_jx_i+f^d_{ij|h}(x)\mid i<j  \hskip 0.1in \text{for}
\hskip 0.1in i,j\in\zN\}\tag 8.7
$$
where $f^d_{ij|h}(x)\in  k[[h]]\langle X\rangle$ are linear combinations of
monomials $h^nx^{n_1}_{i_{1}}\ldots x^{n_k}_{i_{k}}$ such that
$i_1>\ldots >i_k$ and $nd+\sum_{s=1}^kn_s(i_s-1)=i+j-d-2$. We shall
say that these monomials are in canonical form. Thus $f^d_{ij|h}(x)$
are linear combinations of monomials in canonical form. Recall the following
definition.
\definition{Definition 8.1} The semigroup algebra $k[[h]]\langle
X\rangle$ has a Poincar\'e-Birkhoff-Witt (PBW) property if every
monomial $A\in \langle X\rangle$ can be reduced to a unique expression as a
linear combination of monomials in canonical form using the
set of relations (8.7) independently of the choice of a reduction procedure.
\enddefinition
We shall use a version of the result
known as the Diamond Lemma [Be], applicable here, which allows us to prove the
PBW
property, by only proving it for the monomials with so called
`overlap' ambiguities [Be].
\definition{Definition 8.2} For each $d\in\zN$ a quantum semigroup
$G^d_{\infty|h}$ is defined as follows. As a `quantum space' $G^d_{\infty|h}$
is defined by its factor semigroup algebra $\Cal A^d_h:=k[[h]]\langle
X\rangle/\Cal I^d_h$, where $\Cal I^d_h\subset k[[h]]\langle X\rangle$ is the
ideal generated
by the set of relations (8.7), and we require that in the quasi-classical limit
$\bigl[x_i,x_j\bigr]=h\{x_i,x_j\}+\text{ O}(h^2)$,
one obtains the Poisson algebra of functions on the Poisson-Lie group
$G^d_{\infty}$ defined by (8.2).
The  multiplication of formal power series in one variable under the operation
of substitution
$\left(xy\right)(u)=x\left(y(u)\right),$
where $x(u)=\sum_{i=1}^{\infty}x_iu^{i}$, and
$y(u)=\sum_{i=1}^{\infty}y_iu^{i}$, gives rise to a comultiplication map
$\Delta\:\Cal A^d_h\to\Cal A^d_h\otimes\Cal A^d_h,$
which is defined on the generators by
$$\Delta(x_n)=\sum_{i=1}^n
x_i\otimes\sum_{\sum_{\alpha=1}^{i}{j_{\alpha}}=n}{x_{j_1}\ldots
x_{j_i}},\qquad n\in \zN,$$
and is required to be an algebra homomorphism.
Also one defines a counit map $c\:\Cal A^d_h\to k[[h]]$ by
$$c(x_i)=\delta_i^1,\qquad i\in \zN.$$
All tensor products are over $k[[h]]$.
This endows $\Cal A^d_h$ with a structure of a bialgebra and the quantum
semigroup $G^d_{\infty|h}$ is defined to be the
bialgebra $\Cal A^d_h$, if in addition $k[[h]]\langle X\rangle$ has the
PBW basis described above.
\enddefinition
\par Does such an object exist? We do not know yet. However, if we consider
the Poisson-Lie factor groups
$G_n^d=G^d_{\infty}\text{mod} \{u^{n+1}\}$, for $n\ge 5$, then there
exist quantum objects that satisfy the above definition for small
values of $d$,
$d\le 5$, and whose quasi-classical limits are the Poisson-Lie groups
$G^d_n$. The definition of the finitely generated quantum semigroups
$G^d_{n|h}$ is the same as the definition above with $\Cal I^d_h$ being an
ideal generated by a finite set of relations $\Cal R^d_h$.
Their defining semigroup algebras turn out to have other interesting
properties.
We shall describe in more detail the construction of the quantum semigroup
$G^2_{5|h}$, while omitting parts of the construction that
consist of lengthy and tedious calculations, and we shall state the results for
the quantum semigroups $G^1_{4|h}$ and $G^3_{5|h}$ without
entering into the details of the calculations.
 The construction of the last two mimics exactly the construction of
$G^2_{5|h}$.
\par Let $G^2_5=G^2_{\infty}\text{mod}\{u^{n+1}\}$, for $n\ge 5$, be
the finite dimensional (dim\ =\ 5, $d=2$) Poisson-Lie group with a Poisson-Lie
 structure defined by
$$\eqalign{\{x_i,x_j\}=&(i-2)jx_jx_{i-2}-i(j-2)x_ix_{j-2}\cr
&\qquad+x_i\sum_{s_1+s_2+s_3=j}x_{s_1}x_{s_2}x_{s_3}
-x_j\sum_{s_1+s_2+s_3=i}x_{s_1}x_{s_2}x_{s_3}.}\tag 8.10$$
The above formulae are obtained from (8.2) with $d=2$, and we adopt the
convention that $x_i=0$ whenever $i<1$. In more detail
we have
$$\allowdisplaybreaks\align\{x_1,x_2 \}&=0,\\
\{x_1,x_3 \}&=-x_1^2+x_1^4,\\
\{x_2,x_3\}&=x_2(x_1^3-2x_1),\\
\{x_1,x_4\}&=x_2(3x_1^3-2x_1),\\
\{x_2,x_4\}&=x_2^2(3x_1^3-4),\\
\{x_3,x_4\}&=x_4(4x_1-x_1^3)+x_3x_2(3x_1^2-6),\tag 8.11\\
\{x_1,x_5\}&=x_3(3x_1^3-3x_1)+3x_2^2x_1^2,\\
\{x_2,x_5\}&=3x_2^3x_1+x_3x_2(3x_1^2-6),\\
\{x_3,x_5\}&=x_5(5x_1-x_1^3)+x_3^2(3x_1^2-9)+3x_3x_2^2x_1,\\
\{x_4,x_5\}&=x_5x_2(10-3x_1^2)+x_4x_3(3x_1^2-12)+3x_4x_2^2x_1.\endalign$$
\proclaim{Theorem 8.3} Let $X=\{x_i\}_{1\le i\le 5}$ be a set of
indeterminatesand let $\langle X\rangle$ be the associative semigroup with
identity generated by $X$.
Consider an ideal $\Cal I_{h}^2$ generated by the set of relations  $\Cal
R^2_h$ in $k[[h]]\langle X\rangle$
$$\align x_1x_2 =x_2x_1\\
x_1x_3 =x_3x_1&+h(-x_1^2+x_1^4)\\
x_2x_3=x_3x_2&+hx_2(x_1^3-2x_1)\\\allowdisplaybreak
x_1x_4=x_4x_1&+hx_2(3x_1^3-2x_1)\\\allowdisplaybreak
x_2x_4=x_4x_2&+hx_2^2(3x_1^3-4)\\\allowdisplaybreak
x_3x_4=x_4x_3&+h[x_4(4x_1-x_1^3)+x_3x_2(3x_1^2-6)]\\ &+
2h^2x_2x_1\tag 8.12\\\allowdisplaybreak
x_1x_5=x_5x_1&+h[x_3(3x_1^3-3x_1)+3x_2^2x_1^2]\\ &+
h^2(-6x_1^4+\frac {9}{2}x_1^6+\frac {3}{2}x_1^2)\\\allowdisplaybreak
x_2x_5=x_5x_2&+h[3x_2^3x_1+x_3x_2(3x_1^2-6)]\\ &+
h^2x_2(6x_1-9x_1^3+\frac {9}{2}x_1^5)\\\allowdisplaybreak
x_3x_5=x_5x_3&+h[x_5(5x_1-x_1^3)+x_3^2(3x_1^2-9)+3x_3x_2^2x_1]\\ &+
h^2x_3(-\frac {15}{2}x_1+6x_1^3+\frac {3}{2}x_1^5)\\ &+
h^3C(x_1^8-x_1^2)\\\allowdisplaybreak
x_4x_5=x_5x_4&+h[x_5x_2(10-3x_1^2)+x_4x_3(3x_1^2-12)+3x_4x_2^2x_1]\\ &+
h^2\bigl[x_4(-24x_1+9x_1^3+\frac {3}{2}x_1^5)+6x_3x_2\bigr]\\ &+
h^3x_2\bigl[-(6+2C)x_1+3Cx_1^7\bigr],\endalign$$
where $C\in k$ is an arbitrary parameter.  Then the semigroup factor algebra
$k[[h]]\langle X\rangle/\Cal I_{h}^2$ defines a quantum
semigroup $G^2_{5|h}$ in
the sense of Definition {\rm 8.2} with a comultiplication defined by {\rm
(8.8)}. Namely,
$$\allowdisplaybreaks\align \Delta x_1&= x_1\otimes x_1\\
\Delta x_2&= x_1\otimes  x_2 + x_2\otimes  x_1^2\\
\Delta x_3&= x_1\otimes  x_3 + x_2\otimes  x_1 x_2 + x_2\otimes  x_2
x_1+x_3\otimes  x_1^3 \\
\Delta x_4&= x_1\otimes  x_4 + x_2\otimes  x_1 x_3 + x_2\otimes  x_2^2
+x_2\otimes  x_3 x_1 + x_3\otimes  x_1^2x_2+
      x_3\otimes  x_1 x_2 x_1 + x_3\otimes  x_2 x_1^2 + x_4\otimes  x_1^4\\
\Delta x_5&= x_1\otimes  x_5 + x_2\otimes  x_1 x_4 + x_2\otimes  x_2 x_3+
x_2\otimes  x_3 x_2 + x_2 x_4 x_1 + x_3 x_1^2x_3 \\ &\qquad+
   x_3\otimes  x_1 x_2^2 + x_3\otimes  x_1 x_3 x_1+ x_3\otimes  x_2 x_1 x_2 +
x_3\otimes  x_2^2 x_1+x_3\otimes  x_3 x_1^2+ x_4\otimes  x_1^3 x_2\\ &\qquad+
   x_4\otimes  x_1^2 x_2 x_1 + x_4 \otimes x_1 x_2 x_1^2 +x_4\otimes  x_2 x_1^3
+ x_5\otimes  x_1^4.\tag 8.13\endalign$$
Moreover the semigroup algebra $k[[h]]\langle X\rangle$ has the
Poincar\'e-Birkhoff-Witt property.\endproclaim
\demo{Proof} The proof is constructive. We look for a set of relations $\Cal
R^2_h$ in $k[[h]]\langle X\rangle$ in the following form
$$\align x_1x_2 =x_2x_1\\
x_1x_3=x_3x_1&+h(-x_1^2+x_1^4)\\
x_2x_3=x_3x_2&+hx_2(x_1^3-2x_1)\\\allowdisplaybreak
x_1x_4=x_4x_1&+hx_2(3x_1^3-2x_1)\\\allowdisplaybreak
x_2x_4=x_4x_2&+hx_2^2(3x_1^3-4)+h^2f_1(x_1)\\\allowdisplaybreak
x_3x_4=x_4x_3&+h[x_4(4x_1-x_1^3)+x_3x_2(3x_1^2-6)]\\
&+h^2x_2f_2(x_1)\\\allowdisplaybreak
x_1x_5=x_5x_1&+h[x_3(3x_1^3-3x_1)+3x_2^2x_1^2]\\ &+
h^2f_3(x_1)\\\allowdisplaybreak
x_2x_5=x_5x_2&+h[x_3x_2(3x_1^2-6)+3x_2^3x_1]\\ &+
h^2x_2f_4(x_1)\\\allowdisplaybreak
x_3x_5=x_5x_3&+h[x_5(5x_1-x_1^3)+x_3^2(3x_1^2-9)+3x_3x_2^2x_1]\\ &+
h^2\bigl[x_3f_5(x_1)+x_2^2f_6(x_1)\bigr]\\ &+
h^3f_7(x_1)\\\allowdisplaybreak
x_4x_5=x_5x_4&+h[x_5x_2(10-3x_1^2)+x_4x_3(3x_1^2-12)+3x_4x_2^2x_1]\\ &+
h^2\bigl[x_4f_8(x_1)+x_3x_2f_9(x_1)+x_2^3f_{10}(x_1)\bigr]\\ &+
h^3x_2f_{11}(x_1),\tag 8.14\endalign$$
where $\{f_i(x)\}_{1\le i\le 11}$ is a set of arbitrary polynomials. Since the
degree of $h$ is $|h|=2$ and the degree of $x_1$ is $|x_1|=0$
this is the most general form of the set of relations that one can have, such
that their quasi-classical limit gives the
Poisson-Lie structure on $G^2_5$.
\par We now require $k[[h]]\langle X\rangle/\Cal I_{h}^2$ to be a
bialgebra, i.e., that $\Delta$ be an algebra homomorphism. This leads
to some restrictions on the polynomials $f_i(x)$. Using the formulae
for comultiplication
(8.13) one can see that the first four relations are compatible with
the coalgebra structure. One obtains functional equations for $f_i(x)$ from
the remaining six relations after a reduction to a canonical form.
We shall analyse first the equations that arise from terms of order $h^2$.
\par (a) The compatibility with comultiplication induces that
$$\Delta (x_2x_4)=\Delta (x_4x_2)+h(-4\Delta x_2^2+3\Delta
(x_2^2x_1^2))+h^2f_1(\Delta x_1).$$
After reducing both sides of this relation to a canonical form,  using the
comultiplication formulae (8.13), we obtain the following linear functional
equation for $f_1(x)$:
$-f_1(x_1\otimes x_1)+x_1^2\otimes f_1(x_1)+x_1^6\otimes f_1(x_1)=0.$
 From now on since all equations for the unknowns $f_i$ will depend only on the
variable $x_1$ we shall use the notation
$x:= x_1$. Therefore we have
$$-f_1(x\otimes x)+x^2\otimes f_1(x)+x^6\otimes f_1(x)=0.\tag  8.15$$
The most general solution of the above equation is
$$f_1(x)=C_1(x^6-x^2),\tag 8.16$$
where $C_1\in k$ is an arbitrary constant. There are no terms of higher order
in $h$ that arise in the analysis of this relation. We
move on to the next.
\par (b) Again after reducing to a canonical form both sides of
$$\Delta(x_3x_4)=\Delta(x_4x_3)+h(\Delta(4x_4x_1)-\Delta(x_4x_1^3)-\Delta(6x_3x_2)+\Delta(3x_3x_2x_1^2))+
h^2\Delta(x_2f_2(x_1))$$
we obtain two equations. One of them arises from a term proportional
to $x_2\otimes 1$, i.e., we have a term of the form
$$(x_2\otimes 1)\bigl[-f_2(x_1\otimes x_1)+(2-2C_1)x_1\otimes
x_1+f_2(x_1)\otimes x_1^5+(-2+2C_1)x_1\otimes x_1^5\bigr],$$
which does not cancel.
It leads to the equation
$$-f_2(x\otimes x)+(2-2C_1)x\otimes x+f_2(x)\otimes x^5+(-2+2C_1)x\otimes
x^5=0.\tag 8.17$$
The term proportional to $1\otimes x_2$ leads to
$$-f_2(x\otimes x)+x\otimes f_2(x)-2C_1x\otimes x^5+2C_1x^5\otimes x^5=0.\tag
8.18$$
Solving (8.17) and (8.18) together we obtain
$$f_2(x)=(2-2C_1)x+2C_1x^5.\tag 8.19$$
There are no terms of higher order in $h$ arising from this relation.
\par At this stage of the calculation we check whether the PBW property is
satisfied in the subalgebra of $k[[h]]\langle X\rangle$ generated
by the set $\{x_1,x_2,x_3,x_4\}$ and subject to the first six relations of
(8.14). By direct calculation, using the Diamond Lemma, one shows that the
monomial $x_2x_3x_4$ can be reduced to a unique canonical form if and only if
$C_1=0$. The other possible monomials of three
variables have a unique canonical form. Therefore we obtain that
$$f_1(x)=0 \qquad \text{and} \qquad f_2(x)=2x.\tag 8.20$$
\par (c) From the next relation one has
$$\Delta (x_1x_5)=\Delta (x_5x_1)+h(-3\Delta (x_3x_1)+3\Delta
(x_2^2x_1^2)+3\Delta (x_3x_1^3))+h^2f_3(\Delta x_1),$$
after a reduction to a canonical one is lead to the equation
$$-f_3(x\otimes x)+x^2\otimes f_3(x)+6x^2\otimes x^4+x^6\otimes
f_3(x)-6x^4\otimes x^4-6x^2\otimes x^6+6x^4\otimes x^6=0.\tag 8.21$$
The most general solution of the above equation is given by
$$f_3(x)=6x^2-6x^4+C_2(x^6-x^2),\tag 8.22$$
where $C_2\in k$ is an arbitrary constant. There are no terms of order $h^3$ or
higher that arise after the reduction to a canonical form
of the above relation.
\par (d) The analysis of the next relation
$$\Delta (x_2x_5)=\Delta (x_5x_2)+h(-6\Delta (x_3x_2)+3\Delta
(x_2^3x_1)+3\Delta (x_3x_2x_1^2))+h^2\Delta (x_2f_4(x_1))$$
leads to two equations. The first one is
$$-f_4(x\otimes x)+(15-2C_2)x\otimes x+f_4(x)\otimes x^5-9x^3\otimes
x^3+(-15+2C_2)x\otimes x^5+9x^3\otimes x^5=0,\tag 8.23$$
and comes from a term proportional to $h^2x_2\otimes 1$.
The second equation comes from a term proportional to $h^2(1\otimes x_2)$ and
reads
$$-f_4(x\otimes x)+x\otimes f_4(x)+9x\otimes x^3-9x^3\otimes x^3-C_2x\otimes
x^5+C_2x^5\otimes x^5=0.\tag 8.24$$
Solving together (8.23) and (8.24) one obtains
$$f_4(x)=(15-2C_2)x-9x^3+C_2x^5.\tag 8.25$$
There are no terms of higher order in $h$ that do not cancel after the
reduction to a canonical form.
\par (e) We move on to the next relation which gives
$$\align \Delta (x_3x_5)=\Delta (x_5x_3)&+h\Delta
(5x_5x_1-9x_3^2+3x_3x_2^2x_1+3x_3^2x_1^2-x_5x_1^3)\\
&+h^2\Delta \Bigl[x_3f_5(x_1)+x_2^2f_6(x_1)\Bigr]+h^3\Delta
f_7(x_1).\endalign$$
Terms of order $h^2$ give rise to five functional equations which we now
describe.\newline
(i) A term proportional to $x_3\otimes 1$ gives rise to
$$-f_5(x\otimes x)+(6-3C_2)x\otimes x+f_5(x)\otimes x^5+6x^3\otimes
x^3-6x^3\otimes x^5+(-6+3C_2)x\otimes x^5=0,\tag 8.26$$
and a term prortional to $1\otimes x_3$ gives rise to the equation
$$-f_5(x\otimes x)+x\otimes f_5(x)-6x\otimes x^3+6x^3\otimes
x^3+(3-C_2)x\otimes x^5+(-3+C_2)x^5\otimes x^5=0.\tag 8.27$$
Solving (8.26) and (8.27) together we obtain for $f_5$,
$$f_5(x)=(6-3C_2)x+6x^3+(-3+C_2)x^5.\tag 8.28$$
(ii) Terms proportional to $x_2^2\otimes 1$ and $1\otimes x_2^2$ give rise to
another two functional equations:
$$\align -f_6(x\otimes x)+x^4\otimes f_6(x)&=0,\quad \text{and}\tag 8.29\\
-f_6(x\otimes x)+1\otimes f_6(x)&=0\tag 8.30\endalign$$
respectively. The only solution that satisfies both (8.28) and (8.29) is
$$f_6(x)=0.\tag 8.31$$
(iii) The last term from the terms of order $h^2$ is a term
proportional to $x_2\otimes x_2$ which gives rise to the following equation
for $f_5$ and $f_6$:
$$-2f_5(x\otimes x)+(12-6C_2)x\otimes x-2(x\otimes x)f_6(x\otimes
x)+12x^3\otimes x^3+(-6+2C_2)x^5\otimes x^5=0.\tag 8.32$$
After substituting the solutions (8.25) and (8.31) into (8.32) we obtain that
it is satisfied identically.
There is one term of order $h^3$ that arises which gives rise to
$$x^2\otimes f_7(x)+f_7(x)\otimes x^8-f_7(x\otimes x)=0.\tag 8.33$$
The most general solution of (8.33) is given by
$$f_7(x)=C_3(x^8-x^2),\tag 8.34$$
where $C_3\in k$ is an arbitrary constant. No terms of higher order in $h$
arise.
\par (f) The last relation to be analyzed is
$$\align \Delta (x_4x_5)=\Delta (x_5x_4)&+h\Delta
(10x_5x_2-12x_4x_3+3x_4x_2^2x_1+3x_4x_3x_1^2-3x_5x_2x_1^2)\\ &+
h^2\Delta\Bigl[x_4f_8(x_1)+x_3x_2f_9(x_1)+x_2^3f_{10}(x_1)\Bigr]+h^3\Delta
(x_2f_{11}(x_1)).\endalign$$
After reducing to a canonical form both sides of the above relation we obtain
ten terms of order $h^2$ that do not cancel and
two terms of order $h^3$. We analyze first the terms of order $h^2$.\newline
(i) Two terms proportional to $x_4\otimes 1$ and $1\otimes x_4$ give rise to
the following two equations
$$\align-f_8(x\otimes x)+(-6-4C_2)x\otimes x+f_8(x)\otimes x^5+9x^3\otimes
x^3+(6+4C_2)x\otimes x^5-9x^3\otimes x^5&=0,\tag 8.35\\
-f_8(x\otimes y)+x\otimes f_8(y)-9x\otimes y^3+9x^3\otimes y^3+(3-C_2)x\otimes
y^5+(-3+C_2)x^5\otimes y^5&=0.\tag 8.36\endalign$$
The most general solution of (8.35) and (8.36) is
$$f_8(x)=(-6-4C_2)x+9x^3+(-3+C_2)x^5.\tag 8.37$$
(ii) Terms proportional to $x_3x_2\otimes 1$ and $1\otimes x_3x_2$ give rise to
the equations
$$\align 6(1\otimes 1)-f_9(x\otimes x)-6(1\otimes x^4)+f_9(x)\otimes x^4=0,\tag
8.38\\
-f_9(x\otimes x)+1\otimes f_9(x)&=0.\tag 8.39\endalign$$
The solution of the system (8.38) and (8.39) is
$$f_9(x)=6.\tag 8.40$$
(iii) A term proportional to $x_3\otimes x_2$ gives rise to
$$-3f_8(x\otimes x)+(-12-12C_2)x\otimes x-(x\otimes x)f_9(x\otimes
x)+27x^3\otimes x^3+(-9+3C_2)x^5\otimes x^5=0.\tag 8.41$$
After substituting (8.37) and (8.40) into (8.41) it yields an
identity. Similarly the term proportional to $x_2\otimes x_3$ leads to
an identity.\newline
(iv) Two terms proportional to $x_2^3\otimes 1$ and $1\otimes x_2^3$ lead to
the equations
$$\align -f_{10}(x\otimes x)+f_{10}(x)\otimes x^3&=0,\tag 8.42\\
1\otimes f_{10}(x)-(x\otimes 1)f_{10}(x\otimes x)&=0.\tag 8.43\endalign$$
The only solution of (8.42) and (8.43) solved together is
$$f_{10}(x)=0.\tag 8.44$$
(v) The terms proportional to $x_2\otimes x_2^2$ and $x_2^2\otimes x_2$ give
rise to
$$-f_8(x\otimes x)+(6-4C_2)x\otimes x-2(x\otimes x)f_9(x\otimes x)-3(x^2\otimes
x^2)f_{10}(x\otimes x)+9x^3\otimes x^3+(-3+C_2)x^5\otimes x^5=0,\tag 8.45$$
and
$$12(1\otimes 1)-2f_9(x\otimes x)-3(x\otimes x)f_{10}(x\otimes x)=0\tag 8.46$$
which are identically satisfied. This becomes obvious after substituting
(8.37), (8.40) and (8.44) into (8.45) and (8.46).\newline
The two terms of order $h^3$ are proportional to $x_2\otimes 1$ and
$1\otimes x_2$ and give rise to
$$\align -f_{11}(x\otimes x)+(3-2C_2-2C_3)x\otimes x+f_{11}(x)\otimes
x+(-3+2C_2+2C_3)x\otimes x^7=0,\tag 8.47\\
-f_{11}(x\otimes x)+x\otimes f_{11}(x)-3C_3x\otimes x^7+3C_3x^7\otimes
x^7=0\tag 8.48\endalign$$
respectively. The most general solution of (8.47) and (8.48) is
$$f_{11}(x)=(3-2C_2-2C_3)x+3C_3x^7.\tag 8.49$$
\par We need one last step in order to complete the construction. We would like
to find whether the so obtained set of relations
define an algebra with the PBW property. After lengthy and tedious calculation
one shows that the requirement that the monomials
$x_1x_3x_5$, $x_1x_4x_5$, $x_2x_4x_5$, and $x_3x_4x_5$ can be reduced to a
unique canonical form imposes the following single
equation on the arbitrary constant $C_2$:
$$-9+2C_2=0.\tag 8.50$$
The monomial $x_1x_2x_5$ is reducible to a canonical form without imposing any
conditions. Thus $C_2=\frac {9}{2}$.
If we introduce $C:= C_3$ we obtain the statement of the Theorem. This
concludes the proof.\qed\enddemo
\remark{Remark} Notice that our construction yields a one-parameter family of
quantum semigroups $G^{2|C}_{5|h}$ parametrized by the
parameter $C$. The following theorem  will describe a family of quantum
semigroups parametrized by even more parameters.
\endremark
\par Let $G^1_4=G^1_{\infty}\text{mod}\{u^{n+1}\}$, for $n\ge 4$, be
the finite dimensional (dim\ =\ 4, $d=1$) Poisson-Lie group with a Poisson-Lie
 structure defined by
$$\eqalign{\{x_i,x_j\}=&(i-1)jx_jx_{i-1}-i(j-1)x_ix_{j-1}\cr
&\qquad+x_i\sum_{s_1+s_2=j}x_{s_1}x_{s_2}-x_j\sum_{s_1+s_2=i}x_{s_1}x_{s_2}.}
\tag 8.51$$
The above formulae are obtained from (8.2) with $d=1$. Writing them explicitly
we have
$$\align\{x_1,x_2 \}=&x_1^3-x_1^2\\
\{x_1,x_3 \}=&2x_2(x_1^2-x_1)\\\allowdisplaybreak
\{x_2,x_3\}=&(3x_1-x_1^2)x_3+x_2^2(2x_1-4)\\\allowdisplaybreak
\{x_1,x_4\}=&x_3(2x_1^2-3x_1)+x_2^2x_1\tag 8.52\\
\{x_2,x_4\}=&x_4(4x_1-x_1^2)+x_3x_2(2x_1-6)\\
\{x_3,x_4\}=&x_4x_2(8-2x_1)+x_3x_2^2+x_3^2(2x_1-9).\endalign$$

\proclaim{Theorem 8.4} Let $X=\{x_i\}_{1\le i\le 4}$ be a set of
indeterminates and let $\langle X\rangle$ be the associative semigroup with
identity generated by $X$.
Consider an ideal $\Cal I_{h}^1$ generated by the set of relations  $\Cal
R^1_h$ in $k[[h]]\langle X\rangle$
$$\align x_1x_2=x_2x_1+&h(x_1^3-x_1^2)\\
x_1x_3=x_3x_1+&h(2x_2x_1^2-2x_2x_1)\\
+&h^2(2x_1^4-3x_1^3+x_1^2)\\
x_2x_3=x_3x_2+&h(3x_3x_1-x_3x_1^2+2x_2^2x_1-4x_2^2)\\
+&h^2(3x_2x_1^2-3x_2x_1)\\
+&h^3(2-2C_3)(x_1^5-x_1^2)\\
x_1x_4=x_4x_1+&h(-3x_3x_1+2x_3x_1^2+x_2^2x_1)\\+&
h^2x_2(3x_1-8x_1^2+5x_1^3)\\+&
h^3\Bigl[(5x_1^5-12x_1^4+7x_1^3)+C_3(x_1^5-x_1^2)\Bigr]\tag
8.53\\\allowdisplaybreak
x_2x_4=x_4x_2+&h(4x_4x_1-x_4x_1^2+2x_3x_2x_1-6x_3x_2+x_2^3)\\+&
h^2(3x_2^2x_1^2-10x_2^2x_1+12x_2^2+12x_3x_1^2-2x_3x_1^3-15x_3x_1)\\+&
h^3x_2\Bigl[(9+2C_3)x_1-17x_1^2+6x_1^3+(5-5C_3)x_1^4\Bigr]\\+&
h^4\Bigl[(22-22C_3)x_1^2+(-4+4C_3)x_1^3+(-18+18C_3)x_1^5+C_4(x_1^6-x_1^2)\Bigr]\\\allowdisplaybreak
x_3x_4=x_4x_3+&h(8x_4x_2-2x_4x_2x_1+x_3x_2^2-9x_3^2+2x_3^2x_1)\\+&
h^2\Bigl[x_3x_2(-x_1^2+16x_1-24)+x_4(-7x_1^2+16x_1)\Bigr]\\+&
h^3\Bigl[(10-10C_3)x_2^2x_1^3+(-9+8C_3)x_2^2\Bigr]\\+&
h^3\Bigl[(-5+5C_3)x_3x_1^4+16x_3x_1^2-(6+9C_3)x_3x_1\Bigr]\\+&
h^4\Bigl[(8-9C_3-2C_4)x_2x_1+(-8C_3+9)x_2x_1^2+(10-10C_3)x_2x_1^4\Bigr]\\+&
h^4(2C_4-18+18C_3)x_2x_1^5\\+&
h^5\Bigl[(C_4+2C_3-2)x_1^2+(4-4C_3-C_4)x_1^6\Bigr]\\+&
h^5\Bigl[(2C_3-2)x_1^5+C_5(x_1^7-x_1^2)\Bigr].\endalign$$
where $C_3,C_4,C_5\in k$ are arbitrary parameters.  Then the semigroup factor
algebra $k[[h]]\langle X\rangle/\Cal I_{h}^1$ defines
 a quantum semigroup $G^1_{4|h}$ in
the sense of Definition {\rm 8.2} with a comultiplication defined by {\rm
(8.13)}.
Also the semigroup algebra $k[[h]]\langle X\rangle$ has the
Poincar\'e-Birkhoff-Witt property.\endproclaim
\remark{Remark} The proof of Theorem 8.4 goes along the same lines as the proof
of Theorem 8.3, i.e\. it is constructive. In the
course of the construction five arbitrary constants $C_1,C_2,C_3,C_4,C_5,$
appear in solving the corresponding functional equations.
The requirement for the existence of a PBW basis fixes two of them. Namely,
$C_1=1$ and $C_2=2-2C_3$. Thus we obtain a 3-parameter family of quantum
semigroups $G^{1|C_3,C_4,C_5}_{4|h}$.\endremark
\par Finally we describe a third quantum semigroup arising after the
quantization of the Poisson algebra of functions
on the finite dimensional (dim\ =\ 5, $d=3$) Poisson-Lie group $G^3_{5}$. The
Poisson-Lie structure on $G^3_{5}$ is given by
$$\eqalign{\{x_i,x_j\}=&(i-3)jx_jx_{i-3}-i(j-3)x_ix_{j-3}\cr &\qquad
+x_i\sum_{s_1+s_2+s_3+s_4=j}x_{s_1}x_{s_2}x_{s_3}x_{s_4}
-x_j\sum_{s_1+s_2+s_3+s_4=i}x_{s_1}x_{s_2}x_{s_3}x_{s_4}.}\tag 8.54$$
The above formulae are obtained again from (8.2) with $d=3$. Writing them
explicitly we have
$$\allowdisplaybreaks\align\{x_1,x_2 \}&=0\\
\{x_1,x_3 \}&=0\\
\{x_2,x_3\}&=0\\
\{x_1,x_4\}&=x_1^5-x_1^2\\
\{x_2,x_4\}&=x_2(x_1^4-2x_1)\\
\{x_3,x_4\}&=x_3(x_1^4-3x_1)\tag 8.55\\
\{x_1,x_5\}&=x_2(4x_1^4-2x_1)\\
\{x_2,x_5\}&=x_2^2(4x_1^4-4)\\
\{x_3,x_5\}&=x_3x_2(4x_1^3-6)\\
\{x_4,x_5\}&=x_4x_2(4x_1^3-8)+x_5(5x_1-x_1^4).\endalign$$
Then we have our last theorem.
\proclaim{Theorem 8.5} Let $X=\{x_i\}_{1\le i\le 5}$ be a set and let $\langle
X\rangle$ be the associative semigroup with identity
generated by $X$.
Consider an ideal $\Cal I_{h}^3$ generated by the set of relations  $\Cal
R^3_h$ in $k[[h]]\langle X\rangle$:
$$\align x_1x_2 =x_2x_1\\
x_1x_3 =x_3x_1&\\
x_2x_3=x_3x_2&\\\allowdisplaybreak
x_1x_4=x_4x_1&+h(x_1^5-x_1^2)\\\allowdisplaybreak
x_2x_4=x_4x_2&+hx_2(x_1^4-2x_1)\\\allowdisplaybreak
x_3x_4=x_4x_3&+hx_3(x_1^4-3x_1)\tag 8.56\\\allowdisplaybreak
x_1x_5=x_5x_1&+hx_2(4x_1^4-2x_1)\\\allowdisplaybreak
x_2x_5=x_5x_2&+hx_2^2(4x_1^4-4)\\\allowdisplaybreak
x_3x_5=x_5x_3&+hx_3x_2(4x_1^3-6)\\\allowdisplaybreak
x_4x_5=x_5x_4&+h\left[x_4x_2(4x_1^3-8)+x_5(5x_1-x_1^4)\right]+h^23x_2x_1.\endalign$$
Then the semigroup factor algebra $k[[h]]\langle X\rangle/\Cal I_{h}^3$ defines
a quantum
semigroup $G^3_{5|h}$ in
the sense of Definition {\rm 8.2} with a comultiplication defined by {\rm
(8.13)},
and the semigroup algebra $k[[h]]\langle X\rangle$ has the
Poincar\'e-Birkhoff-Witt property.\endproclaim
\remark{Remark} The proof is again constructive. Note that no arbitrary
parameters arise in dimension 5 for $d=3$. Arbitrary parameters
arise in higher dimensions though. We have been able to
construct all quantum semigroups $G^d_{n|h}$ for $n\le 7$ and $d\le
5$. However, we refrain from describing more
examples here.\endremark

\subhead Acknowledgements\endsubhead

\par The author would like to thank B. Kupershmidt for encouragement
and help, and I. Gelfand for helpful discussions.

\end